\documentclass[11pt]{article}
\usepackage{amssymb,amsthm,amsfonts,hhline,color}
\usepackage{longtable,tocloft}

\def\C{{\bf C}}

\def\Z{{\bf Z}}

\def\Q{{\bf Q}}
\def\tr{{\rm Tr}}
\def\K32{K3^{[2]}}
\def\M23{M_{23}}
\def\S{${\mathcal{S}}$}
\def\Aut{{\rm Aut}}
\def\rk{{\rm rk}\,}

\def\pf{\noindent{\bf Proof:\ }}
\def\qed{\hfill\framebox[2.5mm][t1]{\phantom{x}}}

\title{Finite groups of symplectic automorphisms \linebreak  of hyperk\"ahler manifolds of type~$\K32$}

\author{Gerald H\"ohn\\
Department of Mathematics, Kansas State University
\\  Geoffrey Mason\thanks{Supported by the NSF}\\
Department of Mathematics, University of California at Santa Cruz}

\date{April, 2016}

\begin{document}

\bibliographystyle{amsalpha}

\theoremstyle{plain}
\newtheorem{thm}{Theorem}[section]
\newtheorem{prop}[thm]{Proposition}
\newtheorem{lem}[thm]{Lemma}
\newtheorem{cor}[thm]{Corollary}
\newtheorem{rem}[thm]{Remark}
\newtheorem{conj}[thm]{Conjecture}

\newtheorem{introthm}{Theorem}
\renewcommand\theintrothm{\Alph{introthm}}

\theoremstyle{definition}
\newtheorem{defi}[thm]{Definition}

\renewcommand{\baselinestretch}{1.2}

\maketitle

\begin{abstract}
\noindent
We determine the possible finite groups $G$ of symplectic automorphisms of
\linebreak  hyperk\"ahler manifolds which are
deformation equivalent to the second Hilbert scheme of a K3 surface.\
We prove that $G$ has such an action if, and only if, it is isomorphic to a subgroup of either the
Mathieu group $\M23$ having at least four orbits in its natural permutation representation on
$24$ elements, or one of two groups $3^{1+4}{:}2.2^2$ and $3^4{:}A_6$ associated to {\S}-lattices
in the Leech lattice.\
We describe in detail those $G$ which are maximal with respect to these properties,
and (in most cases) we determine all deformation equivalence classes of such group actions.\
We also compare our results with the predictions of Mathieu Moonshine.
\end{abstract}

\tocloftpagestyle{plain}

{
\renewcommand{\baselinestretch}{0.8}
\tableofcontents
\small
\renewcommand{\baselinestretch}{0.99}
\listoftables
}

\section{Introduction}

A \emph{hyperk\"ahler manifold} is a $4n$-dimensional compact Riemannian manifold with holonomy group contained in 
${\rm Sp}(n)$.\ Such a manifold is  of \emph{type $\K32$} if it is deformation equivalent
to the second Hilbert scheme of a K3 surface.\ An example of a K3 surface is the Fermat quartic $Y\subset {\bf CP}^3$ given by
the equation $x_0^4+x_1^4+x_2^4+x_3^4=0$.\ An isometry of a hyperk\"ahler manifold fixing the
complex structures is called a {\it symplectic automorphism.\/}\ 
See~\cite{Huy-basic} for a review of basic properties of hyperk\"ahler manifolds.

\medskip

In the present paper, we determine and study those finite groups $G$ which can occur as groups of symplectic automorphisms 
of hyperk\"ahler manifolds of type $\K32$.\
Recent work of Mongardi \cite{Mon-thesis} shows that $G$ is isomorphic to a subgroup of
the Conway group ${\rm Co}_0$, the group of isometries of the Leech lattice $\Lambda$.\
Moreover, the fixed-point sublattice $\Lambda^G$ must have rank at least $4$.\ Mongardi also gave (loc.\ cit.)\ restrictions on the possible automorphisms of prime order.\ In a well-known paper \cite{Mu},
Mukai showed that a finite group of symplectic automorphisms of a K3 surface is isomorphic to a subgroup of the Mathieu  group $M_{23}$  having at least five orbits on its defining action on  $24$ elements.\
Our main result is the following theorem, which may be regarded as a higher-dimensional analog of Mukai's result.
\begin{introthm}\label{thmmain}
Let $G$ be a finite group of symplectic automorphisms of a hyperk\"ahler manifold of type $\K32$.\ 
Then $G$ is isomorphic to  one of the following:
\vspace{-2mm}
\begin{itemize}
\itemsep0em
\item[(a)] a subgroup of  $M_{23}$  with at least four orbits in its natural action on $24$ elements, 
\item[(b)] a subgroup of one of two subgroups 
$3^{1+4}{:}2.2^2$ and $3^4{:}A_6$ of ${\rm Co}_0$ associated to  \S-lattices in the Leech lattice. 
\end{itemize}
\end{introthm}
There are $13$ isomorphism classes of subgroups of type~(a) that are \emph{maximal\/} in the poset 
of all such groups.\  We will describe them and the two maximal groups of type~(b) in detail.\

\smallskip

By explicit construction, Mukai also showed (loc.\ cit.)\ that every subgroup of $M_{23}$ that satisfies the conditions 
of his theorem indeed occurs as a group of symplectic automorphisms 
of a $K3$ surface.\
These groups also act on the corresponding Hilbert schemes, thereby providing
examples of groups $G$ as in part (a) of Theorem \ref{thmmain}, and examples 
explicitly realizing several more of the maximal groups are known.\
We will establish the full analog of Mukai's result, namely:
\begin{introthm}\label{thmreal}
Each group $G$ in Theorem~\ref{thmmain}
can be realized as group of symplectic automorphisms of some hyperk\"ahler manifold of type $\K32$.
\end{introthm}

\smallskip 

Hashimoto has classified~\cite{Ha} all the deformation equivalence classes of finite symplectic
group actions on K3 surfaces.\ He found that for each group permitted by Mukai's theorem there is
\emph{a unique}  such class, with the exception of five cases where there are two such classes.\
We obtain a similar result for $\K32$:
\begin{introthm}\label{thmdeformationclasses}
There are at least $243$ deformation classes of finite symplectic group actions on hyperk\"ahler manifold of type $\K32$.
\end{introthm}
We can deduce Theorem \ref{thmdeformationclasses} from the following purely lattice-theoretic result:
\begin{introthm}\label{thmlatticeclasses}
Let $L$ be the unique even, integral lattice of signature $(3,20)$ and discriminant group of order $2$.\
There are at least $243$ conjugacy classes of subgroups $G$ of the isometry group $O(L)$ of $L$ such
that  the orthogonal complement $L_G$ of the fixed-point lattice $L^G$ in $L$
satisfies the following three properties:
\vspace{-3mm}
\begin{itemize}
\itemsep0em
\item[(i)]  $L_G$ is negative-definite;
\item[(ii)]  $L_G$ contains no vectors of norm $-2$;
\item[(iii)]  $L_G$ contains no vectors $v$ of norm $-10$ such that $v/2$ is contained in the dual lattice $L^*$.
\end{itemize}
\end{introthm}
The different classes can be read off from Tables~\ref{Gconclasses},~\ref{Lisoclasses} and~\ref{complementlattice1}. 

\medskip

The methods used in our paper are based on ideas developed by Nikulin, Mukai, Kond\={o} and Hashimoto
for K3 surfaces~\cite{Nikulin,Mu,Ko,Ha}.\ We also use fundamental results on the geometry of
hyperk\"ahler manifolds obtained by many authors in recent decades, including work on the global Torelli theorem
due to Huybrechts, Markman and Verbitsky.\ Recent results of Mongardi are crucial in allowing us to achieve a complete 
classification.

\smallskip

We also provide a somewhat new and more conceptual proof of Mukai's original result on symplectic automorphisms
of K3 surfaces \cite{Mu}.\ To explain this,  let $N$ denote the $K3$-lattice, 
i.e.,\ the unique even, unimodular, integral lattice of signature $(3,19)$, and let $G$ be a group of symplectic automorphisms
of a $K3$ surface.\ Kond\={o} showed in~\cite{Ko} that the lattice $N_G(-1)\oplus A_1$ can be embedded into one the $23$ Niemeier lattices
with roots, and a case-by-case analysis reveals that $G$ must be a subgroup of $M_{23}$ with at least five
orbits.\ Conversely,  in the appendix of~\cite{Ko}, and again by a case-by-case analysis, Mukai is able to realize 
each group arising from such a lattice construction as symplectic automorphisms.\ 
In~\cite{Ha}, Hashimoto computed all possible embeddings of $N_G(-1)\oplus A_1$ into the 
Niemeier lattices, and it turns out that the $82$ isomorphism types of group lattices $(N_G,G)$ are
in one-to-one correspondence with the combinatorial structure of symplectic group actions as determined previously in~\cite{xiao}.\
In our approach, we embed $N_G(-1)$ into the Leech lattice $\Lambda$.\ 
This essentially reduces the computation of the group lattices $(N_G,G)$
to the {\it group-theoretic problem\/} of enumerating all conjugacy classes of subgroups $G\subseteq {\rm Co}_0$ such that $rk(\Lambda^G)\geq 5$, and
indeed we find that there are just $82$ isomorphism types of such $(N_G,G)$.\ 
In addition, we clarify the result~\cite{Ha} that $(N_G,G)$ together with $N^G$
uniquely determine the conjugacy class of $G$ in $O(N)$.\ This is achieved by an improved {\it group theoretical
analysis\/} of embeddings $N_G\oplus N^G\subseteq N$.\ To a large extent, the above analysis of symplectic automorphisms
of K3 surfaces is contained in the corresponding analysis for $\K32$. We will mention the results for K3 surfaces, and possible
necessary modifications, at relevant points in the paper.

\medskip

Much of our interest in the subject matter of the present paper originates from issues surrounding moonshine.\
Hirzebruch suggested~\cite{Hi-book} that the Witten genus of a hypothetical $24$-dimensional monster manifold could be related 
to monstrous moonshine.\ Furthermore, the equivariant denominator identity of the monster Lie algebra can be interpreted
as the equivariant second quantized Witten genus of a monster manifold as noted by the first author~\cite{Ho-ober}.\
Mathieu Moonshine~\cite{EOT} connects the Mathieu group $M_{24}$ with  the complex elliptic genus of a K3 surface.
It seems natural to investigate geometric questions dealing with the equivariant second quantized complex elliptic genus of a K3 surface.\ 
See also~\cite{C} for a physical interpretation.\ 
Mathieu Moonshine is also closely related to a multiplicative version of Moonshine for $M_{24}$ found by the second author~\cite{MaM24}.\
Important input also came from recent work of Gaberdiel, Hohenegger and Volpato~\cite{GHV3}, 
where the lattice approach of Kond\=o for K3 surfaces was partially generalized to sigma models on K3 surfaces.

\medskip

The paper is organized as follows.\
In Section~\ref{SSK32} we cover required background about integral lattices and hyperk\"ahler manifolds of type $\K32$.\ 
In Section~\ref{slattices} we discuss the conjugacy classes of subgroups $G\subseteq {\rm Co}_0$ with $\rk(\Lambda^G)\geq 4$.\
Building on Mongardi's work, in Section \ref{SSCC} we show that there are exactly
$15$ conjugacy classes of elements in ${\rm Co}_0$ (the \emph{admissible} classes) that can occur as symplectic automorphisms of $\K32$.\
This is achieved by applying the equivariant Atiyah-Singer index formula to Hirzebruch's $\chi_y$-genus.\ 
In Section~\ref{SSGMT}, a group-theoretic analysis based on this conjugacy restriction then shows that the only groups 
satisfying  $\rk(\Lambda^G)\geq 4$ and consisting only of admissible elements  are those described in parts (a) and (b) of Theorem~\ref{thmmain}, 
or certain groups of order $12$, $16$, $32$, $48$ or $64$.\ Apart from a certain (inevitable) amount of computer calculation, 
the methods here are an extension of those used in \cite{Mason} to study the corresponding problem for K3 surfaces.\ 
In Section~\ref{admissibleconj},  we show that there are exactly $198$ conjugacy classes of such groups in ${\rm Co}_0$.\
In Section~\ref{subol} we determine --- apart from a few cases --- the conjugacy classes of groups in $O(L)$ 
that arise from these $198$ conjugacy classes, while in Section~\ref{geometricreal} we determine which of these 
conjugacy classes arise from symplectic group actions on some $\K32$.\
In the final section, we compare the equivariant complex elliptic genus of a $\K32$ with the predictions of Mathieu Moonshine 
applied to the second quantized elliptic genus.\ 
In the appendix, we describe the conjugacy classes of subgroups $G\subseteq {\rm Co}_0$ found in Section~\ref{admissibleconj},
together with additional information about  $G$ and the corresponding
lattices $L_G$. 

{\small
\paragraph{Acknowledgments.}
The first author thanks D.~Huybrechts for discussions which partially motivated this work,
and G.~Mongardi for answering questions about his work.\ 
We are grateful to T.~Creutzig, M.~Gaberdiel, V.~V.~Nikulin and R.~Volpato for useful discussions.\
The first author enjoyed the hospitality of the Hausdorff Research Institute of Mathematics in Bonn
during the early stages of this work, and we are both indebted to the Simons Center for inviting us to 
the workshop on Mock Modular Forms, Moonshine, and String Theory during the Fall 2013, where
part of this work was  done.\ The Simons Foundation also provided us with a license for the computer algebra system
Magma, and G.~Nebe and D.~Lorch helped us by providing certain Magma procedures. }



\section{Background on  $K3^{[2]}$}\label{SSK32}

\subsection{Integral lattices}

We introduce some notation related to integral lattices and record some results that we will need. 

\smallskip

Let $A=(A,q)$ be a {\it finite quadratic space}, i.e.,\ a finite abelian group $A$ together with a 
quadratic form $q: A\longrightarrow \Q/2\Z$. We denote the corresponding orthogonal group by $O(A)$.\ This is the subgroup
of ${\rm Aut}(A)$ that leaves $q$ invariant.

Let $L$ be an even integral lattice, with \emph{dual lattice} $L^*$. The
 {\it discriminant group} $L^*/L$ is equipped with the
{\it discriminant form} $q_L: L^*/L \rightarrow \Q/2\Z$, $x+L \mapsto \langle x,\,x\rangle\ {\rm mod}\ 2\Z$.\ This
turns $L^*/L$ into a finite quadratic space, called the {\it discriminant space\/} of $L$ and denoted by
$A_L:=(L^*/L, q_L)$.\ We let $O(L):={\rm Aut}(L)$ be the automorphism group (i.e.\ group of isometries) of $L$.

Automorphisms in $O(L)$ induce  orthogonal transformations of the discriminant space $A_L$. 
This leads to the short exact sequence
$$1\longrightarrow O_0(L)   \longrightarrow O(L) \longrightarrow \overline{O}(L)  \longrightarrow 1,$$
where $\overline{O}(L)$ is the subgroup of $O(A_L)$ induced by $O(L)$ and $O_0(L)$ consists of the automorphisms
of $L$ which act trivially on $A_L$.

\smallskip
Let $\mathcal{L}$ be the category whose objects are even, integral lattices and whose morphisms are
\emph{injective isometries}.\ The category of \emph{group lattices} consists of objects $(L, G)$ where
$L$ is an object of $\mathcal{L}$ and $G \subseteq {\rm Aut}_{\mathcal{L}}(L)=O(L)$ is a subgroup.\ A morphism
$(L, G)\rightarrow(L', G')$ of group lattices is a pair $(\iota, j)$ where $\iota : L\rightarrow L'$ is a morphism in
$\mathcal{L}$ and $j : G \rightarrow G'$ is an \emph{injective\/} morphism of groups such that the following diagram 
commutes for all $g \in G$:
$$\begin{array}{ccc}
 L &  \stackrel{g}{\longrightarrow} &  L\\
\ \ \downarrow\iota  & & \ \ \downarrow\iota\\
  L' & \stackrel{j(g)}{\longrightarrow}  &  L'
\end{array}_.$$
In particular, if $(L,G)$ is a group lattice and $\iota : L \rightarrow L'$  an
\emph{isomorphism} in $\mathcal{L}$ we set 
\begin{eqnarray*}
\iota[G]:=\{ \iota \circ g \circ \iota^{-1}\mid g \in G\}.
\end{eqnarray*}
Then $(L', \iota[G])$ is a group lattice isomorphic to $(L, G)$.
Upon identifying $L$ and $L'$, this just means that $G$ and $G'$ are conjugate subgroups of $O(L)$.

The \emph{invariant} and \emph{coinvariant} lattices of a group lattice $(L, G)$ are respectively defined as follows:
\begin{eqnarray*}
L^G &=& \{x\in L\mid gx=x \hbox{ for all\ } g\in G\}, \\
L_G &=& \{x\in L\mid (x,y)=0 \hbox{ for all\ } y\in L^G\}.
\end{eqnarray*}
These are both {\it primitive sublattices\/} of $L$, i.e.,
$L/L^G$ and $L/L_G$ are free abelian groups.\ The restriction of the $G$-action to $L_G$ turns it
into a group lattice $(L_G,G)$. 
Moreover, $G$ acts trivially on the discriminant group $A_{L_G}$.

As a matter of notation, by $L(n)$ we will mean a lattice $L$ with norms scaled by an integer $n$.\ 
We also note that the \emph{genus\/} of an even integral lattice $L$ is determined by
the quadratic space $A_L$ together with the signature of $L$ \cite{Nikulin}.


\subsection{Automorphisms of $\K32$} 

In this subsection we fix the following notation:\
$X$ is a hyperk\"ahler manifold of type $\K32$,  ${\rm Aut}(X)$ is the group
of \emph{symplectic} automorphisms of $X$, and 
$G\subseteq {\rm Aut}(X)$  is a \emph{finite} subgroup.

\medskip
The second integral cohomology $L:=H^2(X,\Z)$ admits a non-degenerate symmetric integral bilinear 
form $(\,.\,,\,.\,)$, called the \emph{Beauville-Bogomolov} form,  with respect to which $L$ is isomorphic to 
the lattice $E_8(-1)^2\oplus U^3\oplus \langle -2 \rangle$
of signature $(3,20)$.  
Here, $E_8(-1)$ denotes the unique even, unimodular, negative-definite lattice of rank~$8$, 
$U$ the hyperbolic plane and  $\langle -2 \rangle\cong A_1(-1)$ the $1$-dimensional lattice generated by a vector of norm $-2$.\ 
The discriminant space $(A_L,q_L)\cong (A_{A_1},-q_{A_1})$ has order order~$2$.\
The group $\overline{O}(L)$ is trivial, i.e., $O(L)$ acts trivially on $A_L$.

\smallskip

There is an {\it injective\/} map (\cite{Mon-thesis,HaTsch,Beau}) 
$$\nu:{\rm Aut}(X)\longrightarrow O(L)$$
by which we may, and shall, identify $G$ with its image in $O(L)$.

\medskip

\begin{thm}[Mongardi~\cite{Mon-K32inv}, Lemma~3.5]
The coinvariant lattice $L_G$ has the following properties:
\vspace{-3mm}
\begin{itemize}
\itemsep0em
\item[(i)] $L_G$ is negative definite.
\item[(ii)]  $L_G$ contains no vectors of norm $-2$.
\end{itemize}
\end{thm}
It is also shown that $L_G$ is contained in the Picard lattice of $X$.

\smallskip

Recall that the \emph{Leech lattice} $\Lambda$  is the unique positive-definite, even, unimodular lattice of rank~$24$ without roots.
Its automorphism group is the Conway group ${\rm Co}_0$.
\begin{thm}\label{Coembedd}
There is an embedding of group lattices $(L_G(-1),G)\rightarrow (\Lambda,{\rm Co}_0)$ such
that $(L_G(-1),G)\cong (\Lambda_G,G)$.
\end{thm}
\pf
The discriminant form $q_L$ of $L$ is the negative of the discriminant form $q_{A_1}$ of the
root lattice $A_1=\langle 2 \rangle$.\ This permits us to extend the lattice $L\oplus A_1$ 
by the coset $(x,y)\in A_L\oplus A_{A_1}, (x\not= 0, y\not= 0)$
to the unique even unimodular lattice $M$ of signature $(4,20)$, thus providing a primitive embedding of $L$
into $M$.\
Since $\overline{O}(L)$ is trivial, the $G$-action extends to $M$, thereby fixing
the sublattice $L^G\oplus A_1\subseteq M$.\ Because $L_G$ is negative-definite,
 $L^G\oplus A_1\subseteq M$ has signature $(4,20-{\rm rk}\, L_G)$.\
In particular, we can find a $4$-dimensional positive-definite subspace $\Pi\subset M\otimes {\bf R}$
such that $L_G= \Pi^\perp \cap L$.\
It was shown in~\cite{GHV3} (see also~\cite{Huy-conway}, Prop. 2.2) that because $L_G$ is negative-definite,
$G$ can be embedded into the Conway group (a result first shown in~\cite{Mon-thesis}).\
Indeed, the proof actually shows that this corresponds to an embedding 
$(L_G(-1),G)\rightarrow (\Lambda,{\rm Co}_0)$ of group lattices with $(L_G(-1),G)\cong (\Lambda_G,G)$. $\hfill \Box$

Note that because $\rk L_G\leq 20$ then $\rk\Lambda^G\geq 4$.\ 
Theorem~\ref{Coembedd} allows us to identify $G$ with a subgroup of ${\rm Co}_0$.\
We will see in Section~\ref{admissibleconj} that for a given group lattice
$(L_G,G)$, the resulting embedding of $G \rightarrow {\rm Co}_0$ is \emph{unique} up to conjugation in ${\rm Co}_0$.


\section{{\S}-lattices and subgroups of ${\rm Co}_0$}\label{slattices}

We have seen in the previous section that a finite group $G$ of symplectic automorphisms of a  
hyperk\"ahler manifold $X$ of type $K3^{[2]}$ defines a subgroup $G\subseteq {\rm Co}_0$
with the property that $\rk(\Lambda^G)\geq 4$. 
In this section, we will establish some general results about such $G$.

\medskip

Recall~\cite{Con} that the $2^{24}$ cosets comprising $\Lambda/2\Lambda$ have representatives
$v$ which may be chosen to be \emph{short vectors}, i.e., 
$(v, v)\leq 8$.\ More precisely, if $(v, v)\leq 6$ then $\{v,-v\}$ are the \emph{only\/}
short representatives of $v+2\Lambda$; if $(v, v)=8$ then the short
vectors in $v+2\Lambda$ comprise a \emph{coordinate frame\/} $\{\pm w_1,\, \ldots,\, \pm w_{24}\}$,
where the $w_j$  are pairwise orthogonal vectors of norm~$8$.\ In particular, if $u\in\Lambda$ then
$u=v+2w$ for some $v$, $w\in\Lambda$ and $v$ a short vector, and if $(v, v)\leq 6$ then $v$ is 
\emph{unique\/} up to sign.

It is well-known that ${\rm Co}_0$ acts \emph{transitively} on coordinate frames,  the (setwise)
stabilizer of one such being the \emph{monomial group} $2^{12}{:}M_{24}$.

A sublattice $S\subseteq \Lambda$ is an \emph{\S-lattice} if, for every $u\in S$, the corresponding
short vector $v$ satisfies  $(v, v)\leq 6$ and furthermore $w\in S$.\ This concept, introduced by Curtis
\cite{Curtis}, will be very useful.

\smallskip

The following theorem depends in an essential way on a result of Allcock \cite{All}.\ See also \cite{GHV3}.
\begin{thm}\label{Cosubgps} Suppose that $G\subseteq {\rm Co}_0$ is a subgroup with ${\rm rk}\, \Lambda^G\geq 4$.
One of the following holds.
\vspace{-3mm}
\begin{itemize}
\itemsep0em
\item[(a)] $G$ leaves a coordinate frame invariant,
\item[(b)] $\Lambda^G$ is an {\S}-lattice of rank $4$,
\item[(c)] $\Lambda^G$ is contained in a $G$-invariant {\S}-lattice of rank larger than $4$.
\end{itemize}
\end{thm}

\pf Let $L:=\Lambda^G$.\ We may assume that (a) does not hold.\
Suppose that $u\in L$ with $u=v+2w$ and $(v, v)\leq 8$.\ Then $v+2\Lambda = u+2\Lambda$ is $G$-invariant, and 
since $G$ leaves no coordinate frame invariant then we have $(v, v)\leq 6$.\ Thus $G$ acts on $\{\pm v\}$.\ 
We claim that $v\in L\cup L^{\perp}$.\ For if $v\notin L$ there is a $g\in G$ such that
$g(v)=-v$.\ Then for $x\in L$ we obtain $(x, v)=(g(x), g(v))=(x, -v)$, showing that $v\in L^{\perp}$.

We use results of Allcock \cite{All}, especially (a special case of)  Lemma~4.8 (loc.\ cit.)\ which we state as follows:\
suppose that $L\subseteq \Lambda$ is a primitive sublattice
of rank at least $4$ with the property that if $u\in L$ with $u=v+2w$ and $(v, v)\leq 8$,
then $(v, v)\leq 6$ and  $v\in L\cup L^{\perp}$.\ 
Then $L$ is \emph{contained} in an {\S}-lattice.\ The previous paragraph
establishes that $L=\Lambda^G$ satisfies these properties, so $L$ is contained in an {\S}-lattice.\ 
Because the family of {\S}-lattices containing
$L$ is closed under intersection and $G$-conjugation, there is a \emph{$G$-invariant} {\S}-lattice
that contains $L$.\ Let $S$ be such a lattice.

If ${\rm rk}\, S =4$ then $L\subseteq S$ has finite index because of our assumption that ${\rm rk}\, L \geq 4$.\ 
Then $L=S$ because $L$ is primitive, and we are in case (b) of the Theorem.\ 
Otherwise ${\rm rk}\, S \geq 5$ and (c) holds.\
This completes the proof of the Theorem. \qed 

\medskip

We now draw some more detailed conclusions concerning the subgroups  $G\subseteq {\rm Co}_0$ 
using Theorem~\ref{Cosubgps}.\ 
This depends  on Curtis's classification of {\S}-lattices~\cite{Curtis}.\ 
See also \cite{GV} for a related discussion.
\begin{thm}\label{thmsgp} There are exactly six conjugacy classes of subgroups
$G\subseteq {\rm Co}_0$ such that $\rk\, \Lambda^G \geq 4$,
$G$ is the full (pointwise) stabilizer of $ \Lambda^G $ in ${\rm Co}_0$, and $G$ 
fixes \emph{no\/} coordinate frame.\
The group $G$ and $\Lambda^G$ are described as follows.
\renewcommand{\labelenumi}{(\roman{enumi})}
\vspace{-3mm}
\begin{enumerate}
\itemsep0em
\item $\Lambda^G$ is an {\S}-lattice of rank $6$ and $G\cong 3^{1+4}.2$,
\item $\rk\,\Lambda^G=5$ and $G\cong 3^{1+4}.2.2$,
\item $\Lambda^G$ is an {\S}-lattice of rank $4$ and $G\cong 3^4.A_6$,
\item $\Lambda^G$ is an {\S}-lattice of rank $4$ and $G\cong 5^{1+2}.4$,
\item $\rk\,\Lambda^G=4$ and $G\cong 3^{1+4}.2.2^2$,
\item $\rk\, \Lambda^G=4$ and $G\cong 3^{1+4}.2.2$.
\end{enumerate}
\end{thm}
\pf Set $L:=\Lambda^G$.\  By Theorem \ref{Cosubgps}, 
there is a $G$-invariant {\S}-lattice $S$ with $L\subseteq S$.\ Let $N$ be the (pointwise) stabilizer of $S$.\ 
Because $S$ is $G$-invariant then $G$ normalizes $N$, so since $N$ fixes $L$ pointwise then $N\unlhd G$.\ 
Set $\widetilde{G}:= G/N \subseteq {\rm Aut}_{{\rm Co}_0}(S)$.

The possibilities for $S$ are as follows (\cite{Curtis}, \cite{Atlas}):
$$\begin{array}{cccc}
\rk\, S & N & {\rm Aut}_{{\rm Co}_0}(S) \\ \hline
4 &   3^4.A_6 & 2 \times (S_3\times S_3).2\\
4 &  5^{1+2}.4 &2\times S_5 \\
6 &  3^{1+4}.2 & 2 \times U_4(2).2 
\end{array}$$

If $S/L$ is \emph{finite} then $L=S$ because $L$ is primitive, so we have $G=N$ and cases (i), (iii) or (iv) apply.\
Thus from now on we will assume that $S/L$ is \emph{not\/} finite.\ In particular we have $\rk\, S \geq 5$, whence
$S$ is the {\S}-lattice of rank~$6$ as we see from the table, moreover $\rk\, L=4$ or $5$.\ 
We also have $3^{1+4}\cong O^2(N)\unlhd G$ and $|\widetilde{G}|>1$.

\smallskip

From the table, the group of isometries of $S$ induced within ${\rm Co}_0$ is the group
$\{\pm 1\}\times W(E_6)$ (which is actually the \emph{full group of isometries\/}).\ 
Indeed, a generator of the direct factor $\pm 1$ acts on $S$ as $-1$, and $W(E_6)\cong U_4(2).2$
is the Weyl group of type $E_6$.\ The lattice $\frac{1}{\sqrt{3}}S$ is isometric to the weight lattice of type $E_6$, 
the roots corresponding to short vectors of norm~$6$ in $S$.

\smallskip

If $\rk\, L=5$ then $\widetilde{G}$ fixes a hyperplane pointwise, so that $\tilde{G}\cong \Z_2$ is generated 
by a reflection in a hyperplane of $S$ orthogonal
to a norm $6$ vector, and any two such hyperplanes are conjugate in the Weyl group.\ 
This is case (ii).

\smallskip

The case $\rk\, L=4$ requires more care. If $L$ is an {\S}-lattice then from the table,
it must be that $L$ is of the first kind, i.e., with stabilizer $3^4.A_6$.\
Although there is a containment of {\S}-lattices of this kind (\cite{Curtis} or
Tables \ref{Gconclasses} and \ref{Lisoclasses} below),  in our set-up
we have $3^{1+4}\unlhd G$, whereas $3^4.A_6$ has no such normal subgroup.\ 
Thus $L$ is \emph{not} an {\S}-lattice, and 
in the proof of Theorem \ref{Cosubgps}
we showed that in this situation we can find a nonzero short
vector $v\in L^{\perp}$.\ Then in the orthogonal $GF(2)$-space $\bar{\Lambda}:= \Lambda/2\Lambda$,
$v$ maps onto a nonzero element $\bar{v}\in rad(\overline{L})$, so that
$\overline{L}$ is a 
$4$-dimensional \emph{degenerate\/} subspace of the nondegenerate $6$-dimensional orthogonal space 
$\overline{S}$.\ ($\overline{L}$, $\overline{S}$ are the images of $L$, $S$ respectively in 
$\overline{\Lambda}$.)\
The pointwise stabilizer of such a degenerate subspace in the full isometry group $O^-_6(2)$ of 
$\overline{S}$ is a $2$-group, and as a result
it follows that $\widetilde{G}$ is also a $2$-group.\ Since no element of order $4$ in 
$\Aut_{{\rm Co}_0}(S)$ fixes a rank $4$ sublattice pointwise,
then $\widetilde{G} \cong \Z_2^k$ for some $k\geq 1$.

Suppose that $\widetilde{G}\not\subseteq W(E_6)$.\ There is a unique conjugacy class of involutions $t$ (the product
of the $-1$-involution and an involution in $U_4(2)$ of type $2A$)
in $\{\pm 1\}\times W(E_6)\setminus{W(E_6)}$ fixing a sublattice in $S$ of rank $\geq 4$ (\cite{Atlas}) 
 and the rank is exactly $4$.\ Moreover $2A$ is a central involution in $U_4(2).2$ and its centralizer $\pmod{\langle 2A \rangle}$ acts
 faithfully on its $-1$-eigenspace, which corresponds to the fixed-space of $-2A$.\ It follows that $\tilde{G}$ has order $2$ in this case, which is
 case (vi) of the Theorem.

The remaining possibility is $\widetilde{G}\subseteq W(E_6)$.\ Involutions in the Weyl group fixing a sublattice
of $S$ of rank $\geq 4$ pointwise are those of type $2B$ and $2C$ (\cite{Atlas}), the fixed-point ranks being $4$ and $5$ respectively.\
Moreover, the product of a pair of distinct commuting Weyl reflections (type $2C$) is of type $2B$.\ It follows that $L$ is the
sublattice of $S$ fixed pointwise by an involution of type $2B$ 
(so that all such sublattices are conjugate in $\Aut_{{\rm Co}_0}(S)$),
or equivalently by a pair of commuting Weyl reflections (so that $|\widetilde{G}|\geq 4$).\
Moreover, any subgroup of $\Aut_{{\rm Co}_0}(S)$ strictly containing $\widetilde{G}$ has fixed sublattice 
of rank \emph{no greater\/} than $3$, whence $\widetilde{G}\cong \Z_2^2$.\ 
This is case~(v).

From what we have proved so far, it follows that the pointwise stabilizers of rank $4$ sublattices of $S$ as in cases (v) and (vi)
are \emph{not} conjugate in the (setwise) stabilizer of $S$.\ On the other hand, if they are conjugate
by an element  $g\in {\rm Co}_0$, then $g$ must normalize their common normal subgroup $3^{1+4}=O^2(N)$.\ Since
the normalizer  of this group \emph{is} the setwise stabilizer of $S$, then $g$ must belong to this group, a contradiction.\
It follows that the groups corresponding to cases (v) and (vi), and then even to all six cases, are \emph{not} conjugate in
${\rm Co}_0$, and the proof of the theorem is complete. $\hfill \Box$


\section{Geometric conditions}\label{SSCC}

In this section, we use geometric arguments to obtain strong restrictions on the group-theoretic 
properties of finite groups of symplectic automorphism of a hyperk\"{a}hler manifold of 
type $K3^{[2]}$.

We determine which conjugacy classes of ${\rm Co}_0$
can arise as symplectic automorphisms.\ 
We refer to these conjugacy classes, and the elements in them, as 
the \emph{admissible conjugacy classes} and \emph{admissible elements} respectively.\ 
The remaining conjugacy classes and elements are called \emph{inadmissible}.\ 
The main result (Theorem~\ref{allowedclasses}) asserts that there are just $15$ admissible conjugacy classes.\
We also determine the structure of the fixed-point set of the admissible elements and the action
on the normal bundle.\ Finally, we show that a $2$-group of symplectic automorphisms has order at most $2^7$.

\medskip

As explained in Section~\ref{SSK32}, a finite group $G$ of symplectic automorphisms of a hyperk\"{a}hler manifold of type $K3^{[2]}$ 
can be identified with a subgroup of ${\rm Co}_0$.\ After making this identification, the primitive embedding of 
$L_G(-1)$ into the Leech lattice $\Lambda$ is such that 
$\rk(\Lambda^G)\geq 4$.\ Obviously then, we have $\rk(\Lambda^g)\geq 4$ for every
$g\in G$.

\medskip
We start with Table~\ref{Leechclasses}, which lists the $42$ conjugacy classes $[g]$ of ${\rm Co}_0$ that satisfy 
the condition $\rk(\Lambda^g)\geq 4$, together with some supplementary data.\ It transpires that such
a $[g]$ is uniquely specified by the triple $(\hbox{order of}\ g,\,{\rm Trace} (g),\,{\rm Trace} (g^2))$,
and this is the entry in the first column of the table.\ In what follows, we often identify a conjugacy class using this
triple.\ The second column is the \emph{Frame shape}\footnote{$g$ has \emph{Frame shape} 
$1^{m_1}2^{m_2}\ldots $ if its characteristic polynomial (considered as a linear transformation of $\Lambda \otimes \mathbf{R}$) 
is $(t-1)^{m_1}(t^2-1)^{m_2}\ldots$.\ Each $m_i\in\mathbb{Z}$.} of $g$, the third column gives $\rk(\Lambda^g)$, the fourth column the
\emph{torsion-invariants} of $A_{\Lambda^g}=(\Lambda^g)^*/\Lambda^g$, and the fifth column (`powers') the nontrivial
prime powers of $g$.\ Column six records whether $g$ belongs (up to conjugacy) to the monomial subgroup $2^{12}{:}M_{24}\subseteq {\rm Co}_0$
($*$ indicates that it \emph{does}).\ Finally, in the seventh column (`excluded') the symbol 4.x refers to the Lemma or Theorem 4.x below
by which the inadmissible elements are excluded.

\medskip

Columns $1$, $2$, $3$, $5$ and $6$ in Table \ref{Leechclasses} can be read-off from the Atlas~\cite{Atlas}. 
The structure of the fixed-point lattice has been investigated in~\cite{KT1,KT2,La,HaLa}.
The table was confirmed using Magma~\cite{magma} together with a realization of ${\rm Co}_0$ as a matrix group.

\begin{table}\caption{Conjugacy classes of ${\rm Co}_0$ with at least four-dimensional fixed-point lattice}
\label{Leechclasses}
{\small
$$\begin{array}{llrclcl}
\hbox{class\ }[g] & \hbox{Frame shape} & \rk\, \Lambda^g & A_{\Lambda^g} & \hbox{powers} & 2^{12}{:}M_{24} &  \hbox{excluded}\\ \hline
( 1, 24, 24 ) 
 &  1^{24}
 &  24
 &  1
 &     
 & *
 & -
\\ \hline

( 2, 8, 24 ) 
 &  1^82^8
 &  16
 &  2^8
 &    
 & *
 & -
\\

( 2, 0, 24 ) 
 & 2^{12}
 &  12
 &  2^{12}
 &     
 & *
 & \textcolor{red}{\ref{invtype}}
\\ 

( 2, -8, 24 ) 
 & 2^{16}/1^8
 &  8
 &  2^8
 &    
 & *
 &  \textcolor{red}{\ref{elim1}}, \, \textcolor{red}{\ref{invtype}}
\\ \hline

( 3, 6, 6 ) 
 &  1^63^6
 &  12
 &  3^6
 &     
 & *
 & -
\\

( 3, 0, 0 ) 
 &  3^8
 &  8
 &  3^8
 &     
 & *
 &  \textcolor{red}{\ref{elim1}}
\\

( 3, -3, -3 ) 
 &  3^9/1^3
 &  6
 &  3^5
 &     
 & {\rm No}
 & -
\\ \hline

( 4, 8, -8 ) 
 & 1^84^8/2^8
 &  8
 &  2^8
 &     ( 2, -8, 24 ) 
 & *
 & \textcolor{red}{\ref{argpower}},\,  \textcolor{red}{\ref{elim1}}
\\

( 4, 4, 8 ) 
 &1^42^24^4
 &  10
 &  2^2 4^4
 &    
    ( 2, 8, 24 ) 
 & *
 & -
\\

( 4, 0, 8 ) 
 &2^44^4
 &  8
 &  2^4 4^4
 &    
    ( 2, 8, 24 ) 
 & *
 & \textcolor{red}{\ref{argfix4}}
\\

( 4, 0, -8 ) 
 &4^8/2^4
 &  4
 &  2^2 4^2
 &     ( 2, -8, 24 ) 
 & *
 & \textcolor{red}{\ref{argpower}}
\\

( 4, 0, 0 ) 
 &4^6
 &  6
 &  4^6
 &    
    ( 2, 0, 24 ) 
 & *
 &  \textcolor{red}{\ref{argpower}}, \, \textcolor{red}{\ref{elim1}}
\\ 

( 4, -4, 8 ) 
 &2^64^4/1^4
 &  6
 &  2^2 4^4
 &    
    ( 2, 8, 24 ) 
 & *
 & \textcolor{red}{\ref{argfix4}}
\\ \hline

( 5, 4, 4 ) 
 & 1^45^4
 &  8
 &  5^4
 &     
 & *
 & -
\\

( 5, -1, -1 ) 
 & 5^5/1
 &  4
 &  5^3
 &     
 & {\rm No}
 &  \textcolor{red}{\ref{argfix5}}
\\ \hline

( 6, 5, -3 ) 
 &1^5.3.6^4/2^4
 &  6
 &  3^5
 &    
    ( 2, 8, 24 ) ,\,
    ( 3, -3, -3 ) 
 & {\rm No}
 & -
\\

( 6, 4, 6 ) 
 & 1^4.2.6^5/3^4
 &  6
 &  2^5 6^1 
 &  ( 2, -8, 24 ) ,\,
    ( 3, 6, 6 ) 
 & *
 & \textcolor{red}{\ref{argpower}}, \, \textcolor{red}{\ref{elim1}}
\\

( 6, 2, 6 ) 
 &1^22^23^26^2
 &  8
 &  6^4
 &  ( 2, 8, 24 )  ,\,  
    ( 3, 6, 6 )
 & *
 & -
\\

( 6, 0, 6 ) 
 & 2^36^3
 &  6
 &  2^3 6^3
 &    
 ( 2, 0, 24 ) ,\,
 ( 3, 6, 6 ) 
 & *
 & \textcolor{red}{\ref{argpower}}
\\

( 6, 0, 0 ) 
 & 6^4
 &  4
 &  6^4
 &  ( 2, 0, 24 ) ,\, 
    ( 3, 0, 0 )   
 & *
 &   \textcolor{red}{\ref{argpower}}, \, \textcolor{red}{\ref{elim1}}
\\ 

( 6, -2, 6 ) 
 &2^46^4/1^23^2
 &  4
 &  2^2 6^2
 &  ( 2,-8,24 ), \,
    ( 3, 6, 6 ) 
 & *
 &  \textcolor{red}{\ref{argpower}}, \, \textcolor{red}{\ref{elim1}}
\\

( 6, -1, -3 ) 
 &3^36^3/1.2
 &  4
 &  3^2 6^2
 &      ( 2, 8, 24 )  ,\,
        ( 3, -3, -3 )
 & {\rm No}
 & \textcolor{red}{\ref{elim1}}
\\

( 6, -4, 6 ) 
 &2^53^4.6/1^4
 &  6
 &  2^1  6^5
 &  ( 2, 8, 24 )  ,\,
    ( 3, 6, 6 )     
 & *
 &  \textcolor{red}{\ref{elim1}}
\\ \hline

( 7, 3, 3 ) 
 & 1^37^3
 &  6
 &  7^3
 &     
 & *
 & -
\\ \hline

( 8, 4, 0 ) 
 &1^48^4/2^24^2
 &  4
 &  2^2 4^2
 &    
    ( 4, 0, -8 ) 
 & *
 & \textcolor{red}{\ref{argpower}}
\\

( 8, 0, 8 ) 
 &2^48^4/4^4
 &  4
 &  4^4
 &    
    ( 4, 8, -8 ) 
 & *
 &  \textcolor{red}{\ref{argpower}}, \, \textcolor{red}{\ref{elim1}}
\\

( 8, 2, 4 ) 
 &1^2.2.4.8^2
 &  6
 &  2^1 4^1 8^2
 &     ( 4, 4, 8 )    
 & *
 & -
\\

( 8, 0, 0 ) 
 & 4^2 8^2
 &  4
 &  4^2 8^2
 &   
    ( 4, 0, 8 ) 
 & *
 & \textcolor{red}{\ref{argpower}}, \,  \textcolor{red}{\ref{elim1}}
\\ 

( 8, -2, 4 ) 
 &2^3.4.8^2/1^2
 &  4
 &  2^1 4^1 8^2
 &   
    ( 4, 4, 8 ) 
 & *
 & \textcolor{red}{\ref{argfix8}}
\\ \hline

( 9, 3, 3 ) 
 & 1^39^3/3^2
 &  4
 &  3^2 9^1

 &     ( 3, -3, -3 ) 
 & {\rm No}
 & - 
\\ \hline

( 10, 3, -1 ) 
 & 1^3.5.10^2/2^2
 &  4
 &  5^3
 &  ( 2, 8, 24 )  ,\, 
    ( 5, -1, -1 )   
 & {\rm No}
 & \textcolor{red}{\ref{argpower}}
\\

( 10, 2, 4 )
 &1^2.2.10^3/5^2
 &  4
 &  2^3 10^1 
 &    
    ( 2, -8, 24 ) ,\,
    ( 5, 4, 4 ) ,
 & *
 &  \textcolor{red}{\ref{argpower}}, \,  \textcolor{red}{\ref{elim1}}
\\

( 10, 0, 4 ) 
 &2^210^2
 &  4
 &  2^2 10^2
 &  ( 2, 0, 24 )  ,\,
    ( 5, 4, 4 )
 & *
 & \textcolor{red}{\ref{argpower}}
\\ 

( 10, -2, 4 ) 
 & 2^35^210/1^2
 &  4
 &  2^1 10^3
 &  ( 2, 8, 24 ) ,\, 
    ( 5, 4, 4 )
 & *
 &  \textcolor{red}{\ref{elim1}}
\\ \hline

( 11, 2, 2 ) 
 & 1^211^2
 &  4
 &  11^2
 &   
 & *
 & -
\\ \hline

( 12, 2, 2 ) 
 & 1^2.4.6^212/3^2
 &  4
 &  2^2 4^1 12^1 
 &     ( 4, -4, 8 ) ,\,
    ( 6, 2, 6 ) 
 & *
 & \textcolor{red}{\ref{argpower}}
\\

( 12, 2, -2 ) 
 & 1^23^24^212^2/2^26^2
 &  4
 &  2^2 6^2
 &    
    ( 4, 8, -8 )  ,\,  
    ( 6, -2, 6 )
 & *
 & \textcolor{red}{\ref{argpower}},\,  \textcolor{red}{\ref{elim1}}
\\

( 12, 1, 5 ) 
 & 1.2^2.3.12^2/4^2
 &  4
 &  3^1 6^2
 &    
    ( 4, 4, 8 ) ,\,
    ( 6, 5, -3 ) 
 & {\rm No}
 & -
\\

( 12, 0, 2 ) 
 &2.4.6.12 
 &  4
 &  2^2 12^2
 &    ( 4, 0, 8 )  ,\,
    ( 6, 2, 6 )
 & *
 & \textcolor{red}{\ref{argpower}}
\\

( 12, -2, 2 ) 
 &2^23^2.4.12/1^2
 &  4
 &  2^1 6^1 12^2
 &     ( 4, 4, 8 ) ,\,
    ( 6, 2, 6 ) 
 & *
 &  \textcolor{red}{\ref{argfix12}}
\\  \hline

( 14, 1, 3 ) 
 & 1.2.7.14
 &  4
 &  14^2
 & 
    ( 2, 8, 24 ) ,\,
    ( 7, 3, 3 ) 
 & *
 & -
\\ \hline

( 15, 1, 1 ) 
 & 1.3.5.15
 &  4
 &  15^2
 &       ( 3, 6, 6 ) ,\,
         ( 5, 4, 4 ) 
 & *
 & -
\end{array}
$$}
\end{table}

Since the triple $(\hbox{order of}\ g,\,{\rm Trace} (g),\,{\rm Trace} (g^2))$ uniquely determines the ${\rm Co}_0$
conjugacy class, the ${\rm Co}_0$ conjugacy class $[g]$ associated to a finite symplectic automorphism $g$
is uniquely determined.

\medskip

The following observation is clear:
\begin{rem}\label{argpower}
If a conjugacy class $[g]$ is inadmissible, then so is $[h]$ whenever 
$g$ is a power of $h$. 
\end{rem}

The lattice-theoretic set-up leads to a condition on the discriminant group. 
\begin{lem}\label{argdiscriminant}
Let $A_{\Lambda^g}=(\Lambda^g)^*/\Lambda^g$ be the discriminant group of 
$\Lambda^g$ with quadratic form $q_{\Lambda^g}:A_{\Lambda^g}\rightarrow \Q/2\Z $.
Suppose that $\rk(A_{\Lambda^g})=\rk(\Lambda^g).$\ Then the following hold:
\begin{itemize}
\itemsep0em
\item[(a)] The discriminant group of $L^g$ has index $2$ in $A_{\Lambda^g}$.
\item[(b)] $q_{\Lambda^g}$ has one of the values $\frac{1}{2}$, $\frac{3}{2}$ in its image. 
\end{itemize}
\end{lem}
We defer the proof until Section~\ref{subol}.

\begin{lem}\label{elim1} Conjugacy classes of type  $(2, -8, 24)$, $(3,0,0)$, $(4,0,0)$, $(4, 8, -8)$, 
$ (6,-1,-3)$, $(6,0,0)$, $(6, -4, 6)$, $(6, -2, 6)$, $(6, 4, 6)$, $(8,0,8)$,
$(8,0,0)$, $(10, -2, 4)$, $(10, 2, 4)$ and $(12, 2, -2)$ are inadmissible.
\end{lem}
\pf
Inspection of the structure of $A_{\Lambda^g}$ in Table~\ref{Leechclasses}, together with 
Lemma~\ref{argdiscriminant}~(a) eliminates types 
$(3, 0, 0)$, $(4, 0, 0)$, $(6, -1, -3)$, $(6, 0, 0)$, $(8,0,8)$ and $(8,0,0)$.\ 
The additional types $(2, -8, 24)$, $(4, 8, -8)$,  $(6, -4, 6)$, $(6, -2, 6)$,
$( 6, 4, 6 )$, $( 10, -2, 4 )$, $(10, 2, 4)$, and $( 12, 2, -2 )$ are excluded
by Lemma~\ref{argdiscriminant}~(b) since a computer calculation shows that
$q_{\Lambda^g}$ has neither the value $\frac{1}{2}$ nor $\frac{3}{2}$ in its image.
\qed

\medskip
Suppose that $g$ is a finite symplectic automorphism of a hyperk\"{a}hler manifold of type $K3^{[2]}$.\ 
We say that $g$ is of \emph{$K3$-type} if there is a symplectic automorphism $h$ of a $K3$ surface
such that $g$ is conjugate in ${\rm Co}_0$ to the element defined by $h$.\  We also say that the elements and conjugacy classes
in ${\rm Co}_0$ corresponding to $g$ are themselves of $K3$-type.\  For symplectic automorphisms of K3 surfaces we have $\rk(\Lambda^g)\geq 5$ \cite{Mu}, while
the  analog of part a)  Lemma~\ref{argdiscriminant} is $\rk(A_{\Lambda^g})\leq \rk(\Lambda^g)-2$.\
Then examination of Table \ref{Leechclasses} establishes the next Remark:
\begin{rem}
There are $8$ conjugacy classes in ${\rm Co}_0$ of $K3$-type, namely  
$(1,24,24)$, $(2,8,24)$, $(3,6,6)$, $(4,4,8)$, $(5,4,4)$, $(6,2,6)$, $(7,3,3)$ and $(8,2,4)$. \qed
\end{rem}
Nikulin \cite{Nikulin}  first  proved  that the order of a (finite order) symplectic automorphism of $K3$ is
at most $8$.\ See also \cite{Mason}, \cite{Mu}.

\medskip

To exclude further cases beyond Lemma \ref{elim1},  we apply the equivariant Atiyah-Singer theorem to the Hirzebruch $\chi_y$-genus 
of the hyperk\"{a}hler manifold $X$ of type $K3^{[2]}$.\ Let $g\in {\rm Aut}(X)$ have finite order $n$ and let 
$$\chi_y(g;X):=\sum_{p,\,q=0}^4(-1)^q\,\tr(g|H^{p,q}(X))\,y^p$$ 
be the equivariant $\chi_y$-genus.

\begin{lem}\label{chiy}
Let $t=\tr(g|H^{1,1}(X))$ and $s=\tr(g^2|H^{1,1}(X))$.\ Then
$$\chi_y(g;X)=3-2t\, y+\frac{6+t^2+s}{2}\,y^2-2t\, y^3 + 3\,y^4.$$
\end{lem}
\noindent
{\bf Proof:}\ Inspection of the Hodge diamond of $X$ shows that the only
nontrivial contributions one has to know are those coming from $H^{1,1}(X)$ and $H^{2,2}(X)$.\ 
The remainder then follow from the symmetries of the Hodge diamond, which holds equivariantly.\ But 
$H^{2,2}(X)\cong \C\oplus S^2H^{1,1}(X)$.\ Together with the formula
for the character of a symmetric square, this gives the result.\ 
For further details, see Camere~\cite{Camere}. $\hfill \Box$

\medskip
We note that $\tr(g| \Lambda)=\tr(g|H^{1,1}(X))+3$.\ Moreover 
from Table~\ref{Leechclasses}, we see that $t$ and $s$ are rational integers.

\medskip

There is a basic result regarding the structure of 
the fixed-point set $X^g$.
\begin{thm}[cf.\ Camere \cite{Camere}, Proposition 3]
The fixed-point set of a finite symplectic automorphism of a compact hyperk\"ahler manifold
is the disjoint union of finitely many components, which are themselves hyperk\"ahler manifolds.\
The centralizer of such an automorphism acts by symplectic automorphisms on the fixed-point set. \qed
\end{thm}
Since the only connected $4$-dimensional hyperk\"ahler manifolds are K3-surfaces and complex $2$-tori,
it follows that the fixed-point set of a non-trivial finite symplectic automorphism on an $8$-dimensional
hyperk\"ahler manifold consists of isolated points,  complex $2$-tori and K3 surfaces. 

The occurrence of $2$-tori can sometimes be excluded by the following  
geometric result of Mongardi.
\begin{thm}[Mongardi~\cite{Mon-thesis}, Proposition 5.1.4]\label{argmongardi}
Let $g$ be a symplectic automorphism of finite order of a hyperk\"ahler manifold $X$ of 
type $\K32$, and suppose  that $X^g$ contains a torus.\ Then $\rk(L^g)\leq 6$.
\end{thm}
\begin{cor}\label{Corinvfp} Suppose that $g$ lies in one of the Conway classes
$(2,8,24)$, $(2,0,24)$, $(3,6,6)$, $(3,0,0)$, $(4,8,-8)$, $(4,4,8)$, $(4,0,8)$, $(5,4,4)$ or $(6,2,6)$.\
Then the components of $X^g$ are isolated fixed-points or $K3$ surfaces.
\end{cor}
\noindent
\pf For all of these choices of $g$ one has  $\rk L^g =\rk \Lambda^g-1\geq 7$ 
(cf.\ column 3 of Table~\ref{Leechclasses}) and the Theorem applies. \qed

\smallskip

To compute  $\chi_y(g;X)$ using the equivariant  Atiyah-Singer index theorem, we have to know
the possible eigenvalues for the action of $g$ on the normal bundle in $X$ of a component $F$ of $X^g$.\
Let $\zeta=e^{2\pi i/n}$ (where $g$ has order $n$).

Since the structure group of the tangent bundle of $X$ can be reduced to ${\rm Sp}(2)\subset {\rm SU}(4)\subset {\rm U}(4)$
there are the following possibilities:
\begin{itemize}
\item[] $F$ is an isolated fixed-point $p$.\ The possible eigenvalues for $g$ acting on $T_pX$ are 
$(\zeta^i,\zeta^{-i},\zeta^j,\zeta^{-j})$, $0<i\leq j \leq n/2$.
\item[] $F$ is a K3 surface or $2$-torus.\ The possible eigenvalues for $g$ acting on the normal
bundle $N$ of $F$ in $X$ are $(\zeta^i,\zeta^{-i})$, $0<i\leq n/2$.
\end{itemize}
By the equivariant fixed-point theorem~(cf.~\cite{Hi-book}) one has the following formula:
\begin{equation}\label{ASI}
\chi_y(g;X)=\sum_{F\subset X^g}\prod_{k=1}^{\dim_{\C} F}\frac{x_k(1+y e^{-x_k})}{1- e^{-x_k}}
\prod_{k=1}^{4-\dim F_{\C}}\frac{1+y\lambda_k e^{-x_k'}}{1-\lambda_k e^{-x_k'}}\,[F],
\end{equation}
the sum running over  the components $F$ of $X^g$.\ 
The $x_k$ and $x_k'$ are the formal roots of the total Chern classes of the tangent and
normal bundle of $F$ respectively, and the $\lambda_k$ are the eigenvalues of $g$ acting on the
normal bundle.\ We will evaluate the right-hand-side for each type of fixed-point component 
and $g$-action on the normal bundle.

For an isolated  fixed-point $p$ of type $(\zeta^i,\zeta^{-i},\zeta^j,\zeta^{-j})$, 
the contribution is 
$$f_{i,j}:=\frac{1+y\zeta^i }{1-\zeta^i}\cdot \frac{1+y\zeta^{-i} }{1-\zeta^{-i}}\cdot 
  \frac{1+y\zeta^j }{1-\zeta^j}\cdot \frac{1+y\zeta^{-j} }{1-\zeta^{-j}}.$$

For a complex $2$-dimensional surface $F$, the total Chern class of $F$ is $c(TF)=(1+x_1)(1+x_2)=1+c_2(TF)$
since $c_1(TF)=0$.\  A short calculation expanding  $e^{-x_k}$ up to
order $2$ shows that the contribution from the tangent bundle in~(\ref{ASI}) is the factor 
$$\label{tangent}
h_0:=\frac{x_1(1+y e^{-x_1})}{1- e^{-x_1}} \cdot \frac{x_2(1+y e^{-x_2})}{1- e^{-x_2}}
=(1+y)^2+(2-20y+2y^2)\cdot \frac{c_2(F)}{12}.
$$
For the normal bundle $N$ with $g$ acting with eigenvalues
$(\zeta^i,\zeta^{-i}$), we have to express
$$h_i:=\frac{1+y\zeta^i e^{-x_1'}}{1-\zeta^i  e^{-x_1'}}\cdot
\frac{1+y\zeta^{-i}e^{-x_2'}}{1-\zeta^{-i} e^{-x_2'}}$$
in the total Chern class $c(N)=(1+x_1')(1+x_2')$ of the normal bundle.\
One obtains
$$\label{normal}
h_i=  -\frac{\zeta ^i+y+  \zeta^{2 i} y + \zeta^i y^2 }{\left(\zeta ^i-1\right)^2}
   -\frac{ \zeta ^i (\zeta ^i+1)(y+1)^2}{\left(\zeta ^i-1\right)^3}\cdot x_1'
   -\frac{ \zeta ^i (\zeta ^{2 i}+4 \zeta ^i+1)(y+1)^2}
         {2 \left(\zeta ^i-1\right)^4}\cdot (x_1')^2 ,
$$
where we have also used that $x_1'+x_2'=c_1(N)=c_1(TX)|_F-c_1(TF)=0$.\
Note that if $g$ is an involution, the linear term for $x_1'$ vanishes.\
Otherwise, $N$ splits canonically into two eigenspace bundles and $x_1'$
is well-defined in this case.

Assume that there are $a_{i,j}$ isolated fixed-points of type $(\zeta^i,\zeta^{-i},\zeta^j,\zeta^{-j})$
and $b_i$ fixed-point components which are K3 surfaces of type $(\zeta^i,\zeta^{-i})$.\ 
The right hand side of~(\ref{ASI}) equals
$$\sum_{i,\,j} a_{i,j} \cdot f_{i,j}
\ +\  \sum_{i} \sum_{ F\subset \Phi_i} h_0 h_i[F] 
\ +\  \sum_{i} \sum_{F\subset \Psi_i } h_0 h_i[F] ,$$
where $\Phi_i$ and $\Psi_i$ denote the union of fixed-point components $F\subset X^g$ which are
K3-surfaces resp.\ $2$-tori of type $(\zeta^i,\zeta^{-i})$.\
Using $-(x_1')^2=c_2(N)=c_2(TX)|_F-c_2(TF)$ and $c_2(TF)[F]=24$ for $F$ a K3 surface resp.~$c_2(TF)[F]=0$ for $F$ a $2$-tori,
we see that for fixed $i$, the sum $ \sum_{F\subset \Phi_i\cup \Psi_i} h_0 h_i[F]$ depends only on $b_i$
and the sum  $C_i:=\sum_{F\subset \Phi_i\cup \Psi_i}c_2(TX|_F)[F]$.

\smallskip

Thus (\ref{ASI})
becomes
\begin{equation}\label{ASIisolated}
\chi_y(g;X)=\sum_{ 0<i\leq j \leq n/2}a_{i,j}\cdot f_{i,j} + \sum_{0<i\leq n/2} b_i\cdot\beta_i + \sum_{0<i\leq n/2} C_i\cdot\gamma_i,
\end{equation}
with explicit polynomials $\beta_i$ and $\gamma_i$ in $y$ and $\zeta$.\
This gives $3\,\varphi(n)$ linear equations (possibly trivial and linearly dependent)
for the integers $a_{i,j}$, $b_i$ and $C_i$ since there  are $3$ independent rational coefficients in the
palindromic polynomial $\chi_y(g;X)$, and the right-hand-side is a polynomial in $y$ with coefficients in 
the cyclotomic field $\Q(\zeta)$ of degree $\varphi(n)$ over $\Q$.\
In addition, the $a_{i,j}$ and $b_i$ are non-negative.

\smallskip

If $h$ is a power of  $g$ then the fixed-point configurations of $g$ and $h$ and
the actions on the normal bundles are related.\ As before, let $n$ be the order of $g$ and
let $h=g^k$ for $k|n$, $k<n$.\ Consider an isolated fixed-point $p$ for which $g$
acts with eigenvalues $(\zeta^i,\zeta^{-i},\zeta^j,\zeta^{-j})$ in the tangent space.\
Then $p$ is also a fixed-point of $h$ and $h$
acts with eigenvalues $(\zeta^{ik},\zeta^{-ik},\zeta^{jk},\zeta^{-jk})$
in the tangent space.\ If both, $\zeta^{ik}$ and $\zeta^{jk}$, are different
from $1$ then $p$ is also an isolated fixed-point of $h$.\ If one of them 
is equal to~$1$, then $p$ belongs to a $4$-dimensional fixed-point set, i.e.~a
K3 surface or a $2$-torus.\ The case that  $\zeta^{ik}=\zeta^{jk}=1$ is
impossible, since otherwise $h=1$.\ 
If $p$ is a fixed-point of $g$
belonging to a higher-dimensional fixed-point component $F$ of $g$, then $h$ acts with the $(\zeta^{ik},\zeta^{-ik})$
in the normal bundle and $\zeta^{ik}$ must be different from $1$, 
i.e.,  $\zeta^{i}$ is necessarily a primitive $n$-th root of unity.
Note that $h$ can have additional fixed-point components besides the one described above.

It will turn out that for all classes of Table~\ref{Leechclasses}, there is at most
one possible fixed-point configuration.\ Thus by analyzing the fixed-point 
structure and the action of $g$ on the normal bundle, we can apply the information
previously obtained to all  non-trivial powers of $g$.\ This gives several restrictions on
the possible fixed-point components and the eigenvalues.

\smallskip

For a given conjugacy class, our approach now is to solve the resulting system of linear equations and inequalities.\
This can be done in a straightforward way with the help of a computer, although the number of equations and variables
will become quite large.\ 
From these calculations, together with some additional geometric results, we obtain
the following main theorem:
\begin{thm}\label{allowedclasses}
A symplectic automorphism  $g$ of finite order of a hyperk\"ahler manifold of type $\K32$ belongs
to one of the $15$ ${\rm Co}_0$ conjugacy classes $(1,24,24)$,
$(2,8,24)$, $(3,6,6)$, $(3,-3,-3)$, $(4,4,8)$, $(5,4,4)$, $(6,2,6)$, $(6,5,-3)$, $(7,3,3)$, $(8,2,4)$, $(9,3,3)$,
$(11,2,2)$, $(12,1,5)$, $(14,1,3)$, $(15,1,1)$.\ If $g$ is of type $(2,8,24)$, the fixed-point set contains a unique $K3$-surface;
if $g$ is of type $(3,-3,-3)$, the fixed-point set contains a unique $2$-torus;
for all other $g\not=1$, the fixed-point set consists of  isolated fixed-points.\ The complete
description of the fixed-point sets and the action on the normal bundles is given
in Table~\ref{Fixpointset}.
\end{thm}
\begin{table}[t]
\caption{Admissible classes and corresponding fixed point configurations}\label{Fixpointset}
$$\begin{array}{lcl}
\hbox{class of $g$} & \hbox{\# of components of a certain type} & \hbox{prime powers}  \\ \hline\hline
(1,24,24)^* &  X &  \\ \hline
(2,8,24)^*  & 28 \times(-1,-1),\  K3 &  \\\hline
(3,6,6)^*  & 27\times (\zeta_3,\zeta_3) & \\
(3,-3,-3)^{\dagger}  & T^2 & \\ \hline
(4,4,8)^*  & 8 \times [(i,i),\, (i,-1)] & (2,16,8)^* \\ \hline
(5,4,4)^*   & (\zeta_5,\zeta_5),\, (\zeta_5^2,\zeta_5^2),\, 12\times (\zeta_5,\zeta_5^2) \\ \hline
(6,2,6)^*   & (\zeta_6,\zeta_6),\,  6\times (\zeta_6,\zeta_6^2) & (2,8,24)^*,\,(3,6,6)^* \\ 
(6,5,-3)^\dagger   & 10\times (\zeta_6,\zeta_6^3),\, 6\times (\zeta_6^2,\zeta_6^3) & (2,8,24)^*,\,(3,-3,-3)^\dagger  \\ \hline
(7,3,3)^*   & 3 \times [(\zeta_7,\zeta_7^2)\,(\zeta_7,\zeta_7^3)\,(\zeta_7^2,\zeta_7^3)]\\ \hline
(8,2,4)^*   & 2 \times [(i,\zeta_8),\,(i,\zeta_8^3),\, (\zeta_8,\zeta_8^3)] & (4,4,8)^*\\ \hline
(9,3,3)^{\dagger}  & 3 \times [(\zeta_9,\zeta_9^3),\,(\zeta_9^2,\zeta_9^3),\,(\zeta_9^3,\zeta_9^4)]  & (3,-3,-3)^\dagger  \\ \hline
(11,2,2)  &  (\zeta_{11},\zeta_{11}^3),\,(\zeta_{11},\zeta_{11}^4),\,(\zeta_{11}^2,\zeta_{11}^3),
                 \,(\zeta_{11}^2,\zeta_{11}^5),\,(\zeta_{11}^4,\zeta_{11}^5)\\ \hline
(12,1,5)^\dagger   &  (\zeta_{12},\zeta_{12}^3),\, (\zeta_{12}^3,\zeta_{12}^5),\,  2 \times (\zeta_{12}^2,\zeta_{12}^3) & (6,5,-3)^\dagger,\,(4,4,8)^*  \\ \hline
(14,1,3)    &  (\zeta_{14},\zeta_{14}^4),\,  (\zeta_{14}^2,\zeta_{14}^3),\,  (\zeta_{14}^5,\zeta_{14}^6) & (2,8,24)^*,\,(7,3,3)^*   \\ \hline
(15,1,1)    &  (\zeta_{15},\zeta_{15}^4),\, (\zeta_{15}^2,\zeta_{15}^7)       & (3,6,6)^*,\,(5,4,4)^*   \\ \hline
\end{array}$$

{\it Notation:\/} ${}^*$ means element is of K3-type and contained in $M_{24}$; ${}^\dagger$ means element is not in $2^{12}{:}M_{24}$.
\end{table}

In the following, we  discuss the proof  in more detail.\

\smallskip

For involutions, the fixed-point formula was first used by Camere.\
She obtained the following result, which we verified with our program.
\begin{thm}[Camere~\cite{Camere}]\label{argcamere}
Let $g$ be a symplectic involution of a hyperk\"ahler manifold $X$ of type $\K32$.\
Then $g$ is of type $(2,0,24)$, $(2,6,24)$ or $(2,8,24)$ and the corresponding fixed-point sets are as follows:

\vspace{-2mm}
\begin{itemize}
\itemsep0em
\item[] $(2,0,24)$: $12$ isolated fixed-points and at least one complex torus, 
\item[] $(2,6,24)$: $36$ isolated fixed-points and at least one complex torus,
\item[] $(2,8,24)$: $28$ isolated fixed-points, one K3 surface and an undetermined number of complex tori.  $\hfill \Box$
\end{itemize}
\end{thm}
\noindent
Combining this with the other information, we obtain
\begin{lem}[Mongardi~\cite{Mon-thesis}, Theorem 6.2.3]\label{invtype} 
The symplectic automorphisms of order~$2$ have type $(2, 8, 24)$.
\end{lem}
\pf
The case $(2,-8,24)$ of Table~\ref{Leechclasses} cannot occur by Camere's theorem,
$(2,6,24)$ cannot occur since it is absent from Table~\ref{Leechclasses}, and 
$(2,0,24)$ is excluded by Corollary~\ref{Corinvfp}.\
The only remaining possibility from Table~\ref{Leechclasses} is $(2,8,24)$.
$\hfill \Box$

Another application of
Corollary~\ref{Corinvfp} shows if $g$ is an involution, then the components of $X^g$ are 
either isolated fixed-points or $K3$ surfaces.\
Therefore we have:
\begin{lem}\label{argpoweriso}
If $g$ is a symplectic automorphism of $X$ of even order, then the components of $X^g$ are either isolated fixed-points
or $K3$ surfaces. $\hfill \Box$
\end{lem}

Next we consider symplectic automorphisms of order $4$.
\begin{prop}\label{argfix4}
Let $g$ be an order $4$ symplectic automorphism of a hyperk\"ahler manifold $X$ of type $\K32$.\
Then $g$ is of type $(4,4,8)$,
and $X^g$ consists of $16$ isolated fixed-points.\ 
There are $8$ fixed-points with eigenvalues $(i,-i,i,-i)$ and 
$8$ fixed-points with eigenvalues  $(i,-i,-1,-1)$.
\end{prop}
\pf Since $g^2$ has order $2$, we  know from Lemma \ref{invtype}  that $g^2$ has type $(2,8,24)$, whence $g$ is of
type $(4,4,8)$, $(4,0,8)$ or $(4,-4,8)$ (cf.\ Table \ref{Leechclasses}). By Lemma~\ref{argpoweriso}, the only fixed-point
components are isolated fixed-points or K3 surfaces. Since $g^2$ is not the identity we have 
$a_{2,2}=b_{2}=0$.\ Thus we must solve eqn.\ (\ref{ASIisolated}) 
for the four variables $a_{1,1}$, $a_{1,2}$, $b_1$ and $C_1$.
For the right-hand-side of (\ref{ASIisolated}) we obtain
$$   
 \frac{a_{1, 1}}{4}(1 + y^2)^2 + 
\frac{a_{1, 2}}{8}(-1 + y)^2(1 + y^2) -
b_1 (11(1+y^4) + 58(y+y^3) + 70y^2 ) + \frac{C_1}{8} (1 + y)^4.
$$
The left-hand-side is given  by Lemma~\ref{chiy}, where the value for $t$ and $s$ can
be read off from the type $(n,t+3,s+3)$ of $g$.\ 

For $g$ of type $(4,4,8)$ one obtains  three linear equations
\begin{eqnarray*}\textstyle
\frac{1}{4}\, a_{1,1}+\frac{1}{8}\, a_{1,2}-11\, b_1+\frac{1}{2}\,C_1 &= &\phantom{-}3, \\ \textstyle
-\frac{1}{4}\, a_{1,2}-58\, b_1+2\, C_1 & =  & -2,\\ \textstyle
\frac{1}{2}\, a_{1,1}+\frac{1}{4}\, a_{1,2}-70\, b_1 + 3\, C_1 & = &\phantom{-} 6
\end{eqnarray*}
with the solutions 
$a_{1,1}=\frac{2}{3} (C_1+12)$, $a_{1,2}= \frac{1}{3} (24-5 C_1)$,
and $b_1= \frac{1}{24}C_1$.\
Using the inequalities $a_{1,1}\geq 0$, $a_{1,2}\geq 0$, $b_1\geq 0$ for 
integral $C_1$ shows that $C_1\in\{0,1,2,3,4\}$.\ Only for $C_1=0$ we
obtain integer solutions  $a_{1,1}=a_{1,2}=8$ and $b_1=0$.\ In particular,
$g$ must have isolated fixed-points.

The same approach for $g$ of type $(4,0,8)$ or $(4,-4,8)$ gives no solutions.
(Moreover, this also holds for the cases of type $(4,8,-8)$, $(4,0,-8)$,
and $(4,0,0)$, though they are already excluded.)\ This completes the proof of
the proposition. $\hfill \Box$

We remark that the $8$ fixed-points of an order four element $g$ with eigenvalues $(i,-i,-1,-1)$
necessarily lie on the $K3$ fixed-point component of $g^2$, and $g$ acts
on this K3 surface with $8$ isolated fixed-points.\ This is in agreement with results of Nikulin
and Mukai for symplectic automorphisms of K3 surfaces.

\medskip
A similar approach will handle the elements of order $8$:
\begin{prop}\label{argfix8}
Let $g$ be a symplectic automorphism of a hyperk\"ahler manifold $X$ of type $\K32$
of order~$8$.\ Then $g$ is of K3-type $(8,2,4)$, acting with $6$ isolated fixed-points and
eigenvalues as in Table~\ref{Fixpointset}.
\end{prop}
\pf By Proposition \ref{argfix4}, $g^2$ has type $(4,4,8)$, whence $g$ is of
type $(8,2,4)$ or $(8,-2,4)$ by Table \ref{Leechclasses}.\ 
Moreover, since $g^2$ has to act with isolated fixed-points, the same 
is true for $g$, and $a_{1,4}=a_{2,4}=a_{3,4}=a_{4,4}=0$.\ Since $g^4\not=1$
then $a_{2,2}=0$.\ Solving eqn.\ (\ref{ASIisolated}) for the remaining variables
$a_{1, 1}$, $a_{1, 2}$, $a_{1, 3}$, $a_{2, 3}$, and $a_{3, 3}$ gives
for $g$ of type $(8,2,4)$ a unique solution 
$a_{1,1}=0$, $a_{1, 2}=2$, $a_{1, 3}=2$, $a_{2, 3}= 2$, $a_{3, 3} = 0$, 
and for $g$ of type $(8,2,-4)$  a unique solution
$a_{1,1}=1$, $a_{1, 2}=2$, $a_{1, 3}=-4$, $a_{2, 3}= 2$, $a_{3, 3} = 1$.\ The latter case is impossible
since $a_{1, 3}$ is negative, and the Proposition is proved. $\hfill \Box$

Note that for an element $g$ of order $8$,  $g^4$ is an involution which has a K3 surface as
a fixed-point component on which $g$ acts with $4$ fixed-points
with eigenvalues $(i,-i)$, as required. 

\medskip

The cases when $g$ has order $3$, $5$, $7$ or $11$ are covered by
Mongardi \cite{Mon-thesis} (the fixed-point sets are also described) and we merely state the result in these cases.\
By our method based on equation~(\ref{ASIisolated}), we can determine the eigenvalues for the action of $g$ on the
normal bundle in all cases.\ Mongardi's result for $g$ of order~$3$ is as follows:
\begin{prop}[Mongardi \cite{Mon-thesis}, Theorem 6.2.4, Proposition 6.2.8]\label{argfix3}
The admissible elements $g$ of order $3$ are the K3-type $(3,6,6)$ acting with $27$ isolated fixed-points, 
and type $(3,-3,-3)$ acting with a complex $2$-torus as fixed-point set and  
eigenvalues as in Table \ref{Fixpointset}. $\hfill \Box$
\end{prop}
Note that  type  $(3,0,0)$ is inadmissible by Lemma \ref{elim1}.\
Mongardi's proof for the type $(3,-3,-3)$ (and also if $g$ has order $7$ or $11$) uses Theorem~1.2 of~\cite{BNS}.\
The eigenvalues are all uniquely determined.

At this point, application of Remark~\ref{argpower}, Lemmas~\ref{elim1} and~\ref{invtype}, 
and Propositions~\ref{argfix4} and~\ref{argfix8} leaves only types $(5, -1, -1)$, 
$(10, 3, -1)$ and $(12, -2, 2)$ from Table~\ref{Leechclasses} to be eliminated.

To deal with the prime $5$ we use:
\begin{prop}[Mongardi \cite{Mon-thesis}, Thm.~6.2.9]\label{argfix5}
Let $g$ be a symplectic automorphism of a hyperk\"ahler manifold of type $\K32$
of order~$5$.\ Then $g$ is of type $(5,4,4)$ acting with $14$ isolated fixed-points and
eigenvalues as in Table~\ref{Fixpointset}. 
\end{prop}
The eigenvalues for $g$ have been implicitly determined in the proof of \cite{Mon-thesis}, Thm.~6.2.9.
We confirmed the calculation with our computer program.

Mongardi's result means that elements of type $(5, -1, -1)$ are inadmissible, and  those of type $(10, 3, -1)$
are then also inadmissible because they have squares of type $(5, -1, -1)$ (cf.\ Table \ref{Leechclasses}).
\begin{prop}\label{argfix12}
Let $g$ be a symplectic automorphism of a hyperk\"ahler manifold of type $\K32$
of order~$12$.\ Then $g$ is of type $(12,1,5)$ acting with $4$ isolated fixed-points and
eigenvalues as in Table~\ref{Fixpointset}.
\end{prop}
\pf
The only conjugacy classes for $g$ having $g^3$ of type $(4,4,8)$ and $g^2$ of type
$(6,2,6)$ or $(6,5,-3)$ are those of type $(12,-2,2)$ and $(12,1,5)$, respectively.\ Since $g^3$ acts
by isolated fixed-points the same holds for $g$.\ Since $g^2$ acts by isolated fixed-points 
(see Proposition~\ref{argfix6} below) we have $a_{i,6}=0\ (1\leq i\leq 6)$.\
Since $g^6\not= 1$,  $a_{2,2}$, $a_{2,4}$, $a_{4,4}$ also vanish,
and since $g^4\not= 1$ then $a_{3,3}=0$.\
In addition, if $g$ is of type $(12,-2,2)$ then $g^4$ has only isolated fixed-points
and $a_{1,3}=a_{2,3}=a_{3,4}=a_{3,5}=0$.

For $g$ of type $(12,-2,2)$, equation (\ref{ASIisolated}) has no solution such that the remaining
seven variables are non-negative.\ For $g$ of type  $(12,1,5)$, the only solution such that the remaining
eleven variables are non-negative and integral is $a_{1, 3}=1$, $a_{3, 5}=1$, $a_{2, 3}=2$ and $0$ 
for the other eight variables.\ The proposition follows.  $\hfill \Box$

For an order $12$ element $g$, its power $g^6$ is an involution with a K3 surface as
one of the  fixed-point components.\  On this, $g$ acts with $2$ fixed-points, as required.\
The power $g^4$  is an order $3$ element and has a $2$-torus as
one fixed-point component on which $g$ acts with $4$ fixed-points. 

This completes the proof that the admissible classes of symplectic automorphisms are those listed 
in Theorem~\ref{allowedclasses}.

\smallskip

Next we discuss the fixed-point configuration for the classes in
Table~\ref{Fixpointset} not yet covered.\

\begin{prop}\label{argfix6}
The admissible $g$ of order $6$ are the $K3$-type $(6,2,6)$ and 
type $(6,5,-3)$, each acting with isolated fixed-points and 
eigenvalues as in Table \ref{Fixpointset}.
\end{prop}

\pf We have already shown that the only possibilities are types $(6,2,6)$ and  $(6,5,-3)$.

In the first case, $g^2$ is of type $(3,6,6)$ acting with isolated fixed points.\
Hence $g$ must act by isolated fixed points, too. In the second case,  $g^2$ is of type $(3,-3,-3)$
and has a $2$-torus as fixed-point set. Furthermore,  $g^3$ is 
of type $(2,8,24)$ and contains only isolated fixed points or a K3 surface in its fixed-point set. 
Thus $g$ can only have isolated fixed points.

We also see that the numbers $a_{2,2}=a_{3,3}=0$ for $g$
(otherwise $g^2$ or $g^3$ would be the identity). 

\smallskip

For $g$ of type $(6,2,6)$, 
equation (\ref{ASIisolated}) has the solution 
$a_{1, 1}=  1 + a_{1, 3}/8$, $a_{1, 2} = 6 - 9\,a_{1, 3}/8$, $a_{2, 3}=0$,
implying that $a_{1, 3}=0$, $a_{1,1}=1$, $a_{1,2}=6$.\
For type $g$ of type $(6,2,6)$, we get the solution
$a_{1, 1}= -5/4 + a_{1, 3}/8$, $a_{1, 2}=  45/4 - 9\, a_{1, 3}/8$ and $a_{2, 3}=  6$,
which implies $a_{1, 3}=10$, $a_{1,1}=0$ and $a_{1,2}=0$. $\hfill \Box$

In both cases,  $g^3$ has a K3 surface as a fixed-point component
on which $g$ acts with $6$ isolated fixed-points, as required.\
For  $g$ of type $(6,5,-3)$,  $g^2$ has a $2$-torus as
one fixed-point component on which $g$ acts with $16$ isolated fixed-points.

\begin{prop}\label{argfix9}
An admissible elements of order $9$ has type $(9,3,3)$, acting with~$9$ isolated fixed-points and
eigenvalues as in Table \ref{Fixpointset}.
\end{prop}
\pf From Table~\ref{Leechclasses} there is only one possible type of element $g$ of order $9$.

Since $g^3$ is of type  $(3,-3,-3)$ with a $2$-torus as fixed-point set, we have to consider the 
nine variables $a_{1,3}$, $a_{2,3}$, $a_{3,4}$, $b_1$,  $b_2$,  $b_4$, $C_1$, $C_2$ and $C_4$
which might be nonzero.\  
Equation~(\ref{ASIisolated}) provides us with the unique solution $a_{1, 3} = a_{2, 3} = a_{3, 4} = 3$,
$b_1=b_2=b_4=0$ and $C_1=C_2=C_4=0$ for non-negative $a_{i,j}$ and $b_i$. $\hfill \Box$

\begin{prop}\label{argfix15}
An admissible element of order $15$ has type $(15,1,1)$, acting with $2$ isolated fixed-points and
eigenvalues as in Table \ref{Fixpointset}.
\end{prop}

\pf The statement about the type is clear.\ By considering $g^5$ of type $(3,6,6)$ and 
$g^3$ of type $(5,4,4)$ it is also clear that $g$ acts with isolated fixed-points and 
$a_{1,3}$, $a_{2,3}$, $a_{3,3}$, $a_{3,4}$, $a_{3,5}$, $a_{3,6}$, $a_{3,7}$, $a_{1,6}$,
$a_{2,6}$, $a_{4,6}$, $a_{5,6}$, $a_{6,6}$, $a_{6,7}$
and $a_{1,5}$, $a_{2,5}$, $a_{3,5}$, $a_{4,5}$, $a_{5,5}$, $a_{5,6}$, $a_{5,7}$
all vanish.\ Equation (\ref{ASIisolated}) has the general solution 
$a_{1, 1} =(-a_{4, 7} - a_{7, 7})/5$, 
$a_{1, 2} = a_{4, 7}$, 
$a_{1, 4} = (5 - 4\,a_{4, 7} + 26\,a_{7, 7})/5$, 
$a_{1, 7} = (a_{4, 7} - 24\,a_{7, 7})/5$, 
$a_{2, 2} = a_{7, 7}$, 
$a_{2, 4} = (a_{4, 7} - 24\,a_{7, 7})/5$, 
$a_{2, 7} = (5 - 6\,a_{4, 7} + 14\,a_{7, 7})/5$, 
$a_{4, 4} = (-a_{4, 7} - a_{7, 7})/5$.\
The only non-negative integral solution occurs when
$a_{4, 7} = a_{7, 7} = 0$, leading to the eigenvalues as in Table~\ref{Fixpointset}. \qed

The two fixed-points of $g$ are the two distinguished fixed-points of $g^3$ with
eigenvalues $(\zeta_5,\zeta_5^4,\zeta_5,\zeta_5^4)$ and $(\zeta_5^2,\zeta_5^3,\zeta_5^2,\zeta_5^3)$.

\begin{prop}[Mongardi \cite{Mon-thesis}, Proposition 6.2.15]\label{argfix7}
An admissible element
of order $7$ has type $(7,3,3)$, acting with $9$ isolated fixed-points with
eigenvalues as in Table \ref{Fixpointset}. $\hfill \Box$
\end{prop}

\begin{prop}\label{argfix14}
An admissible element $g$
of order $14$ has type $(14,1,3)$, acting with $3$ isolated fixed-points and
eigenvalues as in Table \ref{Fixpointset}.
\end{prop}
\pf 
The statement about the type is again clear.\ By considering $g^7$ of type $(2,8,24)$ and 
$g^2$ of type $(7,3,3)$ it follows that equation~(\ref{ASIisolated}) contains the variables
$a_{1, 1}$, $a_{1, 2}$, $a_{1, 3}$, $a_{1, 4}$, $a_{1, 5}$, $a_{1, 6}$, $a_{2, 3}$,  
$a_{2, 5}$, $a_{3, 3}$, $a_{3, 4}$, $a_{3, 5}$, $a_{3, 6}$, $a_{4, 5}$, $a_{5, 5}$, $a_{5, 6}$.\
The solution of the corresponding system of linear equations
depends on the variables
$a_{1, 5}$, $a_{1, 6}$, $a_{3, 4}$, $a_{3, 5}$, $a_{3, 6}$, $a_{4, 5}$, $a_{5, 5}$, $a_{5, 6}$.\
The only non-negative integral solution is obtained for
$a_{1, 4}=a_{2, 3}=a_{5, 6}=1$, with the remaining variables vanishing.\ This
verifies the entry for $g$ in Table~\ref{Fixpointset}. $\hfill \Box$

\medskip
If $g$ has order $14$ then  $g^7$ is an involution with a K3 surface as
a fixed-point component on which $g$ acts with $3$ fixed-points, as required.

\begin{prop}[Mongardi \cite{Mon-thesis}, Proposition 6.2.16]\label{argfix11}
An admissible element
of order $11$ has type $(11,2,2)$, acting with $5$ isolated fixed-points and
eigenvalues as in Table~\ref{Fixpointset}. $\hfill \Box$
\end{prop}

This completes the proof of Theorem \ref{allowedclasses}.

In Section~\ref{geometricreal}, we will show that there are hyperk\"ahler manifolds $X$ of type $\K32$
realizing symplectic automorphisms of each of the
fifteen admissible types described in Theorem \ref{allowedclasses} and Table~\ref{Fixpointset}.

\bigskip

For K3 surfaces, it follows from \cite{Mu} that a finite $2$-group of symplectic automorphisms is 
isomorphic to a subgroup of the $2$-Sylow subgroup of $M_{23}$ and thus has order at most $2^7$.\
A short proof is given in \cite{Mason}:
starting from the fact that a symplectic
involution on a $K3$ surface has exactly $8$ isolated fixed-points, it is shown
that the centralizer of an involution in a group of symplectic
automorphisms group must be a subgroup of $\hat A_8$, the unique non-split 
extension of the alternating group $A_8$ by $\Z_2$.\
The result then follows courtesy of the fact that the $2$-Sylow subgroups of $\hat A_8$ and 
$M_{23}$ are \emph{isomorphic}.

\smallskip

Let $G$ be a $2$-group of symplectic automorphisms acting on a hyperk\"{a}hler manifold $X$ of type $\K32$.\
Let $t\in G$ be an involution, and consider the fixed-point set $X^t$.\
By Theorem \ref{allowedclasses} we know
that $X^t$ has a unique component that is a $K3$ surface, call it $Y$.\ We claim that
$C_G(t)/\langle t \rangle$ acts \emph{faithfully} on $Y$.\ 

By way of contradiction, assume that $C_G(t)/\langle t \rangle$ does \emph{not} act faithfully on $Y$.\ Then there is
$t\in A\subseteq G$ such that $|A|=4$ and $A$ acts trivially on $Y$.\ 
The group $A$ cannot be cyclic, because elements of order $4$
have discrete fixed-points by Theorem \ref{allowedclasses}.\ 
Therefore, $A=\langle s, t\rangle \cong \Z_2^2$.  Since both $s$ and $t$ leave $Y$
fixed pointwise, they each act on the tangent space at a point  $p\in Y$, 
and this consists of the tangent space of $Y$ at $p$ and the restriction
of the normal bundle.\  Both $s$ and $t$ act trivially on the first piece, and as $-1$ on the second.\ 
Therefore $st$ acts trivially on both,
and hence is the identity.\ As $s$ and $t$ are involutions then $s=t$, which is the required contradiction.

\smallskip

\begin{thm}\label{thm27}
Let $G$ be a $2$-group such that for every involution $t$, the quotient $C_G(t)/\langle t\rangle $ 
is isomorphic to a subgroup of $M_{23}$. 
Then $|G|\leq 2^7$, and $G$ belongs to one of $70$ possible isomorphism types.
\end{thm}
\pf  We verified this by checking the condition for all $2$-groups of order at most $256$~\cite{OBr91}
using Magma.  $\hfill \Box$


\section{{Isomorphism classes of admissible subgroups in ${\rm Co}_0$}}\label{SSGMT}

In this section we describe the isomorphism classes of subgroups $G\subseteq {\rm Co}_0$ with the 
following two properties:
\vspace{-2mm}
\begin{itemize}
\itemsep0em
\item[1.] $G$ consists of admissible elements,
\item[2.] $\rk \Lambda^G\geq 4$.
\end{itemize}
We call these subgroups {\it admissible.} 
We note that an admissible element generates an admissible cyclic subgroup.

The admissible subgroups $G$ are of three distinct types according to whether
$G$ is a $2$-group; all conjugacy classes of $G$ meet the
monomial group $2^{12}{:}M_{24}$ and $G$ is \emph{not} a $2$-group; or $G$ contains  elements in one of the four admissible classes 
$(3,-3,-3)$, $(6,5-3)$, $(9,3,3)$ or $(12,1,5)$ which do \emph{not\/} meet $2^{12}{:}M_{24}$.\ 
One of the main results of this section is a sufficient condition for $G$ to be isomorphic to a subgroup of 
$M_{23}$: namely, that $G$ is of the second type.\
If $G$ \emph{is\/} a $2$-group then it lies in 
$2^{12}{:}M_{24}$ just because  the monomial group contains a Sylow $2$-subgroup of ${\rm Co}_0$.\

The maximal subgroups in the first two cases classes will be explicitly described.
In the third case, $G$ has to be one of the four groups described
in Theorem~\ref{thmsgp} (i), (ii), (iii) or (v).


\subsection{Admissible subgroups related to $M_{23}$}\label{admissiblem23}

In this subsection we consider subgroups $G\subseteq {\rm Co}_0$ with the following properties:
\vspace{-2mm}
\begin{itemize} 
\itemsep0em
\item[1.] $G$ is admissible, 
\item[2.] $G$ contains \emph{no} elements of type $(3, -3, -3)$, 
\vspace{-11mm}
       \begin{equation} \label{co1props} \end{equation}
\item[3.] $G$ is \emph{not} a $2$-group.
\end{itemize}

As a main result we have:
\begin{thm}\label{thm57911}\ Assume that $G$ satisfies the hypotheses of (\ref{co1props}).\ Then $G$ is isomorphic to a subgroup of one of the following $13$ groups:
\begin{eqnarray*}
&(a)&  L_2(11),\\
&(b)&  \Z_2\times L_2(7),\\
&(c)&  \Z_2^3{:}L_2(7),\\
&(d)&  A_7,\\
&(e)&  L_3(4),\\
&(f)& (\Z_3\times A_5){:}\Z_2,\\
&(g)&  \Z_2^4{:}A_6,\\
&(h)&  \Z_2^4{:}S_5,\\
&(i)& M_{10},\\
&(j)& S_6,\\
&(k)&  \Z_3^2{:}QD_{16},\\
&(l)& \Z_2^4{:}(S_3\times S_3),\\
&(m)& Q(\Z_3^2{:}\Z_2),\ |Q|=2^6. 
\end{eqnarray*}
\end{thm}

Of the $15$ admissible classes of elements enumerated in Theorem \ref{allowedclasses} (cf.\ Table~\ref{Leechclasses}), those of types
$(3,-3,-3)$, $(6,5-3)$, $(9,3,3)$ and $(12,1,5)$ all have some power of type $(3, -3, -3)$.\ Since this latter type is excluded by assumption,
all four of these classes are excluded.\ There remain $11$ admissible conjugacy classes, these being precisely the admissible classes
that meet the monomial group $2^{12}{:}M_{24}$.\ Indeed, each of these classes meets $M_{23}$ and the class is characterized by 
the order $n$ of the elements that it contains ($n=1$, $\ldots$, $8$, $11$, $12$, $15$).

\medskip
 The upshot is that 
for $g\in G$, the character $\chi(g)$  of the representation of $G$ on $V:=\Lambda\otimes \Q$
is equal to the trace of the corresponding element in $M_{23}$ in the usual permutation representation of degree $24$.\
Explicitly, if $g\in G$ has order $n$ then 
\begin{eqnarray}\label{traceform}
\chi(g)=\varepsilon(n):=24\biggl(n\prod_{p|n}\bigl(1+{\frac{1}{p}}\bigr)\biggr)^{-1}.
\end{eqnarray}

\medskip

Following Mukai~\cite{Mu}, we call a $24$-dimensional representation $\rho$ of a finite group $H$ a {\it Mathieu representation\/}
if $\tr\rho(g)=\varepsilon(n)$ whenever $g\in H$ has order $n$.\
Consequently, if $G$ is as in the statement of Theorem \ref{thm57911} then $V$ furnishes a Mathieu representation of $G$ over $\Q$ with $\dim V^G \geq 4$.

The following general result  was obtained in~\cite{Mu}, Thm.\ (3.22):
\begin{thm}[Mukai]\label{cardG}
Suppose that $H$ has a Mathieu representation over $\Q$.\ Then
$|H|$ divides $|M_{23}|=2^7\cdot 3^2\cdot 5 \cdot 7 \cdot 11 \cdot 23$.    
\end{thm}
In our case there are no admissible elements of order $23$, whence
$|G| | 2^7\cdot3^2\cdot5\cdot 7 \cdot 11$.\ We also recover the upper bound of $2^7$ of Theorem~\ref{thm27}.

\medskip

The following formula is well-known:
\begin{eqnarray}\label{fpf}
\dim V^G = |G|^{-1}\sum_{g\in G} \chi(g).
\end{eqnarray}
We refer to this result as the \emph{invariant subspace formula} (i.s.f.).

We turn to the proof of Theorem \ref{thm57911}, which will be divided into several cases.\ Let $G$ be as in the statement of the Theorem.\ In
particular, since $G$ is admissible then $\dim V^G\geq 4$.\ We use this inequality frequently in tandem with the i.s.f.\ 
to eliminate a number of possible groups.

\paragraph{\bf Case $1$.} $11||G|$.\ We will show in this easiest case that
\begin{eqnarray*}
\mbox{$G$ is isomorphic to a subgroup of $L_2(11)$}.
\end{eqnarray*}

Let $P\subseteq G$ be a Sylow $11$-subgroup.\ Thus $P\cong \Z_{11}$, and since
there are no elements of order $11k\ (k\geq 2)$ then
$C_G(P)=P$.\ The i.s.f.\ shows that there are no dihedral subgroups of order $22$,
whence $|N_G(P)|=11$ or $55$.

In the first case, $G$ has a normal $11$-complement by Burnside's normal $p$-complement theorem, call it $Q$.\
Thus
$G = QP$ with $\gcd(|Q|, 11)=1$.\ Since $C_G(P)=P$ then there must be a prime $p=2$, $3$, $5$ or $7$
and a subgroup $E\subseteq Q$ with $E\cong \Z_p^k$ for some $k\geq 0$ such that $P$
normalizes and acts faithfully on $E$.\ Moreover, if $Q\not= 1$ then we can choose $E\not= 1$.
Setting $H=EP$, the i.s.f.\ applied to $H$ shows that the only possibility is $E=1$, so that $H=1$ and
$G=P$ is a Sylow $11$-subgroup of $L_2(11)$.

The second possibility is  $|N_G(P)|=55$.\ The previous argument shows that
$P$ cannot normalize \emph{any\/} nontrivial subgroups of $G$ of order coprime to $11$.\
So if $G$ is \emph{solvable\/} then  $G=N(P)$ has order $55$ and is the normalizer of a Sylow $11$-subgroup of $L_2(11)$.\
Finally, assume that $G$ is \emph{nonsolvable}.\ A minimal normal subgroup $N \unlhd G$ cannot have order 
coprime to $11$, so it is simple and contains $P$.\ By the Frattini argument,
$G=NN_G(P)$.\ If $N_G(P)$ is not contained in $N$ then $P$ is self-normalizing in $N$, 
so $N$ cannot be simple by the same Burnside theorem, contradiction.\ Thus $N_G(P)\subseteq N$ and $G=N$ is simple.\
Because $|G|$ divides $2^6.3^2.5.7.11$, any one of several classification theorems will tell us that
$G \cong L_2(11)$.\ This completes Case $1$.

\paragraph{\bf Case $2$.} $7||G|$.\ 
We will establish
\begin{eqnarray*}
\mbox{$G$ is isomorphic to a subgroup of $\Z_2\times L_2(7)$, $\Z_2^3{:}L_2(7)$, $A_7$ or $L_3(4)$}.
\end{eqnarray*}
Let $P\cong \Z_7$ be a Sylow $7$-subgroup.
\begin{lem}\label{lem28}
$G$ has no subgroup $H$ of any of the following types:
dihedral of order $14$, order $28$, order $2^4.7$.
\end{lem}
\begin{pf}\ The i.s.f.\ eliminates the dihedral group.\ Suppose that $|H|=28$.\ 
Since there are no elements of order $28$ and no dihedral group of order $14$,  
then $H$ is either $\Z_2\times \Z_{14}$ or $\Z_7{:}\Z_4$.\ The first 
possibility is eliminated by the i.s.f.\ The second possibility does not hold either, 
because calculation (e.g.\ by Magma) shows that $\chi$  \emph{cannot} be the character of 
a  rational representation of such a group
(the Schur index is greater than $1$).\ Suppose that $|H|=2^4.7$.\ Since there are no dihedral groups of order $14$ or
abelian groups of order $28$  then a $2$-Sylow-subgroup $Q\subseteq H$ is normal (Burnside's theorem again), 
moreover $|C_Q(P)|=2$.\ The only possibility is $Q\cong \Z_2^4$, and one verifies by applying the i.s.f.\ to $QP$ 
that this is impossible. $\hfill \Box$
\end{pf}

Now suppose that $Q\not= 1$ is a subgroup of order coprime to $7$ and normalized by $P$.\
Then $P$ leaves invariant a Sylow $r$-subgroup of $Q$ for each prime divisor $r$ of $|Q|$.\
Let $R$ be such a Sylow $r$-subgroup.\ If $r\geq 3$ then $|Q|\leq r^2$, so $P$ acts trivially
on $R$, thereby producing elements of order $7r$, contradiction.\ So $Q=R$ is a $2$-group.\ 
We have $|C_Q(P)|\leq 2$ by Lemma \ref{lem28} (there are no subgroups of order $28$).\ 
Suppose that $Q$ is extra-special.\ Then $P$ acts faithfully on $Q$, and we have $|Q|=2^{1+6k}$.\ 
The (unique) faithful irreducible representation of $Q$ has dimension $2^{3k}$, so we must have $k=1$, and 
the i.s.f.\ applied to $QP$ yields a contradiction.\ Therefore $Q$ is not extra-special.\ So if $|Q|\geq 4$, 
there is a noncyclic characteristic elementary abelian $2$-group of $Q$, call it $E$.\ Because $|C_E(P)|\leq 2$
then $P$ acts faithfully on $E$, so $|E|\geq 8$.\ If $P$ centralizes $Z(Q)$ then $EZ(Q)P$ contains a group of order $2^4.7$, contradiction.\ 
So $P$ does not centralize $Z(Q)$, whence we may, and shall,
take $E\subseteq Z(Q)$.\ The same argument then shows that $C_Q(P)=1$.
If $|Q|\geq 16$ we must then have $|Q|=2^{3k}, k\geq 2$, and the i.s.f.\ again yields a
contradiction.\
So in fact $Q=E\cong \Z_2^3$.\ Thus we have so far shown:
\begin{eqnarray*}\label{lem7signal} 
&&\mbox{Suppose there is $1\not=Q\subseteq G$ normalized by $P$ with $\gcd(|Q|, 7)=1$.}\\
&&\mbox{Then $Q\cong \Z_2$ or $\Z_2^3$, and in the second case
$C_Q(P)=1$}.
\end{eqnarray*}

Now assume that  $\Z_2^3\cong Q\unlhd G$.\ If $G$ is solvable we easily find (there being no dihedral subgroup of order $14$ by
Lemma \ref{lem28}) that
$G\cong \Z_2^3{:}\Z_7$ or $\Z_2^3{:}(\Z_7{:}\Z_3)$.\ Suppose that $G$ is nonsolvable.\ Since
$C_G(Q)$ has order coprime to $7$ we must have $C_G(Q)=Q$.\ Since 
$\Aut(Q)\cong L_3(2)\cong L_2(7)$ is a minimal simple group, the only possibility is
$G/Q \cong L_2(7)$, so that $G$ occurs in a short exact sequence
$1{\rightarrow} \Z_2^3 {\rightarrow} G {\rightarrow} L_2(7) {\rightarrow} 1$.\ 
We assert that this extension splits.\ Indeed, let $H{\subseteq}G$ be the subgroup that commutes with a given nonzero element of $Q$\
so that $H/Q{\cong}S_4$.\ Use the i.s.f.\ to see that $O_2(H){\setminus}{Q}$ contains involutions, from which it follows that $O_2(H){\cong} Q_8*Q_8$.\ Then we can check directly that $H$ is in fact a \emph{split extension}
of $S_4$ by $Q$, and in particular $Q$ splits in a Sylow $2$-subgroup of $G$.\ Therefore, the extension $G$ of $L_2(7)$ by $Q$ also splits.
Note that the possible solvable groups with $\Z_2^3\unlhd G$ are all contained  in this nonsolvable example.

Assume that $|Q|=2$.\ If $G$ is solvable, it follows from what has already been established that
$G$ is isomorphic to a subgroup of  $\Z_2\times (\Z_7{:}\Z_3)$.\ 
Finally, suppose $G$ is nonsolvable.\
Then $G/Q$ is simple of order divisible by $7$, so $G/Q\cong L_2(7), A_7, A_8, L_3(4)$.\ Since there are no elements of order $10$, only the first possibility survives.\
Therefore, $G\cong \Z_2\times L_2(7)$ or $SL_2(7)$.\ In the second case the i.s.f.\ fails,
leaving only the first possibility.\ Note once again that all solvable possibilities are contained in the nonsolvable case.

The last  possibility is that there is \emph{no} nontrivial normal subgroup of $G$ of order coprime to $7$.\ Then a minimal normal subgroup $N\unlhd G$ is simple and contains $P$.\
The only possibilities for $N$ (taking into account that $|G|$ divides $2^6.3^2.5.7$ with no elements
of order $9$) are
\begin{eqnarray*}
N \cong L_2(7),\ A_7,\ A_8,\ L_3(4).
\end{eqnarray*}
We claim that $A_8$ cannot occur.\ Indeed, $A_8\subseteq M_{23}$ acts with only $3$ orbits on the $24$ points.\ 
Therefore, applying the i.s.f.\ to our abstract group $G\cong A_8$ will certainly yield $\dim V^G=3$, contradiction.\ 
We assert that $G=N$ in the other three cases, so assume this is false.
If $N\cong L_3(4)$, we have $V = V^G\oplus U$ where $U$ is $G$-invariant 
and irreducible of dimension $20$.\ ($L_3(4)$ has no irreducible representation of dimension less than $20$ \cite{Atlas}.)\ 
Since $|G|$ is not divisible by $27$, we must have $|G/N||4$, and every nontrivial
coset of $N$ in $G$ contains an involution $t$ such that $\tr_U(t)\not= 4$ (loc.\ cit).\ 
Because $8=\chi(t)=4+\tr_U(t)$, this is a contradiction.\ On the other hand, if
$N\cong L_2(7)$ or $A_7$ then $G\cong {\rm PGL}_2(7)$ or $S_7$ respectively, and
since both of these groups contain a dihedral group of order $14$ they cannot occur either.\ This 
completes the proof that $G=N\cong L_2(7)$, $A_7$ or $L_3(4)$.\ Because 
$L_2(7)\subseteq A_7$ there is no need to include $L_2(7)$ on the list of possible groups,
and the analysis of Case 2 is done.

Before taking up Case 3 we interpolate a Lemma.
\begin{lem}\label{lem2k5} Suppose that $|G|=2^k.5$.\ Then either $G\cong \Z_5{:}\Z_4$ or $O_2(G)\cong \Z_2^4$.
\end{lem}
\begin{pf}\ Let $P$ be a Sylow $5$-subgroup of $G$.\ Since $G$ is necessarily solvable, and assuming that the first stated possibility 
does not hold we have $Q:= O_2(G)\not= 1$.\ Since there are no elements of order $10$ then $Z(Q)$ contains a 
$P$-invariant subgroup $E\cong \Z_2^4$.\  We will obtain a contradiction if $Q\not= E$.\
Indeed, in the contrary case we have $|Q|\geq 2^8$ because $C_Q(P)=1$, and this is impossible because a Sylow $2$-subgroup has order
$\leq 2^7$.\ $\hfill \Box$
\end{pf}

\paragraph{\bf Case 3.}  $5||G|$ and there exists $\Z_2^4\cong N \unlhd G$. 
We will establish
\begin{eqnarray*}
\mbox{$G$ is isomorphic to a subgroup of $\Z_2^4{:}S_5$ or $\Z_2^4{:}A_6$}.
\end{eqnarray*}

First, we can check using the i.s.f.\ that $\dim V^N=9$ and $\dim V^E=10$ for every hyperplane
$E\subseteq N$.\ This means implies that if $V = W\oplus V^N$ is a $G$-invariant decomposition,
then $\dim W=15$ and (considered as $N$-module) $W$ is the sum of the $15$ distinct nontrivial 
irreducible characters of $N$.\ Therefore, we can choose a $1$-dimensional subspace of $V^G$,
call it $V_0$ so that $U:= W\oplus V_0$ is both $G$-invariant and a \emph{free} $N$-module.\ 
There is a unique set of $1$-dimensional spaces
$V_i:= \C v_i$ ($0\leq i \leq 15$) spanning $W\oplus V_0$ and afford the irreducible characters of $N$,
and $G$ permutes the $V_i$ among themselves.\ Furthermore, we readily find that $G = N{:}H$ is a
\emph{split} extension and that $G$ acts as a transitive permutation group on $U$ with point-stabilizer $H$.\ 
Applying Lemma~\ref{lem2k5} with $G$ replaced by $O_2(G)P$ ($P$ a Sylow $5$-subgroup of $G$), we conclude that 
$N=O_2(G)$, i.e., $O_2(H)=1$.\ Furthermore,  the i.s.f.\ shows that $G$ cannot contain any abelian groups of
order $2^5$ and rank at least $4$.\ In particular, $N$ is self-centralizing in $G$
whence $H$ is isomorphic to a subgroup of \Aut$(N)\cong L_4(2)\cong A_8$.\ From the analysis of
Case 2 it also follows that $\gcd(|H|, 7)=1$, and by the i.s.f.\ there is no subgroup isomorphic to 
$\Z_2^4{:}\Z_{15}$.

From these reductions, it follows that if $G$ is solvable then $G$ is isomorphic to a subgroup of $\Z_2^4{:}(\Z_5{:}\Z_4)$, 
which is itself a subgroup of $\Z_2^4{:}S_5$.\ If $G$ is nonsolvable
then $H$ has a normal subgroup $K\cong A_5$ or $A_6$, and in the latter case we have
$H\cong A_6$ or $S_6$.\ We can eliminate the latter possibility much as in the $L_3(4)$ case handled earlier.\ 
Indeed, if $H\cong S_6$ then $V$ decomposes into a transitive permutation representation of $K$ of dimension $15$ 
(corresponding to its action on the involutions of $N$),
plus a $4$-dimensional fixed-point subspace, plus a $5$-dimensional irreducible in $V^N$.\ 
Let $t\in H$ be an involution acting on this $5$-dimensional space with trace $-1$ (such a $t$ always exists).\
 $t$ has trace $7$ or $3$ on the $15$-dimensional space, so that $\chi(t)=10$ or $6$, contradiction.\
Suppose that $K\cong A_5$.\ Because there are no subgroups of the form $\Z_2^4:\Z_{15}$ then  $G$ is isomorphic 
to a subgroup $\Z_2^4{:}S_5$.\ This completes Case 3.

\paragraph{\bf Case $4$.} $5||G|$, $\gcd(|G|, 77)=1$. We show in this case that either Case 3 holds, or 
\begin{eqnarray*}
\mbox{$G$ is isomorphic to a subgroup of $(\Z_3\times A_5){:}\Z_2$, $S_6$ or $M_{10}$}.
\end{eqnarray*}

If $O_2(G)\not= 1$ then by Lemma \ref{lem2k5} we are in Case 3 and there is nothing to  prove.\  
Hence, we may assume that $O_2(G)=1$.\ If there is a 
nontrivial normal subgroup of order prime to $5$ it must then be a $3$-group, call it $R \unlhd G$.\
Then $P\cong \Z_5$ acts trivially on $R$ since $|R|\leq 9$.\ There is no abelian subgroup of
order $45$ by the i.s.f.\, so $R=O_3(G)\cong \Z_3$.\ Let $M/R$ be a minimal normal subgroup of
$C_G(R)/R$.\ If $(\gcd(|M|, 5)=1$ then $M/R$ is a $2$-group and it centralizes $R$ and therefore
descends to a yield a nontrivial normal $2$-subgroup of $G$, contradiction.\ Thus $P\subseteq M$.\
If $M$ is solvable then $M=P\times R\cong \Z_{15}$ and $G\subseteq \Z_{15}{:}\Z_4$.\ If
$M$ is nonsolvable then $M/R\cong A_5, M\cong \Z_3\times A_5$, and $G$ is isomorphic to a subgroup
of  $(\Z_3\times A_5){:}\Z_2$.\ This contains the solvable group described earlier in the paragraph.

Finally, assume that a minimal nontrivial normal subgroup $N$ contains $P$.\ If it is equal to $P$
then $G$ is isomorphic to a subgroup of $\Z_5{:}\Z_4$.\ Otherwise, $N$ is simple, so that
$N\cong A_5$ or $A_6$ and $G=A_5$, $S_5$ or $A_6\subseteq G \subseteq\Aut(A_6)$.\ In the latter case, 
since there are no elements of order $10$ then we cannot have $G=PGL_2(9)$ or
$\Aut(A_6)$.\ Therefore, $G\cong A_6$, $M_{10}$ or $S_6$.\ (The latter two groups correspond to the
two cosets of $\Aut(A_6)/A_6$ distinct from that corresponding to $PGL_2(9)$).\ This completes the analysis of Case~4.

We have now completed the proof of Theorem~\ref{thm57911} in case  $|G|$ is divisible by one of $5$, $7$ or $11$.\ 

\paragraph{\bf Case $5$.} $|G|=2^f3^2$. We will show that 
either $G$ is contained in one of the groups occurring in  Cases 1--4, or else it is  a subgroup of one of
\begin{eqnarray*}\ \Z_3^2{:}QD_{16},\ \Z_2^4{:}(S_3\times S_3)\ \mbox{or}\ Q{:}(\Z_3^2{:}\Z_2)\
\mbox{with}\ |Q|=2^6. 
\end{eqnarray*}

First note that $G$ is solvable of $3$-length $1$.\ Because there are no elements of order $9$ then a
 Sylow $3$-subgroup $R\subseteq G$ is isomorphic to $\Z_3^2$
and $G=RT$ where $T$ is a Sylow $2$-subgroup of $G$.\ Suppose first that $O_2(G)=1$.\
Then $T$ acts faithfully on $R$ and is therefore isomorphic to a subgroup of
$\Aut(R)=GL_2(3)$.\ This latter group has Sylow $2$-subgroup $QD_{16}$ (quasidihedral), 
so $G$ is isomorphic to a subgroup of $\Z_3^2{:}QD_{16}$.

Next assume that $O_2(G)\not= 1$, $S:= O_3(G)\not= 1$.\ Since there is no
abelian subgroup of order $18$ by the i.s.f.\, we have $|Q|=3$.\
Similarly, because there is no subgroup $\Z_3\times A$ with $|A|\geq 16$ by the i.s.f.,\ we have
$|Q|=4$.\ Then $G$ is isomorphic to a subgroup of $(\Z_3\times A_4){:}\Z_2$, and this group is contained in $(\Z_3\times A_5){:}\Z_2$.

Now suppose that $O_3(G)=1$.\ Because $R$ is self-centralizing and there is no $\Z_3\times A$, $|A|=16$ as before, 
then $Q:= O_2(G)$ has order $4^k$ where $k=2$ or $3$  is the number of subgroups 
$U\subseteq R$ of order $3$ satisfying $C_Q(U)\not= 1$ (in which case $C_Q(U)\cong \Z_2^2$).\ 
Since there is no subgroup isomorphic to 
$\Z_2^5$, it follows easily that either $k=2$ and $Q\cong \Z_2^4$, or else $k= 3$
and $Z(Q)\cong \Z_2^2$.\ Moreover, $G=QH$ where $H:=N_G(R)$.

\begin{lem}\label{lem48} If $\Z_3\cong U\subseteq R$ then
either $C_H(U)=R$ or $C_Q(U)=1$.
\end{lem}
\begin{pf} Assume this is false.\ Set $Q_0:=C_Q(U)\cong \Z_2^2$ and $R=U\times U_0.$\ Then 
$U_0$ acts on $Q_0$, and since $C_G(R)=R$ then $C_{Q_0}(U_0)=1$.\ Now
$8||C_G(U)|$.\ Since $C_Q(R)=1$ we can
choose an involution $t\in C_H(U)$.\ Since $R\unlhd H$ then $t$ normalizes $R$ and we may, and shall,
choose $t$ so that it normalizes $U_0$.\ Since $t$ commutes with $U$ but not $R$ then $t$ must act as
the inverting automorphism of $U_0$.\ Now $U_0\langle t \rangle \cong S_3$ acts
on $Q_0$.\ Since  $t$ commutes with $Q_0$ then so does $U_0= [U_0, t]$, that is 
$Q_0=C_{Q_0}(U_0)$.\ This contradicts the earlier statement that $C_{Q_0}(U_0)=1$, and the proof of
the Lemma is complete. $\hfill \Box$\end{pf}

\medskip
We now consider the two possibilities for the integer $k$ defined, as before, by the equality 
$|Q|=4^k$.

\paragraph{\bf Case 5(a).} $k=2$.
Here, $Q\cong \Z_2^4$ and $H$ is isomorphic to a subgroup of
$\Z_3^2{:}QD_{16}$.\ However, because $k=2$ then a Sylow $2$-subgroup $T_0$ of $H$
\emph{cannot} act transitively (by conjugation) on the subgroups of $R$ of order $3$.\ As a consequence, $T_0$ is neither $QD_{16}$ itself, nor is it
$\Z_8$ or $Q_8$.\ It is thus a subgroup of $D_8$.\ If it is $D_8$ then we can find a $\Z_3\cong U\subseteq R$ which satisfies $C_Q(U)\not= 1, C_H(U)\not= R$, against Lemma \ref{lem48}.\
Therefore,  $T_0\cong \Z_2, \Z_4$ or $\Z_2^2.$\ Furthermore, the conditions of 
Lemma \ref{lem48} have to be satisfied.\ If $T_0\cong \Z_4$ then 
$G=\Z_2^4{:}(\Z_3^2:\Z_4)$  is a subgroup of $\Z_2^4{:}A_6$, which occurs in Case 3.\
If $T_0\cong \Z_2^2$ then $G=\Z_2^4{:}(S_3\times S_3)$, while if $T_0\cong \Z_2$,
then $G$ is isomorphic to a subgroup of  $G=\Z_2^4{:}(S_3\times S_3)$.\ This completes Case 5(a).

\paragraph{\bf Case 5(b).} $k=3$. Here, $Q$ is the product of three subgroups
$Q_i:=C_Q(U_i)\cong \Z_2^2$ ($i=0$, $1$, $2$) for three distinct subgroups $U_i\subseteq R$ or order $3$.\
We may, and shall, take $Q_0=Z(Q)$.\ Suppose that $QR$ is a \emph{proper\/} subgroup of $G$.\ Then the $2$-Sylow subgroup
$T_0$ of $H$ is nontrivial.\ $T_0$ normalizes $U_0=C_R(Z(Q))$, and no nonidentity element can
act trivially by Lemma \ref{lem48}.\ Thus $T_0=\langle t \rangle$ has order $2$.\ Now $t$ must also
normalize the unique subgroup $U_3\subseteq R$ of order $3$ satisfying $C_Q(U_3)=1$.\ We claim that $t$ also inverts $U_3$.\ 
Otherwise, $t$ centralizes $U_3$, so $U_3$ acts on
$C_Q(t)$.\ Then $C_Q(t)$ has order $16$ and contains $Z(Q)$, from which it follows that
$C_Q(t)$ is abelian.\ Then $\langle t \rangle \times C_Q(t)$ is abelian of order $32$, and
the fixed-point formula shows that this is not possible.\ This completes the proof that
$t$ inverts $U_3$.\ (It follows from these facts that the Sylow $2$-subgroup
$T=Q\langle t \rangle$ of $G$ is isomorphic to a Sylow $2$-subgroup of
$M_{23}$.) In any case,  $G\cong Q{:}(\Z_3^2{:}\Z_2)$, and Case 5(b) is finished.

\medskip
This completes the discussion of the case when $3^2||G|$. 

\paragraph{\bf Case 6.} $|G|=2^f.3$.
We choose a computational approach for the case when $3||G|$.
By Theorem~\ref{cardG} we have $|G||384$.\ 
We use the library of small groups~\cite{BEO} in Magma to find
those $G$  satisfying $3| |G| |384$ and possessing a $20$-dimensional
character $\chi-4\cdot\mathbf{1}$, $\chi(g)$ being determined by the order of $g$.\
With the exception of a single group $X_{24}$ of order $24$ (library entry $\#7$),
all others are contained in one of the $13$ groups of Theorem~\ref{thm57911}.\ 
See also Table~\ref{number3}.

\smallskip

$X_{24}$ has eight $1$-dimensional and four $2$-dimensional irreducible characters.\ 
In the representation with character $\chi-4\cdot \mathbf{1}$,
all irreducible characters appear with multiplicity $1$ apart from one $2$-dimensional
character which has  multiplicity $3$.\ However, two of the $2$-dimensional irreducible
representations have Schur index $2$, implying that the  representation
affording $\chi-4\cdot \mathbf{1}$ \emph{cannot\/} be rational.\ 
Consequently, $G$ cannot be contained in ${\rm Co}_0$.

\medskip

This completes the proof of the Theorem. $\hfill \Box$

\begin{table}\caption{Isomorphism classes of $2^f.3$-groups}\label{number3}
{
$$
\begin{array}{l|rrrrrrrr|r}
\hbox{order} & 3 & 6 & 12 & 24 & 48 & 96 & 192 & 384  & \hbox{total}\\ \hline
\hbox{no.~of groups}  &  1 & 2 & 5 & 15 & 52 & 231 & 1543 & 20169 & 22018\\
\hbox{and correct representation} &  1 & 2 & 4 & 5 & 8 & 6 & 8 & 4  & 38 \\ \hline
\hbox{contained in $M_{23}$ } &  1 & 2 & 4 & 4 & 8 & 6 & 8 & 4  & 37  \\
\hbox{others contained in ${\rm Co}_0$ } & 0& 0& 0& 0 & 0 & 0 & 0 & 0 &  0
\end{array}
$$}
\end{table}

\bigskip
There is another characterization of the $13$ isomorphism types of groups occurring in
Theorem \ref{thm57911} which we state by
reformulating the result in terms of the usual permutation action of $M_{23}$ on 
a set $\Omega$ of $24$ elements.
\begin{thm}\label{M23subgroups} The following sets coincide:
\vspace{-3mm}
\begin{itemize}
\itemsep0em
\item[(a)] The $13$ isomorphism classes of groups (a)--(m) in Theorem~\ref{thm57911}.
\item[(b)] The  isomorphism classes of subgroups of $M_{23}$ \emph{maximal} 
subject to having \emph{at least} $4$ orbits on $\Omega$.
\item[(c)] The isomorphism classes of subgroups of $M_{23}$ \emph{maximal} 
subject to having \emph{exactly} $4$ orbits.
\end{itemize}
\end{thm}

\pf In the following it is convenient to identify a group with its isomorphism class.\

That each of the $13$ groups $G$ listed in Theorem~\ref{thm57911} is isomorphic to a subgroup of
$M_{23}$ with at least four orbits on $\Omega$ is implicit in the proof of Theorem~\ref{thm57911} 
(and can be verified by a Magma calculation).\
So certainly $G$  is \emph{contained} in a group in~(b).

\medskip

On the other hand, a Sylow $2$-subgroup of $M_{23}$ is strictly contained in the group $Q(\Z_3^2{:}\Z_2)$, 
so \emph{no}  $2$-subgroup $G\subseteq M_{23}$ can be maximal subject to having at least $4$ orbits.\
It is then clear that every group $G$ in (b) satisfies the assumptions, 
and therefore also the conclusions of
Theorem \ref{thm57911}, so $G$ is contained in one of the groups of part~(a).\ Together with the previous paragraph, 
and because there are no containments between the groups (a)--(m) in Theorem~\ref{thm57911}, 
this shows that the sets in parts (a) and (b) coincide. 

\medskip
Finally,  each $G$ in (a)  has \emph{exactly} $4$ orbits, 
for otherwise it would have at least $5$ orbits and would appear in Mukai's list~\cite{Mu} --- but it does not.\ Therefore,
by the same argument as the last paragraph, 
the groups in (a) and (c) also coincide.\ This completes the proof of the Theorem.
$\hfill \Box$

\smallskip
In the following, we let $\mathcal{T}$ denote the isomorphism classes of groups satisfying any (and hence all)
of the properties (a)--(c) in the last theorem.

\begin{rem}\rm \setlength{\baselineskip}{4.7mm}
The $M_{23}$-conjugacy classes of the groups $G$ in $\mathcal{T}$
are \emph{not} unique.\
For the groups $A_7$, $\Z_2^3{:}L_2(7)$ and $\Z_2^4{:}S_5$ there are two conjugacy classes.\
For $\Z_2^4{:}S_5$, the orbit structure on $\Omega$ is different for the two conjugacy classes, whence
they are also not conjugate in $M_{24}$.
\end{rem}

\smallskip

\begin{rem}\rm \setlength{\baselineskip}{4.7mm}
The proof of Theorem~\ref{thm57911} together with Theorem~\ref{M23subgroups} 
actually shows that if a finite group $G$ has a Mathieu representation $V$ over $\Q$
with $\dim V^G \geq 4$ and is not a $2$-group then $G$ is contained in a group
from $\mathcal{T}$.
\end{rem}

\smallskip

\begin{rem}\rm \setlength{\baselineskip}{4.7mm}
The subgroups of $M_{23}$ with at most four orbits also appear in the work of
Dolgachev and Keum~\cite{DoKe} on finite groups $G$ of symplectic automorphism of algebraic K3 surfaces
$Y$ defined over an algebraic closed field $k$ of positive characteristic $p$.\
Dolgachev and Keum show that if ${p}\not|\ |G|$ then $G$  has elements of order no bigger than $8$,
and affords a Mathieu representation on $V=H^*_{\rm et}(Y,\Q_{\ell})$ for $\ell \not = p$ prime, 
with $\dim V^G \geq 4$.

Using further geometric arguments, they found ten isomorphism classes of such groups which could potentially act on
a K3 surface and be maximal, i.e.~not contained in a larger such group (\cite{DoKe}, Thm.~5.2).\
These are precisely the ten groups of our Theorem~\ref{thm57911} containing no elements of order greater than $8$
(a condition which excludes $L_2(11)$, $\Z_2\times L_2(7)$ and $(\Z_3\times A_5){:}\Z_2$).\
Our proof of Theorem~\ref{thm57911} 
assumes only that $G$ has a rational Mathieu representation with $\dim V^G \geq 4$ and is not a $2$-group.\ 
On the other hand, there are nine types of $2$-groups \emph{not} contained in a $2$-Sylow subgroup of $M_{23}$,
as we will show in the next Subsection.

Dolgachev and Keum also list all isomorphism types of groups which can occur in their situation
and which do \emph{not} appear  in Mukai's classification.\ These are necessarily subgroups of the ten maximal types.\
They list $28$ such groups (\cite{DoKe}, Thm.~5.2).\ These are in essential agreement  with our
list of such groups, which can be read off from Table~\ref{Gconclasses}.\ (The only ambiguity is that
our nonisomorphic groups No.~144 and No.~146 of order $192$ are both of type $\Gamma_{13}a_1{:}3$ and hence 
both correspond to the group (xxxiii) from~\cite{DoKe}.)

If  $p||G|$, further groups $G$ can be realized
(cf.\ Section~6 in~\cite{DoKe} and~\cite{Ko2}).\ Some of these,  such as $L_2(11)$, $M_{11}$ and $M_{22}$,
are subgroups of $M_{23}$.\ Others, such as $U_4(3)$, are not.\ It follows from the calculations described 
in Section~\ref{admissibleconj} that $M_{11}$, $M_{22}$ and $U_4(3)$ cannot be groups of
symplectic automorphisms of any hyperk\"ahler manifold of type $K3^{[n]}$, and likewise cannot be groups 
of symplectic autoequivalences for a complex K3 surface as considered in~\cite{Huy-conway}.
\end{rem}



\subsection{Admissible $2$-groups}

Suppose $G$ only has elements of orders $1$, $2$, $4$ or $8$. Then $G$ is necessarily a $2$-group, and
we already know from Theorem~\ref{cardG} that $|G|\leq 2^7$.
To describe the possible $2$-groups we will describe a computational approach, which seem unavoidable if one wants a complete 
proof of the result that is not too long.\

\begin{prop}\label{new2groups} Let $G$ be a $2$-group of exponent at most $8$ 
having a complex Mathieu representation $V$ with $\dim V^G \geq 4$.
Then $G$ is isomorphic either to a subgroup
of $M_{23}$, or to one of $9$ additional
groups of order $16$, $32$ or $64$ described in Table~\ref{new2}.
\end{prop}

\begin{table}\caption{The nine non-excluded $2$-groups} \label{new2}
$$\begin{array}{rrcccc}
\hbox{No.} & \hbox{order} & \hbox{Group library} & \hbox{Symbol} & A & G/A \\  \hline 
   1 & 16 & \#4 & \Gamma_2c_2  & \Z_4\times \Z_2  &  \Z_2\\
   2 & 16 & \#5 &  \Z_8\times \Z_2   & \Z_8\times \Z_2  & 1\\ 
   3 & 16 & \#10 &  \Z_4\times \Z_2^2   & \Z_4\times \Z_2^2  & 1\\ 
  4 & 32 &  \#8 &  \Gamma_7a_3  & \Z_4\times \Z_2  &  \Z_4 \\ 
  5 & 32 &  \#30 &  \Gamma_4c_1  & \Z_4\times \Z_2^2& \Z_2\\ 
  6 & 32 &  \#32 &  \Gamma_4c_3  &\Z_4^2 & \Z_2\\ 
  7 & 32 &  \#35 &  \Gamma_4a_3 &\Z_4^2 & \Z_2\\ 
  8 & 32 &  \#50 & \Gamma_5a_2  &\Z_4\times \Z_2   &  \Z_2^2\\ 
  9 & 64 &  \#36 & \Gamma_{23}a_3  &\Z_4^2&  \Z_4
\end{array}$$
\end{table}

\pf We first use the library of small groups in Magma to find
the $2$-groups $G$ of order at most $256$ and exponent at most $8$~\cite{OBr91}.\ Then we check to 
see if $G$ has a $20$-dimensional representation of the correct type.\footnote{Magma Commands: \tt D:=SmallGroupDatabase();
value:=[20,4,99,0,99,99,99,-2];
[ \#[ G : G in [SmallGroup(o, n): n in [1..NumberOfSmallGroups(o)]] $\mid$
 (\{x[1]: x in Classes(G)\} subset \{1,2,4,8\}) and IsCharacter(CharacterRing(G) ! [value[x[1]] : x in Classes(G)]) ] : 
 o in [1,2,4,8,16,32,64,128,256] ];}

\begin{table}\caption{Isomorphism classes of $2$-groups}\label{number2}
{
$$
\begin{array}{l|rrrrrrrrr|r}
\hbox{order} & 1 & 2 & 4 & 8 & 16 & 32 & 64 & 128 & 256 & \hbox{total}\\ \hline
\hbox{no.~of $2$-groups}  & 1 & 1 & 2 & 5 & 14 & 51 & 267 & 2328 & 56092 & 58761\\
\hbox{of exponent $\leq8$}  & 1 & 1 & 2 & 5 & 13 & 45 & 234 & 2093 & 53529 & 55923\\
\hbox{and correct representation} & 1 & 1 & 2 & 5 & 12 & 12 & 6 & 1 & 0  & 40 \\ \hline
\hbox{contained in $M_{23}$ } &  1 & 1 & 2 & 5 & 9 & 7 & 5 & 1 & 0 & 31\\
\hbox{others contained in ${\rm Co}_0$ } & 0& 0& 0& 0 & 1 & 2 & 1 & 0 & 0 & 4 
\end{array}
$$}
\end{table}

The result is given in Table~\ref{number2}.\
As expected, there are \emph{no} such groups of order $256$,
and a \emph{unique} group of order $128$, namely a Sylow $2$-subgroup of $M_{23}$.\
We also list in Table~\ref{number2} the number of $2$-groups contained in $M_{23}$.

\medskip
To explain Table~\ref{new2}, recall that a normal subgroup $A\unlhd G$ which is maximal with respect to being abelian, 
is necessarily self-centralizing.\ 
Thus $G/A$ is isomorphic to a subgroup of $\Aut(A)$.\
In the present situation, the existence of the complex representation means that $A$ is necessarily isomorphic
to a subgroup of one of $\Z_2^4$, $\Z_4\times\Z_2^2$, $\Z_4^2$ or $\Z_8\times\Z_2$.\ 
Table~\ref{new2} lists a choice of $A$ for each of the nine isomorphism classes of $G$ \emph{not\/} contained in $M_{23}$.\ 
Some choices of $A$ (elementary abelian and cyclic) are not represented in Table~\ref{new2}, meaning that the corresponding $G$ 
is contained in $M_{23}$. The group library number refers to~\cite{BEO}, the symbol for the non-abelian groups is as
in~\cite{HS}.  \qed

\smallskip

Four of the nine groups in Table~\ref{new2} satisfy $\dim V^G \geq 5$ and have been
described in~\cite{Mu}, Prop.~(6.3.).

\begin{lem}\label{aclasses} 
The number of conjugacy classes of abelian subgroups $A\cong  \Z_4\times \Z_2 $,
$\Z_8\times \Z_2$, $  \Z_4^2$ and  $\Z_4\times \Z_2^2$ in $2^{12}{:}M_{24}$ containing only
admissible elements is $8$, $0$, $6$ and~$7$, respectively. 
\end{lem}
\pf We construct all abelian subgroups of rank at most~$3$ by the following procedure.\ We first select representatives $a$ for
each of the $7$ conjugacy classes in $2^{12}{:}M_{24}$ of the correct Conway types.\ For each of those elements, we select 
representatives $b$ for the conjugacy classes of its centralizer  in $2^{12}{:}M_{24}$ and check if
the subgroups $H=\langle a, b\rangle$  generated by $a$ and $b$ contains only elements of admissible Conway types.\ In total,
there are $26$ such conjugacy classes of subgroups $H$.\
Then we select representatives $c$ for the conjugacy classes of the centralizer of
the subgroup $H$ and determine the conjugacy classes of subgroups $K=\langle H,c\rangle$ in $2^{12}{:}M_{24}$.\
There are $41$ such conjugacy classes of subgroups $K$ containing only admissible elements.\ 
The resulting number of conjugacy classes for each choice of $A$
is as stated in the Lemma.
\phantom{xxxx}\qed
 
\medskip

\begin{thm}\label{thm2group}
Suppose $G$ is an admissible $2$-group.
Then either $G$ is isomorphic to a subgroup of $M_{23}$,
or is isomorphic to $\Z_4\times\Z_2^2$, or is contained
in a group of order $32$ or $64$ described  below (Table~\ref{new2}, nos.~8 and~9).
\end{thm}

\noindent
{\bf Proof:} We have to check which of the nine groups $G$ of Proposition~\ref{new2groups} (cf.\ Table~\ref{new2}) are
isomorphic to subgroups of ${\rm Co}_0$ and contain only admissible elements.\ Note that  because the monomial subgroup $2^{12}{:}M_{24}$
contains a Sylow $2$-subgroup of ${\rm Co}_0$, the relevant calculations can be carried out in the monomial group.

\medskip

For each conjugacy class of abelian subgroups $A\subseteq2^{12}{:}M_{24}$ as in Lemma \ref{aclasses},
we compute the normalizer $N(A)$ in $2^{12}{:}M_{24}$ and check if it contains a group $G$ as in Table~\ref{new2} 
containing only admissible elements.\
This is done by selecting a representative $d$ of each conjugacy class of elements in $N(A)$
such that $d^4$ is in $A$ and determining the isomorphism type of $\langle A,d\rangle$. 
In the case of group no.~8 in Table~\ref{new2} we also select in addition all elements $e$ of $N(A)$ 
such that $e^2$ and the commutator $[e,d]$ is in $A$ 
and determine the isomorphism type of $\langle A,d,e\rangle$.\ It transpires that
only groups no.~3, 4, 8, and 9 occur.\ Moreover, group no.~4 is a subgroup of group no.~9. 
\phantom{xxxx}\qed


\section{Conjugacy classes of admissible subgroups in ${\rm Co}_0$}\label{admissibleconj}

In the previous section, we determined by largely theoretical arguments the 
abstract isomorphism types of admissible subgroups of ${\rm Co}_0$.

For the classification of group lattices $(\Lambda_G,G)$, we will enumerate the 
conjugacy classes of admissible subgroups $G\subseteq {\rm Co}_0$. The following main result of this
section depends heavily on computer calculations:
\begin{thm}\label{admissiblegroupclasses}
There are $198$ conjugacy classes of admissible subgroups $G\subseteq {\rm Co}_0$.\ In particular, there are exactly
$22$ classes which are \emph{maximal} (with respect to containment), which are described as follows:
\vspace{-3mm}
\begin{itemize}
\itemsep0em
\item[(a)] thirteen subgroups of $M_{23}$ (Table~\ref{M23groups});
\item[(b)] two groups $3^4.A_6$ and $3^{1+4}{:}2.2^2$ related to ${\cal S}$-lattices;
\item[(c)] two groups of order $48$ and five $2$-groups (Table~\ref{M23groups2}).
\end{itemize}
\end{thm}
Detailed information about all of these classes of groups can be found in Table~\ref{Gconclasses} 
in the appendix.\ The corresponding group lattices $(\Lambda_G,G)$, which we will discuss in the 
next section, are pairwise nonisomorphic (cf.~Table~\ref{Lisoclasses} in the appendix).

\noindent
{\bf Remark:} There are just $82$ admissible conjugacy classes of $G$ containing
only elements of $K3$-type. The corresponding group lattices $(L_G,G)$ were
first determined by Hashimoto~\cite{Ha}, who used the 23 Niemeier lattices with roots rather than the
Leech lattice that we use here.\ These $82$  classes 
are in $1$-$1$ correspondence with the combinatorial structure of symplectic
group actions on a K3 surface, as determined by Xiao~\cite{xiao}.

\medskip

In the following, we describe our method of computing admissible conjugacy classes of subgroups $G\subseteq {\rm Co}_0$ using Magma.\
We first realized ${\rm Co}_0$ as a group of integral $24\times 24$ matrices starting from an explicit 
description of the Leech lattice.\ This realization was used to determine the conjugacy classes of ${\rm Co}_0$
and to compute the dimension of $\Lambda^G$ for a subgroup $G\subseteq  {\rm Co}_0$.\ In addition, 
a realization as permutation group on the $196560$ minimal vectors of $\Lambda$ together with an
explicit isomorphism with the matrix group realization was constructed.\ This realization was used to 
check if two subgroups of ${\rm Co}_0$ are conjugate.\
Computations with both realizations are relatively time consuming and had to be minimized.

\smallskip

By Theorem~\ref{thmsgp}, 
an admissible subgroup $G$ is either conjugate to a subgroup of the monomial 
group $2^{12}{:}M_{24}$ or a subgroup of one of three groups related to ${\cal S}$-lattices.\ We identified $2^{12}{:}M_{24}$ inside the permutation representation of ${\rm Co}_0$ by computing the
stabilizer of $1104$ norm~$4$ vectors $(\pm u \pm v)/2$,  where $\pm u$ and $\pm v$  run through a coordinate frame of $\Lambda$.\
For the ${\cal S}$-lattices, we determined the pointwise stabilizers by using the explicit realizations given
by Curtis~\cite{Curtis}.\ For the ${\cal S}$-lattice of rank $6$, we computed in addition the normalizer of the 
stabilizer in ${\rm Co}_0$.\ For calculations inside $2^{12}{:}M_{24}$, we used a
realization as a permutation group on the $48$ elements.\
We also determined an explicit isomorphism between this permutation realization of $2^{12}{:}M_{24}$
and the above described frame stabilizer inside ${\rm Co}_0$.\
For calculations inside the three  ${\cal S}$-lattice 
groups, we constructed a permutation representation
of much smaller degree together with an
explicit isomorphism with the corresponding subgroups in ${\rm Co}_0$.\
This allowed us to determine their complete subgroup lattice.
For each conjugacy class of subgroups we constructed
the corresponding subgroup of ${\rm Co}_0$ and selected the admissible one.\ 
Finally, we determined the ${\rm Co}_0$-conjugacy classes.

\medskip

For the subgroups $G\subseteq 2^{12}{:}M_{24}$ we distinguished two cases: $2$-groups and 
non $2$-groups.\ For the  non $2$-groups, we were able to construct the corresponding part of the
subgroup lattice of $2^{12}{:}M_{24}$ completely.\ Starting from $2^{12}{:}M_{24}$, we constructed inductively via the order
all maximal subgroups of fixed order of the already found conjugacy classes.\ Then we checked for a given order all 
the found maximal subgroups for conjugacy in $2^{12}{:}M_{24}$.\ Since the number of required conjugacy checks 
is growing almost quadratically with the number of such subgroups, we first determined for each subgroup
the size of each of its conjugacy classes, using this as a numerical indicator for non-conjugacy and thereby reducing
 the number of conjugacy checks in $2^{12}{:}M_{24}$.\ Overall, there are $279,343$ classes of
such subgroups in $2^{12}{:}M_{24}$.\ Just $280$ of these classes contain only admissible elements and
of these, $241$ classes are admissible.\ They belong to $94$ ${\rm Co}_0$-conjugacy classes.

It follows from Theorem~\ref{M23subgroups} that 
the admissible subgroups of $2^{12}{:}M_{24}$ of order divisible by
$3$, $5$, $7$, or $11$ are isomorphic to subgroups of the set 
$\mathcal{T}$ of isomorphism classes of subgroups of $M_{23}$ maximal
among having at least four orbits on~$\Omega$.\
In most cases, there exist several $2^{12}{:}M_{24}$ conjugacy classes for a group
in  $\mathcal{T}$.
With respect to ${\rm Co}_0$-conjugacy however, our calculations show:
\begin{thm}\label{thm57911conj}
For each of the groups in $\mathcal{T}$, there exists a unique 
${\rm Co}_0$-conjugacy class inside $2^{12}{:}M_{24}$. $\hfill \Box$
\end{thm}

However, for groups of order $2^f.3$ there are further possibilities:
\begin{thm}
Suppose that $G\subseteq 2^{12}{:}M_{24}$  is admissible and \emph{not} a $2$-group.\
Then  $G$ is conjugate in ${\rm Co}_0$ to a subgroup of either a group in the set ${\cal T}$ of subgroups of 
$M_{23}$ (cf.\ Theorem~\ref{M23subgroups}),
or  one of two conjugacy classes of groups of order $48$. $\hfill \Box$
\end{thm}
Although the two groups of order $48$ project isomorphically to subgroups in $M_{23}\subseteq M_{24}$, they are
not conjugate.

\smallskip

We collect basic information about the groups in $\mathcal{T}$ in Table~\ref{M23groups} and
the two groups of order~$48$ in Table~\ref{M23groups2}.
The entries of the first four columns is self-explanatory, the column $A_{\Lambda^G}$ gives the 
structure of the discriminant group $(\Lambda^G)^*/\Lambda^G$.\
The last three columns give information
on all $2^{12}{:}M_{24}$-conjugacy classes of $G$.\ To explain this, let $E=\Z_2^{12}$ and
$M=M_{24}$, so that the monomial group is $E{:}M$.\  Then in Table~\ref{M23groups} and~\ref{M23groups2},
$A=G\cap E$ and $P=G/A \subseteq M$.\
The seventh column describes the orbits of $P$ on $\Omega$, and a star in the last column indicates that 
$G\cong A.P$ is {\it not\/} a subgroup $A{:}P$ of $E{:}M$.

\begin{table}\caption{Maximal admissible groups contained in $M_{23}$}\label{M23groups}

$$
\begin{array}{rlrcrcrrc}
&  & &  & &   \multicolumn{3}{c}{2^{12}{:}M_{24}\hbox{-classes}}   \\ \cline{6-8}
\hbox{No.} & \  G &|G| &  \!\!\!  \mbox{rk} \Lambda^G  \! \!\!\!& A_{\Lambda^G}  
&(|P|,|A|) & \mbox{orbits}  \\
\noalign{\hrule height0.8pt}

1   & {\rm L}_2(11)
 &  660
 &  4
 &  11^2
 &
[ 660, 1 ] & 
[ 1, 1, 11, 11 ]
\\  \noalign{\hrule height0.8pt}

2   & {\rm L}_3(4)
 &  20160
 &  4
 &   2^1 42^1 
 &
[ 20160, 1 ] & 
[ 1, 1, 1, 21 ]
\\ \hline

3   & A_7
 &  2520
 &  4
 &  105^1
 &
[ 2520, 1 ] & 
[ 1, 1, 7, 15 ]
\\ \hline

4   & \Z_2^3{:}{\rm L}_2(7)
 &  1344
 &  4
 &  4^1 28^1
 &
[ 1344, 1 ] & 
[ 1, 7, 8, 8 ]
\\ & & & & &  
[ 1344, 1 ] & 
[ 1, 1, 8, 14 ]
\\ & & & & & 
[ 168, 8 ] & 
[ 1, 1, 1, 7, 14 ]
\\ \hline

5   & \Z_2\times {\rm L}_2(7)
 &  336
 &  4
 &  14^2
 &
[ 336, 1 ] & 
[ 1, 2, 7, 14 ]
\\ & & & & &  
[ 168, 2 ] & 
[ 1, 1, 7, 7, 8 ]
\\  \noalign{\hrule height0.8pt}

6   & \Z_2^4{:}A_6
 &  5760
 &  4
 &  4^1 24^1 
 &
[ 5760, 1 ] & 
[ 1, 1, 6, 16 ]
\\ & & & & & 
[ 360, 16 ] & 
[ 1, 1, 1, 6, 15 ]
\\ \hline

7   & \Z_2^4{:}S_5
 &  1920
 &  4
 &  4^1 40^1
 &
[ 1920, 1 ] & 
[ 1, 1, 2, 20 ]
\\ & & & & & 
[ 1920, 1 ] & 
[ 1, 2, 5, 16 ]
\\ & & & & &  
[ 120, 16 ] & 
[ 1, 1, 2, 5, 15 ]
\\ \hline

8   & S_6
 &  720
 &  4
 &  6^1 30^1
 &
[ 720, 1 ] & 
[ 2, 6, 6, 10 ]
\\ & & & & &  
[ 720, 1 ] & 
[ 1, 2, 6, 15 ]
\\ \hline

9   &  M_{10}
 &  720
 &  4
 &  2^1 60^1
 &
[ 720, 1 ] & 
[ 1, 1, 10, 12 ]
\\ \hline

10   & (\Z_3\times A_5){:}\Z_2
 &  360
 &  4
 &  15^2
 &
[ 360, 1 ] & 
[ 1, 3, 5, 15 ]
\\ \noalign{\hrule height0.8pt}

11  &  Q(\Z_3^2{:}\Z_2)
 &  1152
 &  4
 &  8^1 24^1
 & 
[ 1152, 1 ] &
[ 1, 3, 4, 16 ]
\\ &&&&&  
[ 288, 4 ] &
[ 1, 3, 4, 4, 12 ]
& *
\\ &&&&&
[ 72, 16 ] &
[ 1, 1, 3, 3, 4, 12 ]
\\ \hline

12  &  \Z_2^4{:}(S_3\times S_3)
 &  576
 &  4
 &  12^1 24^1
 & 
[ 576, 1 ] &
[ 2, 3, 3, 16 ]
\\ &&&&&
[ 576, 1 ] &
[ 1, 3, 8, 12 ]
\\ &&&&&
[ 36, 16 ] &
[ 1, 2, 3, 3, 6, 9 ]
\\ \hline

13  &  \Z_3^2{:}QD_{16}
 &  144
 &  4
 &  6^1 36^1
 & 
[ 144, 1 ] &
[ 1, 2, 9, 12 ]
\\ \hline
\noalign{\hrule height0.8pt}

\end{array}
$$
\end{table}

\begin{table}\caption{Maximal admissible groups contained in $2^{12}{:}M_{24}$ but not in $M_{23}$}\label{M23groups2}

$$
\begin{array}{rlrcrcrcc}
&  & &  & &   \multicolumn{3}{c}{2^{12}{:}M_{24}\hbox{-classes}}   \\ \cline{6-8}
\hbox{No.} & \ \  G &|G| &  \!\!\!  \mbox{rk} \Lambda^G  \! \!\!\!& A_{\Lambda^G} \ 
&(|P|,|A|)& \mbox{orbits}  \\
\noalign{\hrule height0.8pt}

1  & \#49
 &  48
 &  5
 &  2^3 6^1 12^1
 &
[ 48, 1 ] &
[ 2, 2, 6, 6, 8 ] & *
\\ &&&&&
[ 24, 2 ] &
[ 1, 1, 3, 3, 8, 8 ] & *
\\ &&&&&
[ 24, 2 ] &
[ 1, 1, 2, 6, 6, 8 ] & *
\\ &&&&&
[ 6, 8 ] &
[ 1, 1, 2, 2, 3, 3, 6, 6 ] & *
\\ \hline
2  & \#32
 &  48
 &  4
 &  2^2 4^1 12^1
 & 
[ 48, 1 ] &
[ 2, 6, 8, 8 ] & *
\\ &&&&&
[ 24, 2 ] &
[ 1, 1, 6, 8, 8 ] & *
\\ &&&&&
[ 24, 2 ] &
[ 1, 1, 2, 6, 6, 8 ] & *
\\ &&&&&
[ 24, 2 ] &
[ 1, 1, 2, 6, 6, 8 ] & *
\\  \noalign{\hrule height0.8pt}
3  & \Z_4 \times \Z_2^2   
 &  16
 &  7
 &  2^4 4^3
 & 
[ 16, 1 ] &
[ 2, 2, 2, 2, 4, 4, 8 ] & *
\\ &&&&&
[ 8, 2 ] &
[ 1, 1, 1, 1, 2, 2, 8, 8 ] & *
\\ &&&&&
[ 8, 2 ] &
[ 1, 1, 2, 2, 2, 4, 4, 8 ] & *
\\ &&&&&
[ 8, 2 ] &
[ 1, 1, 2, 2, 2, 2, 2, 2, 2, 8 ] & *
\\ &&&&&
[ 8, 2 ] &
[ 1, 1, 2, 2, 2, 2, 2, 2, 2, 8 ] & *
\\ &&&&&
[ 4, 4 ] &
[ 1, 1, 1, 1, 2, 2, 2, 2, 2, 2, 4, 4 ] & *
\\ &&&&&
[ 2, 8 ] &
[ 1, 1, 1, 1, 1, 1, 1, 1, 2, 2, 2, 2, 2, 2, 2, 2 ] & *
\\ \hline
4  & \Z_4^2
 &  16
 &  6
 &  2^2  4^4
 &
[ 16, 1 ] &
[ 4, 4, 4, 4, 4, 4 ] & *
\\ &&&&&
[ 8, 2 ] &
[ 1, 1, 2, 2, 2, 4, 4, 8 ] & *
\\ &&&&&
[ 4, 4 ] &
[ 1, 1, 1, 1, 2, 2, 2, 2, 2, 2, 4, 4 ] & *
\\ \hline
5  &  \Gamma_5a_2
 &  32
 &  6
 &  2^2 4^4
 &
[ 8, 4 ] &
[ 1, 1, 2, 2, 2, 2, 2, 2, 2, 8 ] & *
\\ \hline
6  & \Gamma_7a_3
 &  32
 &  4
 &  2^2 8^2
 &
[ 8, 4 ] &
[ 1, 1, 2, 2, 2, 4, 4, 8 ] & *
\\ &&&&&
[ 8, 4 ] &
[ 1, 1, 2, 2, 2, 4, 4, 8 ] & *
\\ \hline
7  & \Gamma_{23}a_3
 &  64
 &  4
 &  2^3 16^1
 & 
[ 16, 4 ] &
[ 1, 1, 1, 1, 4, 8, 8 ] & *

\\ \hline

\end{array}
$$
\end{table}


\medskip

For the $2$-subgroups of  $2^{12}{:}M_{24}$, we were unable to compute the corresponding part of the subgroup
lattice because the number conjugacy classes became too large.\ Instead we determined only the admissible one.
First we selected a $2$-Sylow subgroup $P\subseteq 2^{12}{:}M_{24}$ and hence ${\rm Co}_0$.\
All $P$-conjugacy classes of admissible $2$-groups $G\subseteq P$ 
were determined, 
starting with the trivial group and successively adding further elements.
Finally, we tested for $2^{12}{:}M_{24}$-conjugacy and ${\rm Co}_0$-conjugacy. 
The numbers of such conjugacy classes of groups for a given order are listed in Table~\ref{number2conj}. 

\begin{table}\caption{Conjugacy classes of admissible $2$-groups}\label{number2conj}
{
$$
\begin{array}{l|rrrrrrrrr|r}
\hbox{order} & 1  & 2 & 4 & 8 & 16 & 32 & 64 & 128 & 256  & \hbox{total}\\ \hline
\hbox{\# $P$-classes}              & 1 & 27 & 317 & 1312 & 2190 & 803 & 239 & 30 & 0 & 4919  \\
\hbox{\# $2^{12}{:}M_{24}$-classes}& 1 &  2 &   7 &   27 &   63 &  34 &  20 &  3 & 0 & 157 \\ 
\hbox{\# ${\rm Co}_0$-classes}     & 1 &  1 &   2 &    6 &   15 &  10 &   6 &  1 & 0 & 42  
\end{array}
$$}
\end{table}

In particular, we obtained the following result:
\begin{thm}\label{thm2groupconj}
Let $G\subseteq {\rm Co}_0$ be an admissible $2$-group.\
Then $G$ is conjugate to a subgroup of either a group in the set ${\cal T}$ of subgroups of 
$M_{23}$ as described in Theorem~\ref{M23subgroups}, or  
of one of $5$ conjugacy classes of $2$-groups listed in Table~\ref{M23groups2}. $\hfill \Box$
\end{thm}

\noindent
{\bf Remarks}: Theorem~\ref{thm2groupconj} implies Theorem~\ref{thm2group}.\ Since the computations for Theorem~\ref{thm2groupconj}
take much longer and both computations are independent, we have presented both of them.

Although the $2$-group $\Z_4^2$  
is one of the $5$ maximal groups in Theorem~\ref{thm2groupconj}, it
is also isomorphic to a subgroup of $M_{23}$.

\medskip

As a final step, we took all the admissible
conjugacy classes of ${\rm Co}_0$ obtained from the ${\cal S}$-lattices and from $2^{12}{:}M_{24}$,
rechecked for ${\rm Co}_0$-conjugacy, and determined the corresponding subgroup lattice structure.\
The result is Theorem~\ref{admissiblegroupclasses}, and the information given in Table~\ref{Gconclasses}
in the appendix.


\section{Subgroups of $O(L)$}\label{subol}

In this section, we investigate which conjugacy classes of subgroups $G$ and lattices $\Lambda_G$ found
in the previous section can indeed be realized as symmetries of the lattice
$$L=H^2(X,\Z)\cong E_8(-1)^2\oplus U^3\oplus \langle -2 \rangle.$$

For most of the realizable groups $G$, we will also determine the exact number of
corresponding isomorphism classes of embeddings $(L_G,G)\rightarrow (L,O(L))$. 

\medskip

For each of the $198$ ${\rm Co}_0$-classes of admissible groups $G$, we first determine the isomorphism classes 
of coinvariant lattices  $L_G\cong \Lambda_G(-1)$.\ This is easily done with Magma.\ There are $69$ such lattices, listed in
Table~\ref{Lisoclasses} in the appendix.

The lattices $K$ in Table~\ref{Lisoclasses} are naturally divided into three main types depending on the
relation between ${\rm rk}\, K$ and ${\rm rk}\,A_K$.\ Setting
\[\alpha(K):=24-{\rm rk}\, K - {\rm rk}\, A_K, \] 
we find the following possibilities:
\vspace{-2mm}
\begin{itemize}
\item[(a)] $54$ lattices with $\alpha(K)\ge 2$:\ these are the $13$ lattices corresponding to the maximal subgroups of
$M_{23}$ with at least four orbits, and the $41$ lattices coming from symplectic group actions on K3 surfaces as classified in \cite{Ha}.\
\item[(b)] $4$ lattices with $\alpha(K)=1$: two of them correspond to the two maximal ${\cal S}$-lattice groups, while
the other two  correspond to certain subgroups.
\item[(c)] $11$ lattices with $\alpha(K)=0$:\ they correspond to $2$-groups or groups of order $2^f.3$ in $2^{12}{:}M_{24}$ which are 
\emph{not\/} conjugate to a subgroup of $M_{23}$.
\end{itemize}
Apart from the non-canonically defined $2$-adic genus symbol, the table in Section~10.2 of~\cite{Ha} 
seems to be in complete agreement with our Tables~\ref{Gconclasses} and \ref{Lisoclasses}.

For certain lattices $K$, there is more than one group such that $K\cong L_G$.\ For such $K$ we
list all classes of groups $G$ with $K \cong L_G$ in Table~\ref{Lgroups}.\
For the $K3$ cases, this again is in accord with the Table in Section~10.4 of \cite{Ha}.\
Let ${\cal G}(K)$ be the classes of groups $G$ from Table~\ref{Gconclasses} which have the same $K \cong L_G$.
By inspection we find:
\begin{thm}\label{GsetofK}
The set ${\cal G}(K)$ contains a \emph{unique\/} maximal class $G_{\rm max}(K)$.\
The other classes $H\in {\cal G}(K)$ correspond to subgroups of $G_{\rm max}(K)$.\
More precisely, the classes in ${\cal G}(K)$ can be identified with the $O(K)$-conjugacy classes of subgroups $H$
of $O_0(K)$ which have a trivial fixed-point lattice $K^H$ and are contained in a conjugacy class $G_{\rm max}(K)$.\
Furthermore, for the lattices of type (a) and (b) one has $G_{\rm max}(K)=O_0(K)$. \qed
\end{thm}

\medskip

To study embeddings $(L_G,G)\rightarrow (L,O(L))$, 
we have to find a lattice $L^G$ of rank $23-{\rm rk}\, L_G$ and an  extension of $L_G\oplus L^G$ to $L$.\ 
Such an extension is described by a \emph{glue code\/} $C\subseteq A_{L_G}\oplus A_{L^G}$ which has
to be an isotropic subspace with respect to the quadratic form $q_{L_G}\oplus q_{L^G}$ for $A_{L_G}\oplus A_{L^G}$.\
The resulting extension $K_C\supseteq L_G\oplus L^G$ is isomorphic to $L$ precisely when the discriminant form 
on $A_{K_C}$ is isomorphic to the discriminant form on $A_L$, since there is only one lattice in the genus of $L$.\
Two lattices $K_C$ and $K_{C'}$ determine $O(L)$-conjugate sublattices $L_G$ inside $L$ if, and only if,
$C$ and $C'$ are in the same orbit for the action of $\overline{O}(L_G)\times \overline{O}(L^G)$ on  $A_{L_G}\oplus A_{L^G}$
induced by the natural action of $O(L_G)\times O(L^G)$.\ For a fixed lattice $K_C$, the $O(L)$-conjugacy classes of $G\subseteq O_0(L_G)$
are given by the $F$-conjugacy classes where $F\subseteq O(L_G)$ is the image of the projection on the first factor of the stabilizer of $C$ under
the action of $O(L_G)\times O(L^G)$.

\smallskip

Since $L^G$ and $L_G$ are both primitive sublattices of $L$, the code $C$ must have the form $C=\{(x,y)\mid x \in A_{L_G},\ y\in A_{L^G}\}$
and satisfy $(x,0)\in C \Rightarrow x=0$ and $(0,y)\in C \Rightarrow y=0$.\ Since $L$ has a discriminant group $A_L$ of order $2$, this implies
that $2|C|^2=|A_{L_G}||A_{L^G}|$, and either 
\begin{itemize}
\item[(1)] $C=\{(x,\gamma(x))\mid x\in A_{L_G}\}$ where $\gamma:A_{L_G}\rightarrow  A_{L^G}$ is a group monomorphism with $q_{L^G}\circ \gamma
=-q_{L_G}$, or
\item[(2)]  $C=\{(\gamma'(y),y)\mid y\in A_{L^G}\}$ where $\gamma':A_{L^G}\rightarrow  A_{L_G}$ is a group monomorphism with $q_{L_G}\circ \gamma'
=-q_{L^G}$.
\end{itemize}
Assume we are in case (1).\ Then $|A_{L^G}/\gamma(A_{L_G})|=2$.\ Since $\gamma(A_{L_G})$ 
is a non-degenerate subspace of $A_{L^G}$ with respect to $q_{L^G}$, there is an orthogonal
decomposition $A_{L^G}=\gamma(A_{L_G})\oplus  \gamma(A_{L_G})^{\perp}$ and $\gamma(A_{L_G})^{\perp}$ 
is generated by an element $v_\gamma\in A_{L^G}$ of order $2$ satisfying $q\vert_{\langle v_\gamma \rangle} \cong q_L$.\
Thus in order to describe the $\overline{O}(L_G)\times \overline{O}(L^G)$-orbits of codes $C$, we first determine
the $\overline{O}(L^G)$-orbits of splittings $A_{L^G}=\gamma(A_{L_G})\oplus {\langle v_\gamma \rangle}$ 
and then the $\overline{O}(L_G) \times S$-orbits of maps
$\gamma: A_{L_G} \longrightarrow \gamma(A_{L_G})$ as in (1), where $S$ is the stabilizer of 
$v_\gamma$ in $\overline{O}(L^G)$ acting on $\gamma(A_{L_G})$.\
Fixing some $\gamma$ allows us to identify $S$ with a subgroup of $O(A_{L_G})$ and the $\overline{O}(L_G) \times S$-orbits
with the double cosets of the pair $(\overline{O}(L_G), S)$ in  $O(A_{L_G})$.

For case (2), $\gamma'(A_{L^G})$ is a subgroup of index $2$ in $A_{L_G}$ and there is an orthogonal
decomposition $A_{L_G}=\gamma'(A_{L^G})\oplus  \gamma'(A_{L^G})^{\perp}$ where $\gamma'(A_{L^G})^{\perp}$ 
is generated by an element $w_{\gamma'}\in A_{L^G}$ of order $2$ with $q\vert_{\langle w_{\gamma'} \rangle} \cong q_L$.\
So to describe the $\overline{O}(L_G)\times \overline{O}(L^G)$-orbits of codes $C$, we first determine
the $\overline{O}(L_G)$-orbits of splittings $A_{L_G}=\gamma'(A_{L^G})\oplus {\langle w_{\gamma'} \rangle}$ 
and then the $S \times \overline{O}(L^G)$-orbits of maps
$\gamma': A_{L^G} \longrightarrow \gamma'(A_{L^G})$ as in (2), where now $S$ is the stabilizer of 
$w_{\gamma'}$ in $\overline{O}(L_G)$ acting on $\gamma'(A_{L^G})$.\
Fixing a $\gamma'$ again permits us to identify $ \overline{O}(L^G)$ with a subgroup of $O(\gamma'(A_{L^G}))$ and 
the $S \times \overline{O}(L^G) $-orbits with the double cosets of the pair $(S,\overline{O}(L^G))$ in $O(\gamma'(A_{L^G}))$.\

For the above discussion see also~\cite{Nikulin}, especially Proposition 1.5.1.

\smallskip

We are now ready for the proof of Lemma~\ref{argdiscriminant}.\ 
We will apply the preceding discussion with $G=\langle g \rangle$.

\pf We have
\begin{eqnarray*}
\alpha(L_g) & := & 24-{\rm rk}\, L_g - {\rm rk}\, A_{L_g} \ = \ 24-{\rm rk}\, {\Lambda}_g - {\rm rk}\, A_{\Lambda_g}\\
& = &  {\rm rk}\, {\Lambda^g}- {\rm rk}\, A_{\Lambda^g }\  = \ {\rm rk}\, A_{\Lambda^g}- {\rm rk}\, A_{\Lambda^g}\ =\ 0.
\end{eqnarray*}
This means that $L_g$ is a lattice of type (c), and the glueing of $L_g$ with $L^g$ is described by case (2),
i.e., by an injective map $\gamma': A_{L^g}\rightarrow A_{L_g}\cong A_{\Lambda^g}$.\ This proves part (a) of Lemma \ref{argdiscriminant}.

For part (b), note that $A_{L^g}=\gamma'(A_{L^g})\oplus \langle w_{\gamma'} \rangle$.\ Since $q_L(w_{\gamma'})=\frac{3}{2}$,
it follows that $q_{L_g}$,  and therefore also $q_{\Lambda^g}$, has $\frac{3}{2}$ in its image. \qed

\medskip

Returning to a general group $G$, after fixing a pair $(L_G, G)$ we proceed along the lines of
 the preceding discussion according to the following five steps:
\vspace{-2mm}
\begin{enumerate}
\item Fix one of the constructions (1) or (2).
\item Determine all lattices in the genus for $L^G$ (which is uniquely determined by the genus of $L_G$).
\item Determine the $\overline{O}(L^G)$-orbits of splittings  $A_{L^G}=\gamma(A_{L_G})\oplus {\langle v_\gamma \rangle}$ 
       resp.\ the $\overline{O}(L_G)$-orbits of splittings $A_{L_G}=\gamma'(A_{L^G})\oplus {\langle w_{\gamma'} \rangle}$.
\item Determine the $\overline{O}(L_G){-}S$ double cosets in  $O(A_{L_G})$ for construction (1), and the 
$S{-}\overline{O}(L^G)$ double cosets in $O(\gamma'(A_{L^G}))$ for construction (2).
\item Determine the $O(L)$-conjugacy classes of $G$ for each double coset.
\end{enumerate}

\paragraph{Construction (1).}
Note that we only have to consider lattices of type (a) or~(b), i.e., with $\alpha(K)\geq 1$.\ Indeed, for lattices of type~(c)
we have ${\rm rk}\, L^G = 23-{\rm rk}\,  L_G <  24-{\rm rk}\, L_G ={\rm rk}\, A_{L_G}$, whence we cannot
embed $A_{L_G}$ into $A_{L^G}$ using $\gamma$.

\medskip

First  consider the case  $\rk L_G=20$.\
From Table \ref{Lisoclasses} we see that there are $13$ lattices of type (a) corresponding to the $G\subseteq M_{23}$
with four orbits in the usual action on $24$ letters, and the two lattices of type (b) corresponding 
to the two maximal $S$-lattice groups.\
In this case, $L^G$ has to be positive-definite of rank~$3$ and the quadratic form $q_{L^G}$ is equivalent to $q_{L_G}\oplus q_{L}$.\  
Recall (proof of Lemma~\ref{argdiscriminant}) that $q_L(x)=\frac{3}{2}$ for the non-zero element $x$ in $A_L$.\
This uniquely determines the genus of $L^G$, and the corresponding lattices
$L^G$ can be read off from the Brandt-Intrau tables of positive definite ternary forms~\cite{BI}.\ 
The result is listed in Table \ref{complementlattice1}.

\smallskip

If $G\cong M_{10}$ and $L^G$ is the lattice with Gramian matrix ${\rm Diag}(2,4,30)$, there are two
$O(L^G)$-orbits of splittings  $A_{L^G}=\gamma(A_{L_G})\oplus {\langle v_\gamma \rangle}$, whereas in
all other cases there is a unique $O(L^G)$-orbit of splittings.

\smallskip

We assert that there is a \emph{unique} $\overline{O}(L_G){-}S$ double coset in $O(A_{L_G})$.\ 
If $G\not \cong S_6$ then $\overline{O}(L_G)$ coincides with $O(A_{L_G})$, and the assertion follows. 
A computation shows that the same result holds if $G \cong S_6$.

\smallskip

Finally, let $\overline{F}\subseteq\overline{O}(L_G) $ be the projection onto the first factor of the stabilizer in $\overline{O}(L_G) \times S$ 
of the identity element under the double coset action.\
Let $F\subset O(L_G)$ be the inverse image of $\overline{F}$ under the natural projection.\
We consider the conjugation action of $F$ on 
the set of subgroups $H\subseteq O_0(L_G)$ with $L_G^H=\{0\}$.\ 
For lattices of type (a), the orbits agree with the $\overline{O}(L_G)$ classes,
whereas for the two lattices of type (b), there are more orbits.\ 
We list the number of conjugacy classes $H$ for both groups
in the last two columns of Table~\ref{complementlattice1}. 

\medskip

Now assume that $\rk L_G<20$.\ 
In this case, $L^G$ has rank at least $4$ and must be indefinite.\ 
We claim that the lattice $L^G$ exists and is unique.\ Indeed, the $41$ lattices $L_G$ of type~(a) are sublattices
of the K3-lattice $N\cong  E_8(-1)^2\oplus U^3$.\ If $N^G$ is an orthogonal complement of $L^G$ in $N$ then
we let $L^G=N^G \oplus A_1(-1)$.\ If the uniqueness criterion of Theorem~1.7 of~\cite{Ha} (which follows from
Eichler's theory of Spinor genera, cf.~\cite{Nikulin})
applies to $N^G$ then it also applies to  $L^G$.\ Therefore, the explicit verification in Section~6 of~\cite{Ha} establishes
uniqueness  for all $30$ lattices $L^G$ of type (a) and rank $>4$.\
For the $11$ lattices $L^G$ of rank $4$, we apply Theorem 1.7 of~\cite{Ha} directly.\
For the two lattices $L_G$ of type~(b) and rank $19$ and $18$, one gets for $L^G$ the two lattices with Gramian matrix
\[ 
{\footnotesize\renewcommand{\arraystretch}{0.7}\left(\begin{array}{cccc}
 0 & 3 & 0 & 0 \\ 3 & 0 & 0 & 0 \\ 0 & 0 & 6 & 0 \\ 0 & 0 & 0 & 6
\end{array}\right)} \quad \mbox{and}\qquad 
{\footnotesize\renewcommand{\arraystretch}{0.7}\left(\begin{array}{ccccc}
 0 & 3 & 0 & 0 & 0 \\ 3 & 0 & 0 & 0 & 0 \\ 0 & 0 & 0 & 3 & 0 \\ 0 & 0 & 3 & 0 & 0 \\ 0 & 0 & 0 & 0 & 6
\end{array}\right)},
\]
respectively. Again, the uniqueness follows from Theorem~1.7 of~\cite{Ha}.

\begin{conj}\label{fullOA}
For $L_G$ a lattice of type (a) or (b) and ${\rm rk}\,L^G\geq 4$, we have $\overline{O}(L^G)=O(A_{L^G})$.
\end{conj}
By applying the criterion from Theorem 1.14.2 in~\cite{Nikulin}, one checks that the conjecture holds for all $L_G$ 
of type~(a) and~(b) in Table~\ref{Lisoclasses},
except perhaps for the  cases with number 
\[n=18,\, 26,\, 31,\, 39,\, 41,\, 47,\, 48,\, 54,\, 55,\, 61.  \]

Unfortunately, most lattice functions of Magma are presently implemented for definite lattices only,
so we cannot verify this conjecture by a computer calculation without extra programming work.\ 
One could do this more theoretically as in~\cite{Ha}, Section 7, by using the strong approximation theorem,
but we refrain from investigating this here.

If Conjecture~\ref{fullOA} holds, it is clear that there is a single
$O(L^G)$-orbit of orthogonal splittings  $A_{L^G}=\gamma(A_{L_G})\oplus {\langle v_\gamma \rangle}$.\ Furthermore, we would have
$S=O(A_{L_G})$ and $F=O(L_G)$.\ It then follows from Theorem~\ref{GsetofK}
that for each  $(L_G,G)$, there is a \emph{unique} $O(L)$-conjugacy class 
of subgroups $G$.

\noindent{\bf Remark:} Since Conjecture~\ref{fullOA} holds in particular for the three K3-lattices $L_G$ of rank $19$ 
(numbers $44$, $46$ and $52$ in Table~\ref{Lisoclasses}) for which there are two lattices in the genus of $N^G$, 
it follows that the corresponding symplectic group actions on $K3$ are \emph{not} deformation equivalent.\ However, the induced 
symplectic group actions on $\K32$ are so.

\medskip

Our calculations have shown:
\begin{thm}\label{construction1}
Let $G\subseteq {\rm Co}_0$ be an admissible subgroup.\
Then $(L_G,G)$ can be realized as the coinvariant lattice for a subgroup $G\subseteq O(L)$ by construction~$(1)$
if, and only if, $\alpha(L_G)\geq 1$.\
If $\rk L_G=20$, the number of corresponding conjugacy classes of $G\subseteq O(L)$ can be read off from Table~\ref{complementlattice1}.\
If $\rk L_G<20$, there is a \emph{unique} such class, provided that either $L_G$ is not one of the cases no.\
$18$, $26$, $31$, $39$, $41$, $47$, $48$, $54$, $55$, $61$ in Table~\ref{Lisoclasses}, or
Conjecture~\ref{fullOA} holds. \qed
\end{thm}

\begin{table}\caption{Conjugacy classes in $O(L)$: ${\rm rank}\, L^G=3$, construction~(1)}\label{complementlattice1}
$$\footnotesize \renewcommand{\arraystretch}{0.8} 
\begin{array}{rl|rrrr|crrr|rr}
\hbox{No.} & G & |A_{L_G}| & |O(L_G)| & |\overline{O}(L_G)| & |O(A_{L_G})| &   \hbox{$L^G$}  & \!\! |O(L^G)|  & |\overline{O}({L^G})| &  |O(A_{L^G})| & \#H & \#\hbox{orbs} \\
\noalign{\hrule height0.8pt}

1   & {\rm L}_2(11)
& 121 &
15840 & 24 & 24 &
{\scriptsize\left(\renewcommand{\arraystretch}{0.5} \begin{array}{rrr}
 2 & 1 & 0 \\
 1 & 6 & 0 \\
 0 & 0 & 22 
\end{array}\right)}
& 8 & 4  &  24 
& 3 & 3
\\

& & &
& & &

{\scriptsize\left(\renewcommand{\arraystretch}{0.5}\begin{array}{rrr}
 6 & 2 & 2 \\
 2 & 8 & -3 \\
 2 & -3 & 8 
\end{array}\right)} 
& 12 & 12  &  24
& 3 & 3

\\ \hline

2   & {\rm L}_3(4)
& 84 &
483840 & 24 &  24 &
{\scriptsize\left(\renewcommand{\arraystretch}{0.5}\begin{array}{rrr} 
 2 & 0 & 0 \\
 0 & 10 & 4 \\
 0 & 4 & 10 
\end{array}\right)}
& 8 & 4 & 24
& 1 & 1
\\

& & &
& & &

{\scriptsize\left(\renewcommand{\arraystretch}{0.5}\begin{array}{rrr} 
 4 & 2 & 0 \\
 2 & 4 & 0 \\
 0 & 0 & 14 
\end{array}\right)}
& 24 & 24 & 24
& 1 & 1

\\ \hline

3   & A_7
& 105 &
 20160 & 8 & 8 &
{\scriptsize\left(\renewcommand{\arraystretch}{0.5}\begin{array}{rrr} 
 2 & 1 & 0 \\
 1 & 2 & 0 \\
 0 & 0 & 70
\end{array}\right)}
& 24 & 4  &  8 
& 1 & 1

\\
& & &
& & &
{\scriptsize\left(\renewcommand{\arraystretch}{0.5}\begin{array}{rrr} 
 2 & 0 & 1 \\
 0 & 6 & 0 \\
 1 & 0 & 18 
\end{array}\right)}
 & 8 & 4  &  8
& 1 & 1
\\

& & &
& & &
{\scriptsize\left(\renewcommand{\arraystretch}{0.5}\begin{array}{rrr} 
 4 & 2 & 1 \\
 2 & 6 & 3 \\
 1 & 3 & 12 
\end{array}\right)}
&  4 & 4  &  8
& 1 & 1
\\

& & &
& & &

{\scriptsize\left(\renewcommand{\arraystretch}{0.5}\begin{array}{rrr} 
6 & 3 & 1 \\
3 & 6 & 1 \\
1 & 1 & 8 
\end{array}\right)}
& 4 & 4  &  8 
& 1 & 1

\\ \hline

4   & \Z_2^3{:}{\rm L}_2(7)
& 112 &
21504 & 16 & 16 &
{\scriptsize\left(\renewcommand{\arraystretch}{0.5}\begin{array}{rrr} 
 4 & 0 & 0 \\
 0 & 6 & 2 \\
 0 & 2 & 10 
\end{array}\right)}
& 4 & 4  &  16
& 3 & 3

\\ \hline

5   & \Z_2\times {\rm L}_2(7)
& 196 &
 10752 & 32 & 32 &
{\scriptsize\left(\renewcommand{\arraystretch}{0.5}\begin{array}{rrr} 
 2 & 0 & 0 \\
 0 &14 & 0 \\
 0 & 0 &14 
\end{array}\right)}
&  16 & 8  &  32
& 3 & 3
\\

& & &
& & &

{\scriptsize\left(\begin{array}{rrr} 
 4 & 2 & 0 \\
 2 & 8 & 0 \\
 0 & 0 &14
\end{array}\right)}
& 8 & 8  &  32 
& 3 & 3
\\ \hline

6   & \Z_2^4{:}A_6
& 96 &
92160 & 16 & 16 &
{\scriptsize\left(\renewcommand{\arraystretch}{0.5}\begin{array}{rrr} 
 2 & 0 & 0 \\
 0 & 4 & 0 \\
 0 & 0 & 24 
\end{array}\right)}
& 8 & 4 & 16
& 7 & 7
\\

& & &
& & &

{\scriptsize\left(\renewcommand{\arraystretch}{0.5}\begin{array}{rrr} 
4 & 0 & 0 \\
0 & 6 & 0 \\
0 & 0 & 8 
\end{array}\right)}
& 8 & 8 & 16
& 7 & 7

\\ \hline

7   & \Z_2^4{:}S_5
& 160 &
 30720 & 16 & 16 &
{\scriptsize\left(\renewcommand{\arraystretch}{0.5}\begin{array}{rrr} 
 2 & 0 & 0 \\
 0 & 4 & 0 \\
 0 & 0 & 40 
\end{array}\right)}
& 8 & 4  &  16 
& 2 & 2
\\

& & &
& & &

{\scriptsize\left(\renewcommand{\arraystretch}{0.5}\begin{array}{rrr} 
 4 & 0 & 0 \\
 0 & 8 & 0 \\
 0 & 0 & 10 
\end{array}\right)}
&  8 & 8  &  16 
& 2 & 2

\\ \hline

8   & S_6
& 180 &
23040 & 32 & 96 &
{\scriptsize\left(\renewcommand{\arraystretch}{0.5}\begin{array}{rrr} 
 4 & 2 & 0 \\
 2 & 4 & 0 \\
 0 & 0 & 30 
\end{array}\right)}
& 24 & 24 & 96
& 1 & 1

\\ \hline

9   &  M_{10}
&  120 &
 5760 & 8 & 8 &
{\scriptsize\left(\renewcommand{\arraystretch}{0.5}\begin{array}{rrr} 
 2 & 0 &  0 \\
 0 & 4 & 0 \\
 0 & 0 & 30 
\end{array}\right)}
& 8 & 4 & 16
& 1 & 1 

\\

& & &
& & &
& 8 & 4 & 16
& 1 & 1

\\

& & &
& & &

{\scriptsize\left(\renewcommand{\arraystretch}{0.5}\begin{array}{rrr} 
 4 & 2 & 0 \\
 2 & 6 & 0 \\
 0 & 0 & 12 
\end{array}\right)}
& 8 & 8 & 16
& 1 & 1

\\ \hline

10   & (\Z_3\times A_5){:}\Z_2
& 225 &
17280 & 48 & 48 &
{\scriptsize\left(\renewcommand{\arraystretch}{0.5}\begin{array}{rrr} 
 4 & 1 & 0 \\
 1 & 4 & 0 \\
 0 & 0 & 30 
\end{array}\right)}
& 8 & 8  &  48 
& 5 & 5
\\

& & &
& & &

{\scriptsize\left(\renewcommand{\arraystretch}{0.5}\begin{array}{rrr} 
 6 & 0 & 0 \\
 0 & 10 &  5 \\
 0 & 5 & 10 
\end{array}\right)}
& 24 & 24  &  48
& 5 & 5
\\ \hline 

11  &  Q(\Z_3^2{:}\Z_2)
& 192 &
36864 & 32 & 32 &

{\scriptsize\left(\renewcommand{\arraystretch}{0.5}\begin{array}{rrr} 
 6 & 2 & 2 \\
 2 & 6 & -2 \\
 2 & -2 & 14 
\end{array}\right)}
& 8 & 8  &  64
& 10 & 10

\\ \hline

12  &  \Z_2^4{:}(S_3\times S_3)
& 288 &
 36864 & 64 & 64 &

{\scriptsize\left(\renewcommand{\arraystretch}{0.5}\begin{array}{rrr} 
 4 & 0 & 0 \\
 0 & 6 & 0 \\
 0 & 0 & 24 
\end{array}\right)}
& 8 & 8  &  64 
& 2 & 2

\\ \hline

13  &  \Z_3^2{:}QD_{16}
& 216 &
3456 & 24 & 24 &

{\scriptsize\left(\renewcommand{\arraystretch}{0.5}\begin{array}{rrr}
 4 & 2 & 0 \\
 2 & 10 & 0 \\
 0 & 0 & 12
\end{array}\right)}
&  8 & 8  & 48
& 2 & 2 
  
\\ \hline
\noalign{\hrule height0.8pt}

14  &  3^{1+4}:{2}.2^2
& 108 &
186624 & 96 & 96 &

{\scriptsize\left(\renewcommand{\arraystretch}{0.5}\begin{array}{rrr} 
 6 & 0 & 0 \\
 0 & 6 & 0 \\
 0 & 0 & 6 
\end{array}\right)}
& 48 & 48  &  288
& 7 & 8

\\ \hline

15   & 3^4:A_6
& 81 &
4199040 & 144 & 144 &

{\scriptsize\left(\renewcommand{\arraystretch}{0.5}\begin{array}{rrr} 
 6 & 3 & 0 \\
 3 & 6 & 0 \\
 0 & 0 & 6 
\end{array}\right)}
& 24 & 24 & 144
& 40 & 71

\\ \hline
\noalign{\hrule height0.8pt}

\end{array}
$$
\end{table}

\paragraph{Construction (2).}\phantom{x}

\smallskip

We first study the $\overline{O}(L_G)$-orbits of splittings $A_{L_G}=\gamma'(A_{L^G})\oplus {\langle w_{\gamma'} \rangle}$.\ 
Thus we searched for elements $w_{\gamma'}\in A_{L_G}$ of order $2$ with $q_{L_G}(w_{\gamma'})=3/2 \pmod{2}$.\
Such elements exist only for the eleven lattices $L_G$ with the numbers
$$n=9,\,14,\,18,\,20,\,27,\,28,\,33,\,38,\,40,\,45,\,50$$
in Table~\ref{Lisoclasses}, and in each case there is a \emph{unique} $\overline{O}(L_G)$-orbit. 

The discriminant form of $L^G$ is equal to $-q_{L_G}|_{w_{\gamma'}^\perp}$, so in each case the genus of $L^G$
is uniquely determined.\
For the four cases $n=58$, $45$, $38$ and $27$ with
$\rk L_G=20$, Table~\ref{complementlattice2} lists  $L^G$
together with additional information concerning $O(L_G)$ and $O(L^G)$.\
In the other seven cases $L^G$ is indefinite, and one easily checks that $L^G$ exists and is
unique in these cases.
\medskip

For reasons which will become apparent in the next section, we will not further analyze the exact number
of $O(L)$-conjugacy classes of realizations of the group lattice $(L_G,G)$ inside $L$.\ We have established:
\begin{thm}\label{construction2}
Let $G\subseteq {\rm Co}_0$ be an admissible subgroup.\
Then $(L_G,G)$ can be realized as the coinvariant lattice for a subgroup $G\subseteq O(L)$ by construction~$(2)$
if, and only if, $L_G$ is one of cases no.\
$9$, $14$, $18$, $20$, $27$, $28$, $33$, $38$, $40$, $45$, $50$
in Table~\ref{Lisoclasses}.
\end{thm}

We note that there are admissible subgroups $G\subseteq {\rm Co}_0$ for
which $(L_G,G)$ {\it cannot\/} be realized in this way.

\medskip

\begin{table}\caption{Conjugacy classes in $O(L)$: ${\rm rank}\, L^G=3$, construction ({\rm 2})}\label{complementlattice2}
$$\footnotesize \renewcommand{\arraystretch}{0.8} 
\begin{array}{rl|rrrr|crrr}
\hbox{No.} & G & |A_{L_G}| & |O(L_G)| & |\overline{O}(L_G)| & |O(A_{L_G})| &   \hbox{$L^G$}  & |O(L^G)|  & |\overline{O}({L^G})| &  |O(A_{L^G})| \\
\noalign{\hrule height0.8pt}
1   & M_{10}
& 120 &
5760 & 8 & 8 &
{\scriptsize\left(\renewcommand{\arraystretch}{0.5} \begin{array}{rrr}
 2 & 0 & 1 \\
 0 & 4 & 0 \\
 1 & 0 & 8  
\end{array}\right)}
& 8 & 4 & 8
\\

& & &
& & &

{\scriptsize\left(\renewcommand{\arraystretch}{0.5}\begin{array}{rrr}
 4 & 2 & 2 \\
 2 & 4 & 1 \\
 2 & 1 & 6 
\end{array}\right)} 
& 8 & 8 & 8
\\ \hline

2   & \Z_3^2{:}QD_{16}
& 216 &
3456 & 24 & 24 &
{\scriptsize\left(\renewcommand{\arraystretch}{0.5} \begin{array}{rrr}
 4 & 2 & 2 \\
 2 & 4 & 1 \\
 2 & 1 & 10  
\end{array}\right)}
& 8 & 8 & 24
\\ \hline
\noalign{\hrule height0.8pt}

3 & 32\,\Gamma_7a_3
& 256 &
196608 & 256  &  256 &

{\scriptsize\left(\renewcommand{\arraystretch}{0.5}\begin{array}{rrr} 
 2 & 0 & 0 \\
 0 & 8 & 0 \\
 0 & 0 & 8
\end{array}\right)}
& 16 & 8 & 32

\\ \hline

4 & 64\, \Gamma_{23}a_3
& 128 &
36864 & 96  & 96 &

{\scriptsize\left(\renewcommand{\arraystretch}{0.5}\begin{array}{rrr} 
 2 & 0 & 0 \\
 0 & 2 & 0 \\
 0 & 0 & 16
\end{array}\right)}
& 16 & 4 & 16

\\ \hline
\noalign{\hrule height0.8pt}

\end{array} $$

\end{table}

\noindent{\bf Remark:} 
We also have applied this method to prove the uniqueness 
of the conjugacy class of $G$ in the isometry group $O(N)$ of the K3-lattice $N$
for the $11$ K3-lattices $(L_G,G)$ of rank $19$ and the corresponding $14$ rank $3$ lattices $N^G$.\
This provides a somewhat more systematic approach for these cases 
compared to Section~8.2 of~\cite{Ha}.


\section{Symplectic actions on $\K32$}\label{geometricreal}

In this section, we will show that the groups $G\subseteq O(L)$ that can be realized as isometries 
of a hyperk\"ahler manifold of type $\K32$ are exactly those 
for which $L_G$ is a lattice of type (a) or (b). 

\medskip

Regarding birational maps of a hyperk\"ahler manifold of type $\K32$, there is the following
characterization by Mongardi~(\cite{Mon-K32inv}, Thm.~3.6):
\begin{thm}
Let $G\subseteq O(L)$ be a finite subgroup, and suppose that the coinvariant lattice
$L_G$ is negative-definite and contains no roots of norm $-2$.\ Then $G$ is induced
by a group of birational transformations of some hyperk\"ahler manifold of type $\K32$.
\end{thm}
Thus we have:
\begin{cor}
The conjugacy classes of groups $G \subseteq O(L)$ found in Section~\ref{subol}
arise as finite group 
of birational transformations of some hyperk\"ahler manifold of type $\K32$.
\end{cor}
We note that there are additional subgroups of $G \subseteq O(L)$ arising from birational transformations.\
The groups $G \subseteq O(L)$ in Section~\ref{subol} have only been
constructed from admissible subgroups $G\subseteq {\rm Co}_0$.

\medskip
Based on work on the global Torelli theorem for hyperk\"{a}hler manifolds,
Mongardi gave the following criterion 
regarding symplectic automorphisms of type $K3^{[n]}$ (\cite{Mon-K3n}, Thm.~4.1) which we state 
here for the case $\K32$:
\begin{thm}\label{realcriterion}
Let $G\subseteq O(L)$ be a finite group.\ Then $G$ is induced by a
group of symplectic automorphisms for some hyperk\"ahler manifold $X$ of type $\K32$
if, and only if, the following holds:
\vspace{-3mm}
\begin{itemize}
\itemsep0em
\item $L_G$ is negative-definite;
\item $L_G$ contains no numerical wall divisor.
\end{itemize}
\end{thm}

The numerical wall divisors for $K3^{[n]}$ have been 
discussed in~\cite{BHT} and \cite{Mon-wall} based on work in~\cite{BM}. 
For manifolds of type $K3^{[2]}$ one has (\cite{Mon-wall}, Prop.~2.12):
\begin{thm}\label{walldivisors}
Let $X$ be a hyperk\"ahler manifold of type $\K32$. Then the numerical wall divisors
are the following vectors in the Picard sublattice of $L$:
\vspace{-3mm}
\begin{itemize}
\itemsep0em
\item vectors $v$ of norm $v^2=-2$;
\item vectors $v$ of norm $v^2=-10$ and $v/2\in L^*$.
\end{itemize}
\end{thm}

Since $L_G$ is a sublattice of the Picard lattice, we have to check for which of the groups $G\subseteq O(L)$
from Section~\ref{subol} the lattice $L_G$ contains no vectors $v$ of norm~$-10$ such that
$v/2\in L^*$.
\begin{prop}\label{wallcriterion}
If $L$ is obtained from $L_G\oplus L^G$ by glueing construction (1) in Section~\ref{subol},
then $L_G$ contains \emph{no} numerical wall divisor.\ If $L$ is obtained by glueing construction (2) then $L_G$ contains
a numerical wall divisor if, and only if, it contains a vector $v$ of norm $-10$ such that
$\overline{v/2}=w_{\gamma'}$ in $A_{L_G}$.
\end{prop}

\pf By Theorem~\ref{walldivisors}, we have to check if $L_G$ contains 
a vector $v$ of norm $-10$ such that $v/2\in L^*$.

In case (1), $L$ contains vectors of the form $(a,b)\in L_G\oplus L^G$ and the cosets 
$(x,\gamma(x))$, $x\in A_{L_G}$.\ Thus a vector $(v/2,0)\in L_G^*\oplus (L^G)^*$ is contained in $L^*$ if and only
if $(v/2,L_G^*)\subset \Z$, i.e.~$v/2\in L_G$.\ But this is impossible since the norm of
$v/2$ is $-5/2$ if $v$ has norm $-10$.

In case (2), $A_{L_G}=\gamma'(A_{L^G})\oplus \langle w_{\gamma'} \rangle $,  $L/(L_G\oplus L^G)=\{
(\gamma'(y),y) \mid y\in A_{L^G}\}$ so that $(w_{\gamma'},0)$ generates $A_L=L^*/L$.\ 
Thus a vector $(v/2,0)\in L_G^*\oplus (L^G)^*$ is contained in $L^*$ if and only
if $(\overline{v/2},\gamma'(A_{L^G}))=0$, i.e.~$\overline{v/2}\in \langle w_{\gamma'}\rangle $ in $A_{L_G}$. 
The case $v/2\in L_G$ is again impossible. $\hfill \Box$

\medskip

\begin{thm}\label{groupreal}
Let $G\subseteq {\rm Co}_0$ be an admissible subgroup.\ Then $G$ is induced by a group of symplectic automorphisms
for some hyperk\"ahler manifold of type $\K32$ such that $(\Lambda_G(-1),G))\cong (L_G,G)$ if, and only if, 
$\alpha(L_G)\geq 1$.
\end{thm}
\pf If there exists a realization of $(L_G,G)$ as coinvariant lattice by a group $G\subseteq O(L)$ using construction~(1), 
then by Proposition~\ref{wallcriterion} $L_G$ contains \emph{no} numerical wall divisor.

We will show that any realization of $(L_G,G)$ using construction~(2) contains a numerical wall divisor.\
With this in mind, for each of the eleven lattices $L_G$ as in Theorem~\ref{construction2} we searched randomly 
for vectors $v$ in the dual lattice $L_G^*$  of norm $-5/2$ such that $2v\in L_G$.\ We always found such a $v$.\
As discussed in the last section, in each case there is a \emph{unique} $O(L_G)$-orbit of norm $3/2$ vectors 
$w_{\gamma'}$ in $A_{L_G}$, i.e.\ we can assume that $\overline{v/2}=w_{\gamma'}$.\ Thus by 
Proposition~\ref{wallcriterion},  $L_G$ contains a numerical wall divisor.

It  therefore follows from Theorem~\ref{realcriterion} that $G$ is induced by a group of symplectic automorphisms
if, and only if, there is a realization of $(L_G,G)$ using construction~(1).\ By Theorem~\ref{construction2}, these are
exactly the admissible groups for which $\alpha(L_G)\geq 1$. $\hfill \Box$

Combining Theorem~\ref{groupreal} with Theorem~\ref{admissiblegroupclasses} we have:
\begin{thm}\label{sympauto}
The finite groups $G$ arising as symplectic automorphisms of a hyperk\"ahler manifold of type $\K32$
are:
\vspace{-2mm}
\begin{itemize}
\itemsep0em
\item[(a)] subgroups of  $M_{23}$  with at least four orbits in the natural action on $24$ elements, 
\item[(b)] subgroups of 
$3^{1+4}{:}2.2^2$ and $3^4{:}A_6$ associated to the corresponding \S-lattices.   
\end{itemize}
\end{thm}
At this juncture, we have established Theorem~\ref{thmmain} and Theorem~\ref{thmreal}.

\bigskip

Let $X_i$ be hyperk\"ahler manifolds of type $\K32$ with  finite groups of
symplectic automorphism $G_i\subseteq{\rm Aut}\,X_i$ ($i=1$, $2$).  We say that $(X_1,G_1)$ is 
\emph{deformation equivalent\/} to $(X_2,G_2)$
if there exists a flat family $\chi\longrightarrow B$ over a connected base $B$ with hyperk\"ahler manifolds of type $\K32$ as fibers,
together with a fiberwise symplectic action of a group $G$ such that $(X_1,G_1)$ and $(X_2, G_2)$ are 
\emph{isomorphic\/} to the action of $G$ at certain fibers.

The following result follows from~\cite{Mon-natural}; cf.\ the proof of Corollary~5.2 in~\cite{Mon-K3n}.
\begin{thm}\label{deformequi}
Two hyperk\"ahler manifolds $X_i$ of type $\K32$ with symplectic automorphism groups $G_i$ ($i=1$, $2$), are
equivariantly deformation equivalent if, and only if, the associated group lattices $(L_i,G_i)$ are isomorphic.
\end{thm} 
Together with Theorem~\ref{groupreal} and the enumeration results from the previous section,
this proves Theorem~\ref{thmdeformationclasses}.\ Indeed,
there are $198-13=185$ group lattices $(L_G,G)$ that can be realized.\ For the $88$ groups $G$
with $\rk\, L_G=20$, there are $146$ conjugacy classes in $O(L)$ (see the last column
in Table~\ref{complementlattice1}), giving at least $185- 88 + 146= 243$ conjugacy classes in all.\
If Conjecture~\ref{fullOA} holds, there are also exactly $243$ deformation classes.

We finally note that Theorem~\ref{thmlatticeclasses} is equivalent to Theorem~\ref{groupreal} by Theorem~\ref{realcriterion},
Theorem~\ref{walldivisors} and Theorem~\ref{deformequi}. 

\smallskip

\noindent{\bf Remarks:}
The admissible conjugacy classes can also be determined with the help of Theorem~\ref{realcriterion} and Theorem~\ref{walldivisors}
instead of the method we used in Section~\ref{admissibleconj}.\ However, we will need 
the exact fixed-point configuration (described in Table~\ref{Fixpointset}) in the final section.\ 
Also, Theorem~\ref{realcriterion} became  only available after main parts of the paper had been written.\ 
Theorem~\ref{groupreal} was conjectured in~\cite{Mon-K3n}.

\paragraph{Explicit examples of group actions.}
It is known that certain maximal groups $G$ can be realized as symplectic actions
on a $\K32$ through induced actions on Fano schemes of lines of 
certain cubic fourfolds $S\subset {\bf C}P^5$, cf.~\cite{Fu}, \cite{Mon-thesis} Ch.~4 and the references therein.\
We collect these examples in Table~\ref{Fano}.\ 
In addition, we consider two apparently new examples, with 
$G\cong 3^{1+4}{:}2.2^2$ and $M_{10}$.\ These are discussed in the next few paragraphs.

Concerning $G\cong 3^{1+4}{:}2.2^2$, first note that the cubic polynomial $f$ in Table~\ref{Fano} 
is invariant under an obvious action of $H=(3^2.S_3\times 3^2.S_3).\Z_2\subseteq {\rm GL}(6,\C)$  
given by permutations and multiplication of the coordinates by cube roots of unity.\
In addition, $f$ is invariant under
the unimodular matrix
{\footnotesize \[ \alpha:= \frac{1}{\sqrt{3}}\left(\renewcommand{\arraystretch}{0.7} 
\begin{array}{cccccc}
\omega & \omega^2 & 1 &  &  & \\
1 & 1 & 1 &  &  & \\ 
\omega^2 & \omega & 1 &  &  & \\ 
& &  & \omega^2 & \omega & 1\\ 
 &  & & \omega^2 & \omega^2 & \omega^2\\ 
&  &  & \omega^2 & 1 & \omega 
\end{array}\right)\qquad
\begin{array}{c}
\\ \\ \\ \\ (\omega=e^{2\pi i/3}).
\end{array}
\]}
An element in ${\rm GL}(6,\C)$ leaving $f$ invariant acts symplectically on the Fano scheme of $S$
if, and only if, it is in ${\rm SL}(6,\C)$ (cf.~\cite{Fu}).\
Calculations show that the projections of $H\cap {\rm SL}(6,\C)$ and $\alpha$ into ${\rm PSL}(6,\C)$ generate
a group isomorphic to $G$.\
Finally, it can be verified that the resulting cubic fourfold $S$ is smooth by 
solving the equations
\[f=\partial f/\partial x_0= \cdots = \partial f/\partial x_5=0,\]
confirming that $0$ is an isolated singularity of $f$.

For the case $G\cong M_{10}$, we start with the containment
$3.A_6\subseteq {\rm SL}(6,\C)$ given by the generators
{\footnotesize \[\left(\renewcommand{\arraystretch}{0.7} 
\begin{array}{cccccc}
 1 &&&&& \\&& 1 &&& \\& 1 &&&& \\&&&& 1 & \\&&& 1 && \\&&&&& 1
\end{array}\right), \qquad
\left(\renewcommand{\arraystretch}{0.7} \begin{array}{cccccc}
& 1 &&&& \\&& \omega &&& \\&&& 1 && \\ \omega^2 &&&&& \\&&&&& 1\\&&&& 1 & 
\end{array}\right) \]}
$\!\!\!$from the Atlas of finite group representations~\cite{Wilson}.\
In addition, we choose the unimodular matrix 
{\footnotesize \[\beta:= \frac{1}{\sqrt{6}}\left(\renewcommand{\arraystretch}{0.7} 
\begin{array}{cccccc}
1 & \omega & \omega^2 & \omega & 1 & \omega \\ \omega^2 & 1 & 1 & \omega & \omega & \omega\\ 
 \omega & 1 & \omega^2 & \omega & \omega^2 & \omega^2\\ 
\omega^2 & \omega^2 & \omega & \omega & 1 & \omega^2\\  1 & \omega^2 & 1 & \omega & \omega^2 & 1\\ 
 \omega^2 & \omega^2 & \omega^2 & \omega^2 & \omega^2 & \omega
\end{array}\right)_. \]}
which normalizes $3.A_6$.\ We found the following cubic polynomial $f$ invariant under the action of $3.A_6$ and $\beta$:
\begin{eqnarray*}
 f &= & (x_1^3+x_2^3+x_3^3+x_4^3+x_5^3+x_6^3) +\frac{1}{5}(-3\zeta^7 - 3\zeta^5 + 3\zeta^4 - 3\zeta^3 + 6\zeta - 3) \times \\
&&  \Bigl[x_1x_2x_3 + x_1x_2x_4 + (\zeta^4 - 1)x_1x_2x_5 + x_1x_2x_6 + (\zeta^4 - 1)x_1x_3x_4 + x_1x_3x_5  \\ 
&&\ \ + x_1x_3x_6 + (\zeta^4 - 1)x_1x_4x_5 - \zeta^4x_1x_4x_6 - \zeta^4x_1x_5x_6 + (\zeta^4 - 1)x_2x_3x_4  \\
&&\ \ + (\zeta^4 - 1)x_2x_3x_5 - \zeta^4x_2x_3x_6 + x_2x_4x_5 + x_2x_4x_6 - \zeta^4x_2x_5x_6 + x_3x_4x_5 \\
&&\ \ - \zeta^4x_3x_4x_6 + x_3x_5x_6 + x_4x_5x_6\Bigr] 
\end{eqnarray*}
with $\zeta =e^{2\pi i/24}$.\
The projection of $3.A_6\langle \beta\rangle$ into ${\rm PSL}(6,\C)$
is isomorphic to $M_{10}$. Again, we verified that $S$ is smooth.

\begin{table}[t]\caption{Fano schemes of cubic fourfolds}\label{Fano}
{$$\renewcommand{\arraystretch}{1.4}
\begin{array}{cl} 
\mbox{Group} & \mbox{Equation for $S$} \\ \hline
 L_2(11) & x_0^3+x_1^2x_5+x_2^2x_4+x_3^2x_2+x_4^2x_1 + x_5^2x_3 \\
  A_7    & x_0^3+x_1^3+x_2^3+x_3^3+x_4^3+x_5^3-(x_0+x_1+x_2+x_3+x_4+x_5)^3 \\
 (\Z_3\times A_5){:}\Z_2  & x_0^2x_1+x_1^2x_2+x_2^2x_3+x_3^2x_0+x_4^3+x_5^3 \\ 
 3^4{:}A_6 & x_0^3+x_1^3+x_2^3+x_3^3+x_4^3+x_5^3 \\
 3^{1+4}{:}2.2^2  & x_0^3+x_1^3+x_2^3+x_3^3+x_4^3+x_5^3 +3(i - 2e^{\pi i/6} - 1)\,( x_0x_1x_2 + x_3x_4x_5) \\
  M_{10} & x_0^3+x_1^3+x_2^3+x_3^3+x_4^3+x_5^3 + \lambda\cdot \tilde\sigma_3(x_0,\,\ldots,\,x_5).
\end{array}$$} 
\end{table}

It would be of interest  to find explicit realizations for the nine remaining 
 maximal groups of Theorem \ref{thmmain} not treated here.


\section{Connections with Mathieu Moonshine}

Suppose that $g\in M_{24}$ belongs to one of the $11$ admissible classes.\
In this section we compare the equivariant complex elliptic genus $ \chi_y(g;q,{\cal L}X)$ of a hyperk\"ahler manifold $X$ of
type $\K32$ with the prediction of  Mathieu Moonshine applied to the second quantized complex elliptic genus
of a $K3$ surface.\ See Theorem~\ref{moonshinegeometric2} below for a precise statement.\
This connection was one of our main motivations for studying  symplectic
automorphisms of $\K32$.


\subsection{Mathieu Moonshine}

Recall \cite{Hirzebruch-Habil} that a \emph{complex genus} in the sense of Hirzebruch is a graded ring homomorphism
from the complex bordism ring into some other graded ring $R$.\ 
For a $d$-dimensional complex manifold $X$ and
a holomorphic vector bundle $E$ on $X$, the $\chi_y$-genus twisted by $E$ is
\[
\chi_y(X, E) := \sum_{p=0}^d \chi(X, \Lambda^pT^*\otimes E)\,y^p.
\]
where $\chi(X, E)=\sum_{q=0}^{d}(-1)^q \dim H^q(X,{\cal O}(E))$.

The complex elliptic genus can be formally defined as the $S^1$-equivariant $\chi_y$-genus 
of the loop space of a manifold:\
$$\chi_y(q,{\cal L}X):= (-y)^{-d/2}\chi_y\bigl(X,\bigotimes_{n=1}^{\infty}\Lambda_{yq^n}T^*
\otimes \bigotimes_{n=1}^{\infty}\Lambda_{y^{-1}q^n}T
\otimes \bigotimes_{n=1}^{\infty}S_{q^n}(T^*\oplus T)\bigr) $$
taking values in $\Q[y^{1/2},y^{-1/2}][[q]]$.\
Here, we let $S_tE=\bigoplus_{i=0}^\infty S^iE\cdot t^i$ and $\Lambda_t E=\bigoplus_{i=0}^\infty \Lambda^iE\cdot t^i$.\
If the first Chern class \emph{vanishes}, this is the Fourier expansion of a 
Jacobi form of weight~$0$ and index equal to one half of the complex dimension of $X$ \cite{Ho-Diplom}.\ 
For automorphisms $g$ of $X$ one has the corresponding equivariant elliptic genus $\chi_y(g;q,{\cal L}X)$
where we let $\chi(g;X,E)=\sum_{q=0}^{d}(-1)^q \, \tr (g | H^q(X,{\cal O}(E)))$.

For a K3 surface $Y$, the elliptic genus $\chi_{-y}(q,{\cal L}Y)$ is the Jacobi form 
\begin{equation}\label{phi01}
2\,\phi_{0,1} (z; \tau)=  8\left(\left(\frac{\theta_{10}(z;\tau)}{\theta_{10}(0;\tau)}\right)^2+
\left(\frac{\theta_{00}(z;\tau)}{\theta_{00}(0;\tau)}\right)^2+
\left(\frac{\theta_{01}(z;\tau)}{\theta_{01}(0;\tau)}\right)^2\right)
\end{equation}
of weight~$0$ and index~$1$.
 
\smallskip

Eguchi, Ooguri and Tachikawa observed \cite{EOT} that 
the  decomposition of  $2\,\phi_{0,1}(z;\tau)$
into characters of an ${\cal N}\!=\! 4$ super algebra 
at central charge $c=6$ has multiplicities which are 
sums of dimensions of irreducible representations of  $M_{24}$.

The observation of Eguchi, Ooguri and Tachikawa suggests the existence of a graded
$M_{24}$-module $K=\bigoplus_{n=0}^{\infty} K_n\, q^{n-1/8}$ whose graded character is given by the decomposition of 
the elliptic genus into characters of the  ${\cal N}\!=\! 4$ super algebra.\   
Subsequently, analogues of McKay-Thompson series in monstrous moonshine were proposed in
several works~\cite{C, EH1, GHV, GHV2, CD}.\ 
The corresponding McKay-Thompson series for $g\in M_{24}$ are of the form
\begin{equation}\label{thompsonseries}
\Sigma_g(q)\ =\ q^{-1/8}\,\sum_{n=0}^{\infty} {\rm Tr}(g|K_n)\, q^{n} \ = \
\frac{e(g)}{24}\, \Sigma(q) - \frac{f_g(q)}{\eta(q)^3}.
\end{equation}
Here, $\Sigma=\Sigma_e$ is the graded dimension of $K$ (an explicit mock modular form), $e(g)$ is the character of
the $24$-dimensional permutation representation of $M_{24}$, $f_g$ is a certain explicit modular form of weight~$2$ on
a congruence subgroup $\Gamma_0(N_g)$, and $\eta$ is the Dedekind eta function.\
Gannon has shown~\cite{Gannon}  that these McKay-Thompson series 
indeed determine a graded $M_{24}$-module.

\smallskip

In~\cite{CrHo}, Creutzig and the first author have shown that for symplectic automorphisms of a K3 surface, 
the McKay-Thompson series of Mathieu Moonshine determines the equivariant elliptic genus:
\begin{thm}\label{moonshinegeometric}
Let $g$ be a finite symplectic automorphism of a K3 surface $Y$.\ Then
the equivariant elliptic genus and the character determined by the McKay-Thompson series of 
Mathieu moonshine agree, i.e.~one has
$$\chi_{-y}(g;q,{\cal L}Y)=\frac{e(g)}{12}\,\phi_{0,1}+f_g\,\phi_{-2,1},$$
where $g$ is considered on the right-hand-side as an element in $M_{24}$.
\end{thm}
Here, $\phi_{-2,1}$ is the Jacobi form
\begin{equation}\label{phi-21}
\phi_{-2,1} = y^{-1}(1-y)^2 \prod_{n=1}^\infty  \frac{(1-q^ny)^2(1-q^ny^{-1})^2}{(1-q^n)^4} 
\end{equation}
of weight $-2$ and index~$1$.


\subsection{The second quantized elliptic genus and its relation to Hilbert schemes of K3}

The elliptic genus of an orbifold $X/G$ for a finite group $G$
acting on a complex manifold $X$ is defined by 
$$\chi_y(q,{\cal L}(X/G)):=\frac{1}{|G|}\sum_{{\scriptsize g,h\in G\atop \scriptsize [g,h]=1}}\chi_y(g;q,{\cal L}_hX),$$
where ${\cal L}_h(X)$ is the $h$-twisted loop space and $\chi_y(g;q,{\cal L}_hX)$ is determined by
formally applying the equivariant Atiyah-Singer index theorem.

For a space $X$, let $\exp(pX):=\sum_{n=0}^{\infty}X^n/S_n\cdot p^n$ be the generating
series of its symmetric powers.\ It follows from calculations by Verlinde, Verlinde, Dijkgraaf and Moore~\cite{DMVV} that
the second quantized elliptic genus
$\chi_{-y}(q,{\cal L}\exp(pX))$ is, up to an automorphic correction factor, the Borcherds lift
of $\chi_{-y}(q,{\cal L}X)$.\ Explicitly one has
$$ \chi_{-y}(q,{\cal L}\exp(pX))=\prod_{n>0,\, m\geq 0,\,\ell} \exp\left(\sum_{k=1}^\infty \frac{1}{k} c(4nm-\ell^2) \,(p^n q^m y^\ell)^k \right), $$
where 
$\chi_{-y}(q,{\cal L}X)= \sum_{n,\,\ell\in\Z }  \,c(4n-\ell^2)\, q^n y^\ell  $.

For K3 surfaces, there is the following connection between the orbifold 
elliptic genus of symmetric powers and the Hilbert schemes conjectured by~\cite{DMVV}.
\begin{thm}[Borisov and Libgober~\cite{BL}]\label{thmBL}
Let $Y$ be a K3 surface.\ Then
$${\chi_{y}(q,{\cal L}\exp(pY))=\sum_{n=0}^\infty \chi_{y}(q,{\cal L}Y^{[n]})\,p^n. } $$ 
$\hfill\Box$

\end{thm}

If $g$ acts on $X$ then there is an induced action of $g$
on $X^n/S_n$ since the diagonal action of $g$ on $X^n$ commutes with the $S_n$-action.\
There is then the following equivariant generalization (cf.~\cite{Ho-ober,C,EH2}):
$$ \chi_{-y}(g;q,{\cal L}\exp(pX))= 
\prod_{n>0,\, m\geq 0,\,\ell} \exp\left(\sum_{k=1}^\infty \frac{1}{k} c_{g^k}(4nm-\ell^2) \,(p^n q^m y^\ell)^k \right), $$
where 
$  \chi_{-y}(g;q,{\cal L}X)= \sum_{n,\,\ell\in\Z }  \,c_g(4n-\ell^2) q^n y^\ell  $.

Alternatively, we may take this directly as the definition of the equivariant second  quantized elliptic genus.

For a K3 surface $Y$, we can use the McKay-Thompson series of Mathieu Moonshine to define
$ \chi_{-y}(g;q,{\cal L}Y)$ for all $g\in M_{24}$ by the formula in Theorem~\ref{moonshinegeometric}.\
Then the previous formula allows us to also define the equivariant second quantized elliptic genus
for $g\in M_{24}$.\ We prove:

\begin{thm}\label{moonshinegeometric2}
Let $g\in M_{24}$ be a finite symplectic automorphism of order $1$, $3$, $4$, $5$, $7$, $8$ or $11$
acting on a hyperk\"ahler manifold $X$ of type $\K32$.\
Then the equivariant elliptic genus $\chi_{-y}(g;q,{\cal L}X)$
and the coefficient of $p^2$ in the equivariant second quantized elliptic genus 
determined by the McKay-Thompson series of Mathieu moonshine agree. 
\end{thm}

There are $11$ classes $g\in M_{24}$ acting symplectically on a manifold of type $\K32$.\ 
If $g=1$, the theorem follows from Theorem~\ref{thmBL}.\
For $g$ acting by symplectic automorphisms on a K3 surface, there is 
probably an equivariant generalization, however this currently seems to be unknown.\
This would not, in any case, apply to the three classes of 
order $11$, $14$ and $15$ which do not correspond to symplectic automorphisms of K3.\
We have verified the result also for the four cases $2$, $6$, $14$ or $15$ for the first coefficients (up to the order $4$ in $q$).

\smallskip

To prove Theorem \ref{moonshinegeometric2}, we start by showing that both the equivariant elliptic genus 
and the $p^2$ coefficient of the equivariant second quantized elliptic
genus are weak Jacobi forms. 
\begin{prop}\label{elljacobi}
Let $N={\rm ord}(g)$ and assume $N>2$.\
Then  $\chi_{-y}(g;q,{\cal L}X)$ 
is a weak Jacobi form on $\Gamma_0(N)$ of weight $0$ and index $2$.
\end{prop}

\pf The proof is similar to that of Lemma~4.2 in~\cite{CrHo}.\ Because $g$ has order $\geq 3$,
it follows from Table~\ref{Fixpointset} that $X^g$ consists only of isolated fixed-points $\{p_i\}$.\
Set
\[
\varphi(u;\,\tau):=\vartheta_1(u;\,\tau)\eta(\tau)^{-3} = -i(y^{1/2}-y^{-1/2})\prod_{n=1}^\infty  (1-yq^n)(1-y^{-1}q^{n})(1-q^n)^{-2}.
\]
The fixed-point formula gives (cf.~also equation~(\ref{ASI}) in the case of $\chi_{-y}(g;X)$):
$$\chi_{-y}(g;q,{\cal L}X)=\sum_{p_i}
\frac{\varphi\bigl(u+\frac{n_i}{N};\,\tau\bigr)\varphi\bigl(u-\frac{n_i}{N};\,\tau\bigr)}
{\varphi\bigl(\frac{n_i}{N};\,\tau\bigr)\varphi\bigl(-\frac{n_i}{N};\,\tau\bigr)}\cdot
\frac{\varphi\bigl(u+\frac{m_i}{N};\,\tau\bigr)\varphi\bigl(u-\frac{m_i}{N};\,\tau\bigr)}
{\varphi\bigl(\frac{m_i}{N};\,\tau\bigr)\varphi\bigl(-\frac{m_i}{N};\,\tau\bigr)},
$$
where the pair $\{n_i,m_i\}$ for a given $p_i$ can be read off from Table~\ref{Fixpointset}.\ Namely,
$\{\zeta^{n_i},\zeta^{-n_i}, \zeta^{m_i}, \zeta^{-m_i}\}$ with $\zeta$ a primitive $N$-th root of unity 
are the eigenvalues of $g$ acting at the tangent space of $p_i$.

In order to check the Jacobi transformation property, we  consider the action of $(\Z/N\Z)^*$
(see the proof of Lemma~4.2 in~\cite{CrHo}).\
The fixed-point of type $\{\zeta^{n_i}, \zeta^{m_i}\}$ is mapped by $d\in (\Z/N\Z)^*$ to 
$\{(\zeta^{n_i})^d, (\zeta^{m_i})^d\}$, and it is clear from Table~\ref{Fixpointset} that this induces a permutation action of $ (\Z/N\Z)^*$
on $\{p_i\}$.\ Therefore, the above expression for $\chi_{-y}(g;q,{\cal L}X)$ is left invariant. $\hfill \Box$

\begin{prop}\label{Ramm}
Let $N={\rm ord}(g)$ and assume $N\not =6$, $14$, $15$.\ 
Then the coefficient of $p^2$ in the equivariant second quantized elliptic genus 
$\chi_{-y}(g;q,{\cal L}\exp(pX))$
is a weak Jacobi form for $\Gamma_0(N)$ of weight $0$ and index $2$.
\end{prop}

\pf
According to M.~Raum~\cite{Raum} (Theorem~1.2),
$$\Phi_g:=pqy\prod_{(n,m,\ell)>0} \exp\left(-\sum_{k=1}^\infty \frac{1}{k} c_{g^k}(4nm-\ell^2) \,(p^n q^m y^\ell)^k \right) $$
is a Siegel modular form of degree~$2$ of a certain weight $k_g$ on the subgroup $\Gamma_0^{(2)}(N)\subseteq Sp(4, \Z)$.\ 
Here, the product runs over triples of integers $(n,m,\ell)$ and
$(n,m,\ell)>0$ means that $n>0$, or $n=0$ and $m>0$, or $n=m=0$ and $\ell <0$.\
This implies that the coefficient $\psi_{g,n}$ of $p^n$ in the Fourier-Jacobi expansion 
$$\Phi_g^{-1} = \sum_{n=0}^\infty \psi_{g,n} \,  p^n$$
is a Jacobi form of weight $-k_g$ and index $m$ for $\Gamma_0(N)$.

To compare $\Phi_g^{-1}$ with the second quantized elliptic genus, we write
$$\Phi_g^{-1}= \chi_{-y}(g;q,{\cal L}\exp(pY)) \cdot \alpha_g^{-1}$$
so that
\begin{eqnarray*}
\alpha_g & =  & \prod_{(m,\ell)\geq 0} \exp\left(-\sum_{k=1}^\infty \frac{1}{k} c_{g^k}(-\ell^2) \,( q^m y^\ell)^k \right) \\
&& \prod_{(m,\ell)\geq 0} \prod_{d\mid N} \left(1-(q^m y^\ell)^d\right)^{ d^{-1} \sum_{e\mid d} \mu(d/e)c_{g^e}(-\ell^2)}.
\end{eqnarray*}
The only contributions come from terms $c_h(-\ell^2)$ for $\ell\in \{0,\,\pm 1\}$, which are determined by the
Hodge structure of the K3 surface.\ 
One has $c_h(-1)=2$ and $c_h(0)=e(h)-4$ where $e(h)$ is the number of fixed-points of
the element $h \in M_{24}$ acting on $24$ elements.

We claim that 
$$\alpha_g = p\,\eta_g\,\phi_{-2,1}$$
(cf.~\cite{C}, Section~4) where $\eta_g$ is the usual twisted eta-product for $g$.\
(If $g$ has cycle shape $a_1^{b_1}a_2^{b_2}\ldots$,
$\eta_g(q):=\eta(q^{a_1})^{b_1}\eta(q^{a_2})^{b_2}\ldots$.)\
Indeed, the contribution for $m>0$, $\ell=\pm 1$ and the constant term $-4$ of $c_h(0)$ gives the product formula
for $\phi_{-2,1}$ as in equation~(\ref{phi-21}) up to the leading $y$-factors.\ The contribution for $m>0$, $\ell=0$
is the twisted eta-product~$\eta_g$.\ Note that  $d^{-1} \sum_{e\mid d} \mu(d/e)e(g^e)$
counts the number of $d$-cycles of~$g$.\
The contribution for $m=0$ and $\ell=-1$ gives the missing $y$-factors, and the factor $p$ comes from the factor $qyp$ in front of $\Phi_g$.

\smallskip

Since the weight of $\eta_g$ is $k_g+2$~\cite{MaM24} and $\phi_{-2,1}$ has weight $-2$ and index~$1$, it follows that
the coefficient $\psi_{g,n-1} \cdot  \alpha_g$ of $p^n$ in $\chi_{-y}(g;q,{\cal L}\exp(pY))$ 
is a Jacobi form of weight~$0$ and index~$n$ for $\Gamma_0(N)$.  $\hfill \Box$

\smallskip
\noindent
{\bf Proof of Theorem~\ref{moonshinegeometric2}:} 
Let $J_{k,m}(\Gamma_0(N))$ be the space of weak holomorphic Jacobi forms for $\Gamma_0(N)$ of weight $k$ and index $m$.
It follows from~\cite{AI} Prop.~6.1
that
$$J_{0,2}(\Gamma_0(N))\cong M_0(\Gamma_0(N)) \times  M_2(\Gamma_0(N)) \times  M_4(\Gamma_0(N)),$$
where $M_l(\Gamma_0(N))$ is the space of holomorphic modular forms of weight $l$ for $\Gamma_0(N)$.\
The isomorphism is defined by 
$$(f_0,f_2,f_4)\mapsto f_0 \, \phi_{0,1}^2  + f_2 \, \phi_{0,1}\phi_{-2,1} + f_4\, \phi_{-2,1}^2$$
for $(f_0,f_2,f_4)\in  M_0(\Gamma_0(N)) \times  M_2(\Gamma_0(N)) \times  M_4(\Gamma_0(N))$.


By using explicit bases for the spaces of modular forms of weights $0$, $2$ and $4$ on $\Gamma_0(N)$,
 the result  follows  from
Proposition~\ref{elljacobi} and Proposition~\ref{Ramm} by checking the equality for sufficiently many
coefficients.\ This was carried through using Magma.  \phantom{a lot of spacexxx} $\hfill \Box$

\medskip

Our calculation shows that Theorem~\ref{moonshinegeometric2} will also hold for the remaining four classes $g$ of order
$2$, $6$, $14$ and $15$, once the Jacobi form property has been verified. See~\cite{PV} for work in this direction.

\medskip

\noindent{\bf Remark:} In general, a finite symplectic automorphism $g$ of a hyperk\"ahler manifold $X$ of type $K3^{[n]}$
defines a conjugacy class $[g]$ in ${\rm Co}_0$. We conjecture that for $g$ in $M_{24}$ 
the corresponding equivariant elliptic genus $\chi_{-y}(g;q,{\cal L}X)$ equals the coefficient of $p^n$ 
in the equivariant second quantized elliptic genus determined by the corresponding McKay-Thompson series 
of Mathieu moonshine.\
There may even be a generalization of the McKay-Thompson series of Mathieu moonshine for $g$ \emph{not} in $M_{24}$.
  

\appendix

\section{Tables of admissible groups and their lattices}\label{AppendixA}

The appendix contains three tables, listing all admissible subgroups of ${\rm Co}_0$ and information about
the coinvariant lattices $L_G$.

\medskip

Table~\ref{Gconclasses} lists the $198$ classes of admissible subgroups $G$ of ${\rm Co}_0$ sorted by size. 
The entry of the first two columns is clear.\ The third column lists the abstract isomorphism type of the group
either by name or the number of the small group library~\cite{BEO}.\ The fourth column has an entry ``K3'' if  
 $(L_G,G)$ is one of the $82$ group lattices of Table~10.2 of~\cite{Ha}, arising from a symplectic
action on a K3 surface; an ``M'' if the group is inside $M_{23}$ and realized; an  ``S'' if it is 
only realized by an ${\cal S}$-lattice group; and  ``$-$'' if the group is not realized.\ 
The fifth column records the dimension of $L_G$, the sixth column gives the number of the
unique largest group with the same lattice $L_G$, and the last column lists the numbers of all admissible groups
containing $G$ as a maximal subgroup.

\smallskip

Table~\ref{Lisoclasses} lists the $69$ isomorphism types of lattices $L_G$ for admissible groups $G$.\ The second column
gives the number (from Table~\ref{Gconclasses}) of the unique largest group with the same lattice $L_G$, while the next
three columns are self explanatory.\ Column ${\rm Det}$ gives the determinant of $L_G$, and column $A_{L_G}$ 
the structure of the discriminant group.\ In the Genus column we provide the genus symbol of $L_G$ as defined in~\cite{CoSl}
(we omit the signature since is directly determined by the rank).\ The last column provides information about the realization:\
if $L_G$ is one of the $41$ lattices of the table in Section~10.3 of~\cite{Ha}, 
the entry is ``K3 $\#n$'' where $n$ is the K3 group number;
if $L_G$ is one of the $13$ lattices belonging to the groups in Theorem~\ref{thm57911}, the entry is ``max $\#n$'' where
$n$ is the $n$-th group of that theorem; if $L_G$ belongs to (maximal) ${\cal S}$-groups, the entry is ''(max) \S-lattice'';
if $G$ can be realized in $O(L)$ the entry is ``$O(L)$''; in the remaining cases there is a ``$-$''.

\smallskip

Finally, in Table~\ref{Lgroups} we list all lattices $L_G$ for which there are several admissible groups 
with the same lattice $L_G$.\ The first column gives the number of the lattice as in Table~\ref{Lisoclasses}.\
The second column gives the number of groups with the same lattice $L_G$.\ The last column lists all such groups by
their number as in Table~\ref{Gconclasses}.
 

\small
\begin{longtable}{rrrclrl}
\caption{Conjugacy classes of admissible groups}\label{Gconclasses}  \\
\mbox{No.} & \mbox{order} & $G$ & \mbox{Type} & \mbox{dim} & \mbox{fix}  &\mbox{minimal overgroups}    \\ \hline
\endfirsthead 
\caption[]{Conjugacy classes of admissible groups}\\
\mbox{No.} & \mbox{order} & $G$ & \mbox{Type} & \mbox{dim} & \mbox{fix}  &\mbox{minimal overgroups}    \\ \hline
\endhead
1  &  1  &  $\# 1 $  &  K3  &  0  &  1  &  \{ 2, 3, 4, 7, 11, 23 \} \\ \hline
2  &  2  &  $\# 1 $  &  K3  &  8  &  2  &  \{ 5, 6, 8, 9, 10, 22, 31 \} \\ \hline
3  &  3  &  $\# 1 $  &  K3  &  12  &  3  &  \{ 8, 9, 18, 19, 21, 25, 32, 54 \} \\
4  &  3  &  $\# 1 $  &   S  &  18  &  170  &  \{ 10, 18, 20, 21 \} \\ \hline
5  &  4  &  $\# 1 $  &  K3  &  14  &  5  &  \{ 12, 13, 14, 15, 17, 26, 28, 30, 53, 79 \} \\
6  &  4  &  $\# 2 $  &  K3  &  12  &  6  &  \{ 13, 14, 16, 24, 25, 27, 29 \} \\  \hline
7  &  5  &  $\# 1 $  &  K3  &  16  &  22  &  \{ 22, 32, 99, 114, 169 \} \\  \hline
8  &  6  &  $\# 2 $  &  K3  &  16  &  24  &  \{ 24, 27, 29, 30, 48, 51, 52, 55, 57, 58, 67, 83 \} \\
9  &  6  &  $\# 1 $  &  K3  &  14  &  9  &  \{ 24, 48, 49, 50, 51, 56, 102 \} \\
10  &  6  &  $\# 2 $  &  S  &  18  &  170  &  \{ 26, 28, 49, 52 \} \\  \hline
11  &  7  &  $\# 1 $  &  K3  &  18  &  54  &  \{ 31, 54, 100 \} \\  \hline
12  &  8  &  $\# 1 $  &  K3  &  18  &  40  &  \{ 33, 40, 44, 45, 109 \} \\
13  &  8  &  $\# 3 $  &  K3  &  15  &  13  &  \{ 40, 46, 47, 56, 60, 111 \} \\
14  &  8  &  $\# 2 $  &  K3  &  16  &  47  &  \{ 33, 34, 36, 38, 39, 41, 42, 43, 44, 46, 47 \} \\
15  &  8  &  $\# 4 $  &  K3  &  17  &  68  &  \{ 42, 43, 45, 46, 57 \} \\
16  &  8  &  $\# 5 $  &  K3  &  14  &  16  &  \{ 35, 36, 37, 41, 47, 58, 100 \} \\
17  &  8  &  $\# 4 $  &  K3  &  17  &  17  &  \{ 39, 40, 45, 55, 59, 110 \} \\  \hline
18  &  9  &  $\# 2 $  &   S  &  18  &  170  &  \{ 49, 61, 62, 63, 64, 65, 66 \} \\
19  &  9  &  $\# 2 $  &  K3  &  16  &  50  &  \{ 48, 50, 62, 66, 78 \} \\
20  &  9  &  $\# 1 $  &   S &  20  &  198  &  \{ 63, 64, 65 \} \\
21  &  9  &  $\# 2 $  &   S &  20  &  198  &  \{ 51, 52, 62 \} \\  \hline
22  &  10  &  $\# 1 $  &  K3  &  16  &  22  &  \{ 53, 67, 102, 136, 181 \} \\  \hline
23  &  11  &  $\# 1 $  &  M  &  20  &  177  &  \{ 99 \} \\  \hline
24  &  12  &  $\# 4 $  &  K3  &  16  &  24  &  \{ 60, 80, 81, 86, 88, 132, 177 \} \\
25  &  12  &  $\# 3 $  &  K3  &  16  &  25  &  \{ 56, 58, 78, 85, 87, 102, 160 \} \\
26  &  12  &  $\# 1 $  &  S  &  19  &  185  &  \{ 59, 82, 130 \} \\
27  &  12  &  $\# 5 $  &  $-$  &  18  &  27  &  \{ 90, 93 \} \\
28  &  12  &  $\# 2 $  &   S &  20  &  191  &  \{ 59, 127, 131 \} \\
29  &  12  &  $\# 5 $  &  K3  &  18  &  112  &  \{ 60, 78, 80, 89, 92 \} \\
30  &  12  &  $\# 1 $  &  K3  &  18  &  112  &  \{ 60, 82, 84, 91, 94, 101, 129 \} \\  \hline
31  &  14  &  $\# 2 $  &  M  &  20  &  162  &  \{ 83 \} \\  \hline
32  &  15  &  $\# 1 $  &  M  &  20  &  164  &  \{ 67 \} \\  \hline
33  &  16  &  $\# 6 $  &  K3  &  19  &  168  &  \{ 69, 70, 71, 75 \} \\
34  &  16  &  $\# 2 $  &  $-$  &  18  &  34  &  \{\} \\
35  &  16  &  $\# 14 $  &  $-$  &  15  &  35  &  \{ 93 \} \\
36  &  16  &  $\# 10 $  &  $-$  &  17  &  36  &  \{\} \\
37  &  16  &  $\# 14 $  &  K3  &  15  &  37  &  \{ 77, 85, 92, 114 \} \\
38  &  16  &  $\# 2 $  &  K3  &  18  &  150  &  \{ 69, 72, 87 \} \\
39  &  16  &  $\# 12 $  &  K3  &  18  &  150  &  \{ 71, 72, 73, 89 \} \\
40  &  16  &  $\# 8 $  &  K3  &  18  &  40  &  \{ 71, 86, 135, 178 \} \\
41  &  16  &  $\# 3 $  &  K3  &  17  &  77  &  \{ 72, 76, 77, 84 \} \\
42  &  16  &  $\# 12 $  & $-$   &  18  &  42  &  \{ 90 \} \\
43  &  16  &  $\# 12 $  &  $-$  &  18  &  74  &  \{ 70, 74 \} \\
44  &  16  &  $\# 6 $  &  $-$  &  19  &  44  &  \{ 73 \} \\
45  &  16  &  $\# 9 $  &  K3  &  19  &  168  &  \{ 71, 91, 94 \} \\
46  &  16  &  $\# 13 $  &  K3  &  17  &  68  &  \{ 68, 69, 71, 74 \} \\
47  &  16  &  $\# 11 $  &  K3  &  16  &  47  &  \{ 68, 72, 75, 76, 77, 88 \} \\  \hline
48  &  18  &  $\# 3 $  &  K3  &  18  &  81  &  \{ 81, 95, 98, 142, 156 \} \\
49  &  18  &  $\# 3 $  &   S &  18  &  170  &  \{ 95, 96, 97 \} \\
50  &  18  &  $\# 4 $  &  K3  &  16  &  50  &  \{ 79, 81, 96, 98, 112 \} \\
51  &  18  &  $\# 3 $  &   S &  20  &  198  &  \{ 80, 95, 98 \} \\
52  &  18  &  $\# 5 $  &   S &  20  &  198  &  \{ 80, 82, 95 \} \\  \hline
53  &  20  &  $\# 3 $  &  K3  &  18  &  53  &  \{ 101, 132, 158, 178 \} \\  \hline
54  &  21  &  $\# 1 $  &  K3  &  18  &  54  &  \{ 83, 140, 141 \} \\  \hline
55  &  24  &  $\# 3 $  &  K3  &  19  &  86  &  \{ 86, 89, 94 \} \\
56  &  24  &  $\# 12 $  &  K3  &  17  &  56  &  \{ 88, 112, 120, 124, 132, 140, 163, 176 \} \\
57  &  24  &  $\# 3 $  &  K3  &  19  &  143  &  \{ 90, 91, 122 \} \\
58  &  24  &  $\# 13 $  &  K3  &  18  &  88  &  \{ 84, 88, 92, 93, 122, 125, 141 \} \\
59  &  24  &  $\# 4 $  &  S  &  20  &  191  &  \{ 151 \} \\
60  &  24  &  $\# 8 $  &  K3  &  18  &  112  &  \{ 112, 113, 121, 123, 152 \} \\  \hline
61  &  27  &  $\# 3 $  &  S  &  18  &  170  &  \{ 97, 116, 117, 119 \} \\
62  &  27  &  $\# 5 $  &  S  &  20  &  198  &  \{ 95, 98, 117, 118 \} \\
63  &  27  &  $\# 4 $  &  S  &  20  &  198  &  \{ 115, 116 \} \\
64  &  27  &  $\# 4 $  &  S  &  20  &  198  &  \{ 115, 117 \} \\
65  &  27  &  $\# 2 $  &  S  &  20  &  198  &  \{ 115, 116 \} \\
66  &  27  &  $\# 5 $  &   S &  18  &  170  &  \{ 96, 115, 118, 119 \} \\  \hline
67  &  30  &  $\# 2 $  &  M  &  20  &  164  &  \{ 101, 142 \} \\  \hline
68  &  32  &  $\# 49 $  &  K3  &  17  &  68  &  \{ 106, 108, 122 \} \\
69  &  32  &  $\# 11 $  &  K3  &  19  &  168  &  \{ 106, 124 \} \\
70  &  32  &  $\# 8 $  & $-$   &  20  &  70  &  \{\} \\
71  &  32  &  $\# 44 $  &  K3  &  19  &  168  &  \{ 106, 123 \} \\
72  &  32  &  $\# 31 $  &  K3  &  18  &  150  &  \{ 103, 105, 106, 107 \} \\
73  &  32  &  $\# 8 $  &  $-$  &  20  &  107  &  \{ 107 \} \\
74  &  32  &  $\# 50 $  &  $-$  &  18  &  74  &  \{\} \\
75  &  32  &  $\# 7 $  &  K3  &  19  &  168  &  \{ 104, 106 \} \\
76  &  32  &  $\# 6 $  &  K3  &  18  &  108  &  \{ 104, 105, 108 \} \\
77  &  32  &  $\# 27 $  &  K3  &  17  &  77  &  \{ 103, 104, 108, 120, 121, 125, 136 \} \\  \hline
78  &  36  &  $\# 11 $  &  K3  &  18  &  112  &  \{ 112, 134, 142, 184 \} \\
79  &  36  &  $\# 9 $  &  K3  &  18  &  79  &  \{ 109, 110, 111, 129, 130, 131, 163, 172 \} \\
80  &  36  &  $\# 12 $  &  S  &  20  &  198  &  \{ 113, 126 \} \\
81  &  36  &  $\# 10 $  &  K3  &  18  &  81  &  \{ 111, 126, 128, 164, 173 \} \\
82  &  36  &  $\# 7 $  &  S  &  20  &  198  &  \{ 113, 161 \} \\  \hline
83  &  42  &  $\# 2 $  &  M  &  20  &  162  &  \{ 162 \} \\  \hline
84  &  48  &  $\# 30 $  &  K3  &  19  &  157  &  \{ 121, 147, 149 \} \\
85  &  48  &  $\# 50 $  &  K3  &  17  &  85  &  \{ 120, 125, 134, 150 \} \\
86  &  48  &  $\# 29 $  &  K3  &  19  &  86  &  \{ 123 \} \\
87  &  48  &  $\# 3 $  &  K3  &  18  &  150  &  \{ 124, 150 \} \\
88  &  48  &  $\# 48 $  &  K3  &  18  &  88  &  \{ 121, 143, 148, 162, 179 \} \\
89  &  48  &  $\# 32 $  &  M  &  20  &  187  &  \{ 123, 144 \} \\
90  &  48  &  $\# 32 $  &  $-$  &  20  &  90  &  \{\} \\
91  &  48  &  $\# 28 $  &  M  &  20  &  196  &  \{ 149 \} \\
92  &  48  &  $\# 49 $  &  K3  &  19  &  157  &  \{ 121, 134, 144, 145, 146 \} \\
93  &  48  &  $\# 49 $  &  $-$  &  19  &  93  &  \{\} \\
94  &  48  &  $\# 28 $  & M   &  20  &  187  &  \{ 123 \} \\  \hline
95  &  54  &  $\# 12 $  &  S  &  20  &  198  &  \{ 126, 137, 139 \} \\
96  &  54  &  $\# 13 $  &  S  &  18  &  170  &  \{ 130, 131, 138, 139 \} \\
97  &  54  &  $\# 8 $  &   S &  18  &  170  &  \{ 127, 137, 138 \} \\
98  &  54  &  $\# 13 $  & S   &  20  &  198  &  \{ 126, 128, 129, 139 \} \\  \hline
99  &  55  &  $\# 1 $  & M   &  20  &  177  &  \{ 177 \} \\  \hline
100  &  56  &  $\# 11 $  & M   &  20  &  188  &  \{ 141 \} \\  \hline
101  &  60  &  $\# 7 $  &  M  &  20  &  164  &  \{ 164 \} \\
102  &  60  &  $\# 5 $  &  K3  &  18  &  102  &  \{ 132, 142, 163, 177, 182, 183, 195 \} \\  \hline
103  &  64  &  $\# 242 $  &  K3  &  18  &  150  &  \{ 133, 144, 146, 150 \} \\
104  &  64  &  $\# 32 $  &  K3  &  19  &  168  &  \{ 133, 147, 158, 172 \} \\
105  &  64  &  $\# 35 $  &  K3  &  19  &  168  &  \{ 133 \} \\
106  &  64  &  $\# 136 $  &  K3  &  19  &  168  &  \{ 133, 149 \} \\
107  &  64  &  $\# 36 $  &  $-$  &  20  &  107  &  \{\} \\
108  &  64  &  $\# 138 $  &  K3  &  18  &  108  &  \{ 133, 143, 145, 148 \} \\  \hline
109  &  72  &  $\# 39 $  &  M  &  20  &  135  &  \{ 135 \} \\
110  &  72  &  $\# 41 $  &  K3  &  19  &  110  &  \{ 135, 151, 178, 197 \} \\
111  &  72  &  $\# 40 $  &  K3  &  19  &  111  &  \{ 135, 152, 179 \} \\
112  &  72  &  $\# 43 $  &  K3  &  18  &  112  &  \{ 157, 164, 192, 193 \} \\
113  &  72  &  $\# 22 $  &  S  &  20  &  198  &  \{ 175 \} \\  \hline
114  &  80  &  $\# 49 $  &  K3  &  19  &  183  &  \{ 136 \} \\  \hline
115  &  81  &  $\# 13 $  & S   &  20  &  198  &  \{ 153, 155 \} \\
116  &  81  &  $\# 8 $  &  S  &  20  &  198  &  \{ 153 \} \\
117  &  81  &  $\# 7 $  &  S  &  20  &  198  &  \{ 137, 155, 160 \} \\
118  &  81  &  $\# 15 $  & S   &  20  &  198  &  \{ 139, 155, 169 \} \\
119  &  81  &  $\# 12 $  &  S  &  18  &  170  &  \{ 138, 153, 154, 155 \} \\  \hline
120  &  96  &  $\# 227 $  &  K3  &  18  &  120  &  \{ 148, 156, 157, 168, 183 \} \\
121  &  96  &  $\# 195 $  &  K3  &  19  &  157  &  \{ 157, 165, 166, 167, 182 \} \\
122  &  96  &  $\# 204 $  &  K3  &  19  &  143  &  \{ 143, 145, 149 \} \\
123  &  96  &  $\# 190 $  & M   &  20  &  187  &  \{ 167 \} \\
124  &  96  &  $\# 64 $  &  K3  &  19  &  168  &  \{ 168 \} \\
125  &  96  &  $\# 70 $  &  K3  &  19  &  148  &  \{ 146, 147, 148, 156 \} \\  \hline
126  &  108  &  $\# 38 $  &  S  &  20  &  198  &  \{ 152, 159 \} \\
127  &  108  &  $\# 15 $  &  S  &  20  &  191  &  \{ 186 \} \\
128  &  108  &  $\# 40 $  &  S  &  20  &  198  &  \{ 152, 159, 160 \} \\
129  &  108  &  $\# 37 $  &  S  &  20  &  198  &  \{ 152, 161 \} \\
130  &  108  &  $\# 37 $  &  S  &  19  &  185  &  \{ 151, 161, 185 \} \\
131  &  108  &  $\# 36 $  &  S  &  20  &  191  &  \{ 151, 186 \} \\  \hline
132  &  120  &  $\# 34 $  &  K3  &  19  &  132  &  \{ 164, 179, 190, 193 \} \\  \hline
133  &  128  &  $\# 931 $  &  K3  &  19  &  168  &  \{ 165, 166, 167, 168 \} \\  \hline
134  &  144  &  $\# 184 $  &  K3  &  19  &  157  &  \{ 156, 157, 174 \} \\
135  &  144  &  $  \Z_3^2{:}QD_{16} $  & M   &  20  &  135  &  \{\} \\  \hline
136  &  160  &  $\# 234 $  &  K3  &  19  &  183  &  \{ 158, 182, 183 \} \\  \hline
137  &  162  &  $\# 10 $  &  S  &  20  &  198  &  \{ 171, 176 \} \\
138  &  162  &  $\# 46 $  &  S  &  18  &  170  &  \{ 170, 171 \} \\
139  &  162  &  $\# 52 $  &  S  &  20  &  198  &  \{ 159, 161, 171, 181 \} \\  \hline
140  &  168  &  $\# 42 $  &  K3  &  19  &  140  &  \{ 162, 188, 193, 197 \} \\
141  &  168  &  $\# 43 $  & M   &  20  &  188  &  \{ 188 \} \\  \hline
142  &  180  &  $\# 19 $  &  M  &  20  &  164  &  \{ 164 \} \\  \hline
143  &  192  &  $\# 1493 $  &  K3  &  19  &  143  &  \{ 165, 188 \} \\
144  &  192  &  $\# 1024 $  &  M  &  20  &  187  &  \{ 167, 174 \} \\
145  &  192  &  $\# 201 $  &  M  &  20  &  196  &  \{ 165, 182 \} \\
146  &  192  &  $\# 1009 $  & M   &  20  &  187  &  \{ 166, 174 \} \\
147  &  192  &  $\# 184 $  & M   &  20  &  187  &  \{ 166 \} \\
148  &  192  &  $\# 955 $  &  K3  &  19  &  148  &  \{ 166, 173, 188, 190 \} \\
149  &  192  &  $\# 1492 $  &  M  &  20  &  196  &  \{ 165 \} \\
150  &  192  &  $\# 1023 $  &  K3  &  18  &  150  &  \{ 168, 174, 183 \} \\  \hline
151  &  216  &  $\# 161 $  & S   &  20  &  191  &  \{ 191 \} \\
152  &  216  &  $\# 158 $  & S   &  20  &  198  &  \{ 175, 176 \} \\  \hline
153  &  243  &  $\# 57 $  &  S  &  20  &  198  &  \{ 180 \} \\
154  &  243  &  $\# 65 $  &  S  &  18  &  170  &  \{ 170, 180 \} \\
155  &  243  &  $\# 51 $  &  S  &  20  &  198  &  \{ 171, 180, 184 \} \\  \hline
156  &  288  &  $\# 1025 $  & M   &  20  &  173  &  \{ 173 \} \\
157  &  288  &  $\# 1026 $  &  K3  &  19  &  157  &  \{ 172, 173, 187 \} \\  \hline
158  &  320  &  $\# 1635 $  &  M  &  20  &  190  &  \{ 190 \} \\  \hline
159  &  324  &  $\# 167 $  & S   &  20  &  198  &  \{ 175, 184 \} \\
160  &  324  &  $\# 160 $  & S   &  20  &  198  &  \{ 176, 184 \} \\
161  &  324  &  $\# 163 $  & S   &  20  &  198  &  \{ 175, 194 \} \\  \hline
162  &  336  &  $ \Z_2\times L_2(7) $  & M   &  20  &  162  &  \{\} \\  \hline
163  &  360  &  $\# 118 $  &  K3  &  19  &  163  &  \{ 178, 179, 193, 196, 197, 198 \} \\  
164  &  360  &  $ (\Z_3\times A_5){:}\Z_2$  &  M  &  20  &  164  &  \{\} \\  \hline
165  &  384  &  $\# 5603 $  &  M  &  20  &  196  &  \{ 196 \} \\
166  &  384  &  $\# 5678 $  &  M  &  20  &  187  &  \{ 187 \} \\
167  &  384  &  $\# 18133 $  &  M  &  20  &  187  &  \{ 187 \} \\ 
168  &  384  &  $\# 18135 $  &  K3  &  19  &  168  &  \{ 187, 190, 196 \} \\  \hline
169  &  405  &  $\# 15 $  &  S  &  20  &  198  &  \{ 181 \} \\  \hline
170  &  486  &  $\# 249 $  &  S  &  18  &  170  &  \{ 185, 186, 189 \} \\
171  &  486  &  $\# 166 $  &  S  &  20  &  198  &  \{ 189, 192, 195 \} \\  \hline
172  &  576  &  $\# 8652 $  & M   &  20  &  196  &  \{ 196 \} \\
173  &  576  &  $ \Z_2^4{:}(S_3\times S_3)$  &  M  &  20  &  173  &  \{\} \\
174  &  576  &  $\# 5129 $  & M   &  20  &  187  &  \{ 187 \} \\  \hline
175  &  648  &  $\# 722 $  & S   &  20  &  198  &  \{ 192 \} \\
176  &  648  &  $\# 704 $  & S   &  20  &  198  &  \{ 192 \} \\  \hline
177  &  660  &  $ L_2(11) $  & M   &  20  &  177  &  \{\} \\  \hline
178  &  720  &  $ M_{10} $  & M   &  20  &  178  &  \{\} \\
179  &  720  &  $ S_6 $ &  M  &  20  &  179  &  \{\} \\  \hline
180  &  729  &  $\# 321 $  &  S  &  20  &  198  &  \{ 189 \} \\  \hline
181  &  810  &  $\# 101 $  &  S  &  20  &  198  &  \{ 195 \} \\  \hline
182  &  960  &  $\# 11358 $  &  M  &  20  &  196  &  \{ 196 \} \\
183  &  960  &  $\# 11357 $  &  K3  &  19  &  183  &  \{ 190, 196, 197 \} \\  \hline
184  &  972  &  $\# 877 $  & S   &  20  &  198  &  \{ 192, 195 \} \\
185  &  972  &  $\# 776 $  &  S  &  19  &  185  &  \{ 191, 194 \} \\
186  &  972  &  $\# 777 $  &  S  &  20  &  191  &  \{ 191 \} \\  \hline
187  &  1152  &  $ 2^6(\Z_3^2{:}\Z_2)$  &  M  &  20  &  187  &  \{\} \\  \hline
188  &  1344  &  $ \Z_2^3{:}L_2(7)$  &  M  &  20  &  188  &  \{\} \\  \hline
189  &  1458  &  $\# 1229 $  &  S  &  20  &  198  &  \{ 194 \} \\ \hline
190  &  1920  &  $\Z_2^4{:}S_5$  &  M  &  20  &  190  &  \{\} \\ \hline
191  &  1944  &  $3^{1+4}{:}2.2^2 $  &  S  &  20  &  191  &  \{\} \\
192  &  1944  &  $\# 3877 $  &  S  &  20  &  198  &  \{ 198 \} \\ \hline
193  &  2520  &   $A_7$  &  M  &  20  &  193  &  \{\} \\ \hline
194  &  2916  &  $3^4{:}(3^2{:}\Z_4)$ &  S  &  20  &  198  &  \{ 198 \} \\ \hline
195  &  4860  &  $3^4{:}A_5$  & S   &  20  &  198  &  \{ 198 \} \\ \hline
196  &  5760  &  $\Z_2^4{:}A_6$  & M   &  20  &  196  &  \{\} \\ \hline
197  &  20160  & $L_3(4)$ &  M  &  20  &  197  &  \{\} \\ \hline
198  &  29160  & $3^4{:}A_6$  &  S  &  20  &  198  &  \{\} \\ \hline
\end{longtable}

\newpage

\small
\begin{longtable}{rrrcrrlll}
\caption{Isometry types of coinvariant lattices $L_G$}\label{Lisoclasses}\\
\mbox{No.} & $G$-No. & $|G|$ & Symbol &$\!\!\!\!$ Rank & Det & $A_{L_G}$  & Genus & Type  \\ \hline
\endfirsthead
\caption[]{Isometry types of coinvariant lattices $L_G$}\\
\mbox{No.} & $G$-No. & $|G|$ & Symbol &$\!\!\!\!$ Rank & Det & $A_{L_G}$  & Genus & Type  \\ \hline
\endhead
1 & 1 & 1 & 1 & 0 & 1 & $1$ & 1 & K3 \# 0   \\ \hline
2 & 2 & 2 &   $\Z_2$ & 8 & 256 & $2^8$ & $2_{\rm I\!I}^{+8}$ & K3 \# 1    \\ \hline
3 & 3 & 3 & $\Z_3$ & 12 & 729 & $3^6$ & $3^{+6}$ & K3 \# 2   \\ \hline
4 & 5 & 4 & $\Z_4$ & 14 & 1024 & $2^24^4$ & $2_{2}^{+2}4_{\rm I\!I}^{+4}$ & K3 \# 4   \\
5 & 6 & 4 & $\Z_2^2$ & 12 & 1024 & $2^64^2$ & $2_{\rm I\!I}^{-6}4_{\rm I\!I}^{-2}$ & K3 \# 3   \\ \hline
6 & 9 & 6 & $\Z_6$ & 14 & 972 & $3^36^2$ & $2_{\rm I\!I}^{-2}3^{+5}$ & K3 \# 6   \\  \hline
7 & 13 & 8 & $D_8$ & 15 & -1024 & $4^5$ & $4_{1}^{+5}$ & K3 \# 15   \\
8 & 16 & 8 & $\Z_2^3$ & 14 & 1024 & $2^64^2$ & $2_{\rm I\!I}^{+6}4_{2}^{+2}$ & K3 \# 14   \\
9 & 17 & 8 & $Q_8$ & 17 & -512 & $2^38^2$ & $2_{3}^{+3}8_{\rm I\!I}^{-2}$ & K3 \# 12   \\  \hline
10 & 22 & 10 & $D_{10}$ & 16 & 625 & $5^4$ & $5^{+4}$ & K3 \# 16   \\ \hline
11 & 24 & 12 & $D_{12}$ & 16 & 1296 & $6^4$ & $2_{\rm I\!I}^{+4}3^{+4}$ & K3 \# 18   \\
12 & 25 & 12 & $A_4$ & 16 & 576 & $2^212^2$ & $2_{\rm I\!I}^{-2}4_{\rm I\!I}^{-2}3^{+2}$ & K3 \# 17   \\
13 & 27 & 12 & $\Z_2\times \Z_6$ & 18 & 1728 & $2^36^3$ & $2_{\rm I\!I}^{-6}3^{+3}$ &  $-$  \\  \hline
14 & 34 & 16 & $\Z_4^2$ & 18 & 1024 & $2^24^4$ & $2_{6}^{+2}4_{\rm I\!I}^{+4}$ & $O(L)$    \\
15 & 35 & 16 & $\Z_2^4$ & 15 & -1024 & $2^84^1$ & $2_{\rm I\!I}^{+8}4_{1}^{+1}$ &  $-$     \\
16 & 36 & 16 & $\Z_2^2\times \Z_4$ & 17 & -1024 & $2^44^3$ & $2_{\rm I\!I}^{+4}4_{7}^{+3}$ &  $-$    \\
17 & 37 & 16 & $\Z_2^4$ & 15 & -512 & $2^68^1$ & $2_{\rm I\!I}^{+6}8_{1}^{+1}$ & K3 \# 21   \\
18 & 40 & 16 & $\Gamma_3a_2$ & 18 & 512 & $2^14^18^2$ & $2_{5}^{+1}4_{1}^{+1}8_{\rm I\!I}^{+2}$ & K3 \# 26   \\
19 & 42 & 16 & $\Gamma_2a_2$ & 18 & 1024 & $2^24^4$ & $2_{\rm I\!I}^{+2}4_{6}^{+4}$ &  $-$   \\
20 & 44 & 16 & $\Gamma_2d$ & 19 & -512 & $2^38^2$ & $2_{5}^{+3}8_{\rm I\!I}^{+2}$ &   $O(L)$   \\
21 & 47 & 16 & $\Gamma_2a_1$ & 16 & 1024 & $2^24^4$ & $2_{\rm I\!I}^{+2}4_{0}^{+4}$ & K3 \# 22   \\  \hline
22 & 50 & 18 & $A_{3,3}$ & 16 & 729 & $3^49^1$ & $3^{+4}9^{-1}$ & K3 \# 30  \\  \hline
23 & 53 & 20 & ${\rm Hol}(\Z_4)$ & 18 & 500 & $5^110^2$ & $2_{2}^{+2}5^{+3}$ & K3 \# 32  \\  \hline
24 & 54 & 21 & $\Z_7{:}\Z_3$ & 18 & 343 & $7^3$ & $7^{+3}$ & K3 \# 33   \\  \hline
25 & 56 & 24 & $S_4$ & 17 & -576 & $4^112^2$ & $4_{3}^{+3}3^{+2}$ & K3 \#  34  \\  \hline
26 & 68 & 32 & $\Gamma_5a_1$ & 17 & -1024 & $4^5$ & $4_{7}^{+5}$ & K3 \# 40   \\
27 & 70 & 32 & $\Gamma_7a_3$ & 20 & 256 & $2^28^2$ & $2_{6}^{+2}8_{6}^{+2}$ & $O(L)$ \\
28 & 74 & 32 & $\Gamma_5a_2$ & 18 & 1024 & $2^24^4$ & $2_{6}^{+2}4_{\rm I\!I}^{+4}$ & $O(L)$    \\
29 & 77 & 32 & $\Gamma_4a_1$ & 17 & -512 & $2^24^28^1$ & $2_{\rm I\!I}^{+2}4_{6}^{+2}8_{1}^{+1}$ & K3 \# 39   \\  \hline
30 & 79 & 36 & $3^2.\Z_4$ & 18 & 324 & $3^16^118^1$& $2_{2}^{+2}3^{+2}9^{-1}$ & K3 \# 46   \\
31 & 81 & 36 & $S_{3,3}$ & 18 & 972 & $3^26^118^1$ & $2_{\rm I\!I}^{-2}3^{+3}9^{-1}$ & K3 \# 48   \\  \hline
32 & 85 & 48 & $2^4{:}\Z_3$ & 17 & -384 & $2^424^1$ & $2_{\rm I\!I}^{-4}8_{1}^{+1}3^{-1}$ & K3 \# 49   \\
33 & 86 & 48 & $T_{48}$ & 19 & -384 & $2^18^124^1$ & $2_{7}^{+1}8_{\rm I\!I}^{-2}3^{-1}$ & K3 \# 54    \\
34 & 88 & 48 & $\Z_2\times S_4$ & 18 & 576 & $2^212^2$ & $2_{\rm I\!I}^{+2}4_{2}^{+2}3^{+2}$ & K3 \# 51  \\
35 & 90 & 48 & $\Z_2 \times {\rm SL}_2(3)$ & 20 & 192 & $2^24^112^1$ & $2_{\rm I\!I}^{-2}4_{2}^{+2}3^{+1}$ &  $-$   \\
36 & 93 & 48 & $\Z_2^2 \times A_4$ & 19 & -576 & $2^36^112^1$ & $2_{\rm I\!I}^{+4}4_{1}^{+1}3^{+2}$ &  $-$   \\  \hline
37 & 102 & 60 & $A_5$ & 18 & 300 & $10^130^1$ & $2_{\rm I\!I}^{-2}3^{+1}5^{-2}$ & K3 \#  55  \\  \hline
38 & 107 & 64 & $\Gamma_{23}a_3$ & 20 & 128 & $2^316^1$ & $2_{5}^{+3}16_{7}^{+1}$ &  $O(L)$  \\
39 & 108 & 64 & $\Gamma_{25}a_1$ & 18 & 512 & $4^38^1$ & $4_{5}^{+3}8_{1}^{+1}$ & K3 \#  56  \\  \hline
40 & 110 & 72 & $M_9$ & 19 & -216 & $2^16^118^1$ & $2_{3}^{+3}3^{-1}9^{-1}$ & K3 \# 63   \\
41 & 111 & 72 & $N_{72}$ & 19 & -324 & $3^236^1$ & $4_{1}^{+1}3^{+2}9^{-1}$ & K3 \# 62   \\ 
42 & 112 & 72 & $A_{4,3}$ & 18 & 432 & $3^112^2$ & $4_{\rm I\!I}^{-2}3^{-3}$ & K3 \# 61   \\  \hline
43 & 120 & 96 & $2^4.D_6$ & 18 & 384 & $2^24^124^1$ & $2_{\rm I\!I}^{-2}4_{7}^{+1}8_{1}^{+1}3^{-1}$ & K3 \# 65   \\  \hline
44 & 132 & 120 & $S_5$ & 19 & -300 & $5^160^1$ & $4_{3}^{-1}3^{+1}5^{-2}$ & K3 \# 70   \\  \hline
45 & 135 & 144 & $\Z_3^2{:}QD_{16}$ & 20 & 216 & $6^136^1$ & $2_{1}^{+1}4_{1}^{+1}3^{-1}9^{-1}$ & max \# 13   \\  \hline
46 & 140 & 168 & $L_2(7)$ & 19 & -196 & $7^128^1$ & $4_{1}^{+1}7^{+2}$ & K3 \# 74   \\  \hline
47 & 143 & 192 & $T_{192}$ & 19 & -192 & $4^212^1$ & $4_{7}^{-3}3^{+1}$ & K3 \# 77   \\
48 & 148 & 192 & $H_{192}$ & 19 & -384 & $4^224^1$ & $4_{2}^{-2}8_{1}^{+1}3^{-1}$ & K3 \# 76   \\ 
49 & 150 & 192 & $4^2.A_4$ & 18 & 256 & $2^28^2$ & $2_{\rm I\!I}^{-2}8_{6}^{-2}$ & K3 \# 75   \\  \hline
50 & 157 & 288 & $A_{4,4}$ & 19 & -288 & $2^16^124^1$& $2_{\rm I\!I}^{+2}8_{1}^{+1}3^{+2}$ & K3 \# 78  \\  \hline
51 & 162 & 336 & $\Z_2\times {\rm L}_2(7)$ & 20 & 196 & $14^2$ & $2_{\rm I\!I}^{+2}7^{+2}$ & max \# 5   \\  \hline
52 & 163 & 360 & $A_6$ & 19 & -180 & $3^160^1$& $4_{5}^{-1}3^{+2}5^{+1}$ & K3 \# 79   \\  
53 & 164 & 360 & $(\Z_3\times A_5){:}\Z_2$ & 20 & 225 & $15^2$ & $3^{-2}5^{-2}$ & max \# 10   \\  \hline
54 & 168 & 384 & $F_{384}$ & 19 & -256 & $4^18^2$ & $4_{3}^{+1}8_{2}^{+2}$ & K3 \# 80    \\  \hline
55 & 170 & 486 & $3^{1+4}{:}2$ & 18 & 243 & $3^5$ & $3^{+5}$ &  ${\cal S}$-lattice    \\  \hline
56 & 173 & 576 & $\Z_2^4{:}(S_3\times S_3)$ & 20 & 288 & $12^124^1$ & $4_{7}^{+1}8_{1}^{+1}3^{+2}$ & max \# 12   \\  \hline
57 & 177 & 660 & ${\rm L}_2(11)$ & 20 & 121 & $11^2$ & $11^{+2}$ & max \# 1   \\  \hline
58 & 178 & 720 & $ M_{10}$ & 20 & 120 & $2^160^1$ & $2_{5}^{+1}4_{1}^{+1}3^{-1}5^{+1}$ & max \# 9   \\
59 & 179 & 720 & $ S_6$ & 20 & 180 & $6^130^1$ & $2_{\rm I\!I}^{-2}3^{+2}5^{+1}$ & max  \# 8   \\  \hline
60 & 183 & 960 & $ M_{20} $ & 19 & -160 & $2^240^1$ & $2_{\rm I\!I}^{-2}8_{1}^{+1}5^{-1}$ & K3 \# 81 \\  \hline
61 & 185 & 972 & $3^{1+4}{:}2.2$ & 19 & -162 & $3^36^1$ & $2_{1}^{+1}3^{-4}$ &  ${\cal S}$-lattice    \\  \hline
62 & 187 & 1152 & $Q(\Z_3^2{:}\Z_2)$ & 20 & 192 & $8^124^1$ & $8_{6}^{-2}3^{-1}$ & max \# 11    \\  \hline
63 & 188 & 1344 & $\Z_2^3{:}{\rm L}_2(7)$ & 20 & 112 & $4^128^1$ & $4_{2}^{+2}7^{+1}$ & max \# 4   \\  \hline
64 & 190 & 1920 & $\Z_2^4{:}S_5$ & 20 & 160 & $4^140^1$& $4_{3}^{-1}8_{1}^{+1}5^{-1}$ & max \# 7   \\  \hline
65 & 191 & 1944 & $3^{1+4}{:}2.2^2$ & 20 & 108 & $3^16^2$& $2_{2}^{+2}3^{+3}$ &  max ${\cal S}$-lattice  \\  \hline
66 & 193 & 2520 & $A_7$ & 20 & 105 & $105^1$ & $3^{+1}5^{+1}7^{+1}$ & max \# 3    \\  \hline
67 & 196 & 5760 & $ \Z_2^4{:}A_6$ & 20 & 96 & $4^124^1$ & $4_{5}^{-1}8_{1}^{+1}3^{+1}$ &  max \# 6    \\  \hline
68 & 197 & 20160 & ${\rm L}_3(4)$ & 20 & 84 & $2^142^1$ & $2_{\rm I\!I}^{-2}3^{-1}7^{-1}$ &  max \# 2    \\  \hline
69 & 198 & 29160 & $3^4{:}A_6$ & 20 & 81 & $3^29^1$ & $3^{+2}9^{+1}$ & max ${\cal S}$-lattice   \\  \hline
\end{longtable}

\begin{table}
\caption{The $34$ coinvariant lattices $L_G$ with several groups}\label{Lgroups}

\bigskip
$
\begin{array}{rrl}
\mbox{No.} & \#\{G\} & \mbox{$G$-No.} \\ \hline

69 & 40 & \{ 198, 195, 194, 192, 189, 184, 181, 180, 176, 175, 171, 169, 161, 160, 159, 155, 153, 152, 139, \\
   &    &  137, 129, 128, 126, 118, 117, 116, 115, 113, 98, 95, 82, 80, 65, 64, 63, 62, 52, 51, 21, 20 \} \\

67 &  7 & \{ 196, 182, 172, 165, 149, 145, 91 \} \\

65 &  7 & \{ 191, 186, 151, 131, 127, 59, 28 \} \\

64 &  2 & \{ 190, 158 \} \\

63 &  3 & \{ 188, 141, 100 \} \\

62 & 10 & \{ 187, 174, 167, 166, 147, 146, 144, 123, 94, 89 \} \\

61 &  3 & \{ 185, 130, 26 \} \\

60 &  3 & \{ 183, 136, 114 \} \\

57 &  3 & \{ 177, 99, 23 \} \\

56 &  2 & \{ 173, 156 \} \\

55 & 12 & \{ 170, 154, 138, 119, 97, 96, 66, 61, 49, 18, 10, 4 \} \\

54 & 11 & \{ 168, 133, 124, 106, 105, 104, 75, 71, 69, 45, 33 \} \\

53 &  5 & \{ 164, 142, 101, 67, 32 \} \\

51 &  3 & \{ 162, 83, 31 \} \\

50 &  5 & \{ 157, 134, 121, 92, 84 \} \\

49 &  6 & \{ 150, 103, 87, 72, 39, 38 \} \\

48 &  2 & \{ 148, 125 \} \\

47 &  3 & \{ 143, 122, 57 \} \\

45 &  2 & \{ 135, 109 \} \\

42 &  5 & \{ 112, 78, 60, 30, 29 \} \\

39 &  2 & \{ 108, 76 \} \\

38 &  2 & \{ 107, 73 \} \\

34 &  2 & \{ 88, 58 \} \\

33 &  2 & \{ 86, 55 \} \\

31 &  2 & \{ 81, 48 \} \\

29 &  2 & \{ 77, 41 \} \\

28 &  2 & \{ 74, 43 \} \\

26 &  3 & \{ 68, 46, 15 \} \\

24 &  2 & \{ 54, 11 \} \\

22 &  2 & \{ 50, 19 \} \\

21 &  2 & \{ 47, 14 \} \\

18 &  2 & \{ 40, 12 \} \\

11 &  2 & \{ 24, 8 \} \\

10 &  2 & \{ 22, 7 \} 
\end{array}
$
\end{table}

\newpage

\noindent{\bf \Large Postscript}

\smallskip

In our paper \cite{HM} we classified all orbits of fixed-point sublattices 
of the Leech lattice and their respective stabilizers inside ${\rm Co}_0$.\ 
There are $290$ different cases, listed in Table~1 of~\cite{HM} .\
This allows us to obtain the $22$ classes of maximal admissible groups $H$
as described in Theorem~6.1 (a), (b), (c) by selecting those fixed-point lattices $\Lambda^G$
from Table~1 (loc.\ cit) for which $H$ is a subgroup of the stabilizer $G$
with $\Lambda^H=\Lambda^G$ containing only admissible elements.\ 
The thirteen groups of part (a) with $\alpha(\Lambda^G)\geq 2$ correspond to the entries 
\#102, \#106, \#108, \#110, \#111, \#112, \#118, \#119, \#120,  \#121, \#128, \#129, \#134,
the two groups of part (b) with $\alpha(\Lambda^G)=1$ correspond to the entries \#101 and \#109
and the seven groups in part (c) with $\alpha(\Lambda^G)=0$ correspond to the entries 
\#83, \#126, \#27, \#40, \#41, \#116, \#124.\ 
For the groups of parts (a) and (b), the group $H$ is the full stabilizer of
$\Lambda^H$, whereas for the groups in part (c), $H$ is strictly smaller than the full
stabilizer $G$, cf.\ Theorem~7.1.

\medskip
This approach will not lead to much shorter calculations in the present paper, furthermore
the results of~\cite{HM} also use the list of 279,343 conjugacy classes
of non $2$-subgroups of $2^{12}{:}M_{24}$ determined in the present paper.


\end{document}